\documentclass[a4paper,11pt,twoside]{article}
\usepackage[TS1,T1]{fontenc}
\usepackage{enumitem}
\usepackage[toc]{appendix}
\usepackage{float}
\usepackage{textcomp}
\usepackage{amsmath,amsfonts,verbatim,afterpage,euscript,mathrsfs,amssymb,empheq}
\usepackage{amsthm}
\usepackage{dsfont}
\usepackage{hyperref}
\usepackage{pst-fill,pst-grad}
\usepackage{textcomp}
\usepackage{multicol}
\newtheorem{Theorem}{Theorem}[section]
\newtheorem{Lemma}{Lemma}[section]
\newtheorem{Proposition}{Proposition}[section]
\newtheorem{Corollary}{Corollary}[section]
\newtheorem{Remark}{Remark}[section]
\newtheorem{Definition}{Definition}[section]
\def \vu{\vec{u}}
\def \vv{\vec{v}}

\def \vw{\vec{\omega}}

\def \vf{\vec{f}}

\def \vphi{\vec{\phi}}

\def \rot{\grad \wedge}
\def \grad{\vec{\nabla}}
\def \vn{\vec{\nabla}}
\def \div{\operatorname{div}}
\def \supp{\operatorname{supp}}
\def \R{\mathbb{R}^3}
\def \M{\mathcal{M}}

\newcommand{\mysection}{\setcounter{equation}{0} \section}

\title{\bf Partial regularity and $L^3$-norm concentration effects around possible blow-up points for the micropolar fluid equations}
\usepackage{authblk}
\author{Diego Chamorro\footnote{\emph{diego.chamorro@univ-evry.fr} (corresponding author)} }
\author{David Llerena\footnote{\emph{david.llerena@univ-evry.fr}} }
\affil{\footnotesize LaMME, Univ. Evry, CNRS, Universit\'e Paris-Saclay, 91025, Evry, France.}
\begin{document}
\sloppy
\maketitle
%%%%%%%%%%%%%%%%%%%%%%%%%%%%%%%%%%%%%%%%%%%%%%%%%%%
\begin{scriptsize}
\abstract{The micropolar fluid system is a model based on the Navier-Stokes equations which considers two coupled variables: the velocity field $\vu$ and the microrotation field $\vw$. Assuming an additional condition over the variable $\vu$ we will first prove that weak solutions $(\vu, \vw)$ of this system are smooth. Then, we will present a concentration effect of the $L^3_x$ norm of the velocity field $\vu$ near a possible singular time.}\\[3mm]
{\bf \scriptsize Keywords: micropolar equations; regularity; blow-up.}\\
\textbf{\scriptsize Mathematics Subject Classification: 35Q30; 35D30; 35B65; 35B44.}
\end{scriptsize}
%\tableofcontents
%%%%%%%%%%%%%%%%%%%%%%%%%%%%%%%%%%%%%%%%%%%%%%%%%%%
\mysection{Introduction}
%%%%%%%%%%%%%%%%%%%%%%%%%%%%%%%%%%%%%%%%%%%%%%%%%%%
In this paper we are interested in studying some properties of weak solutions of the micropolar fluid equations. Recall that these equations are given by the following coupled system
\begin{empheq}[left=\empheqlbrace]{alignat = 3}
&\partial_t \vu = \Delta \vu -(\vu \cdot \vn)\vu-\vn p +\frac{1}{2}\rot\vw ,\qquad \div(\vu)=0,\label{MicropolarFluidsEquationsEqua1}\\[1mm]
&\partial_t \vw = \Delta \vw +\vn \div(\vw)-\vw -(\vu \cdot \vn)\vw+\frac{1}{2}\rot\vu ,\label{MicropolarFluidsEquationsEqua2}\\[1mm]
&\vu (0,x)=\vu_0(x), \;\;\vw(0,x)=\vw_0(x)\quad \mbox{and}\quad \div(\vu _0)=0,\quad x\in \R.\notag
\end{empheq}
In the previous equations the initial data are $\vu_0$ and $\vw_0$ and the variables are $(\vu, p, \vw)$ where the vector field $\vu:[0, +\infty[\times \R \longrightarrow \R$ is the velocity field of the fluid, the scalar function $p:[0, +\infty[\times \R \longrightarrow \mathbb{R}$ is the internal pressure and the vector field $\vw:[0, +\infty[\times \R \longrightarrow \R$ is the angular velocity or the microrotational velocity. It is worth to remark here that the first equation \eqref{MicropolarFluidsEquationsEqua1} above is related to the incompressible 3D Navier-Stokes system (we have $\div(\vu)=0$) while the second equation \eqref{MicropolarFluidsEquationsEqua2} gives the evolution of the microrotational velocity field $\vw$.\\

This system of PDEs was introduced in 1966 by Eringen in \cite{Eri66} and it has been  studied by many authors, see \emph{e.g.} \cite{SerrinMicropolar}, \cite{ChLl21},  \cite{Gal97}, \cite{LoMelo}, \cite{Yam05}, \cite{Yua08} and the references therein. Apart from the various applications of this model (see for example \cite{BegBahr08}, \cite{Gay13} and \cite{Mit02}), a very interesting feature from the mathematical perspective of this micropolar fluid system is the fact that the variable $\vw$ is not a divergence-free vector field, and this makes its study quite different from other systems of PDEs based on the Navier-Stokes equations (such as the magneto-hydrodynamic equations, see \emph{e.g.} \cite{ChHe21}).\\

Let us start with two simple remarks concerning the system (\ref{MicropolarFluidsEquationsEqua1})-(\ref{MicropolarFluidsEquationsEqua2}). First, it is easy to observe that the equation related to the variable $\vu$ in \eqref{MicropolarFluidsEquationsEqua1} is invariant according to the following scaling 
\begin{align*}
\vu_{\lambda}(t,x)=\lambda \vu(\lambda^2t,\lambda x),\quad p_{\lambda}(t,x)=\lambda^2 p(\lambda^2t, \lambda x)\quad \text{and}\quad \vw_{\lambda}=\lambda^2\vec { \omega}(\lambda^2t,\lambda x) \quad\mbox{ where } \lambda >0,
\end{align*}
however the triplet $(\vu_\lambda, p_\lambda, \vw_\lambda)$ is no longer a solution for the whole micropolar system since the second equation \eqref{MicropolarFluidsEquationsEqua2} does not have a ``natural'' scaling that preserves the structure of the equation (due to the presence of the term $\vw$),\label{Scaling} and his fact reveals one of the major differences between these two equations. \\

We continue by observing that the information about the pressure $p$ can be easily obtained from the variable $\vu$: indeed, by formally applying the divergence operator in the equation \eqref{MicropolarFluidsEquationsEqua1}, since $\div(\vu)=0$ and $\div(\rot \vw)=0$, we obtain the following equation for the pressure:
\begin{equation}\label{Eq_Pression}
-\Delta p=\div((\vu\cdot \vn)\vu),
\end{equation}
so we can write $p=\frac{1}{(-\Delta)}\div((\vu\cdot \vn)\vu)$ and then pressure $p$ is  only related to the velocity field $\vu$, therefore we will consider the pair $(\vu, \vw)$ as the main variables. These two simple remarks will be essential in the sequel.\\

Note now that Leray-type weak solutions of the previous system \eqref{MicropolarFluidsEquationsEqua1}-\eqref{MicropolarFluidsEquationsEqua2} can be easily obtained: indeed, from two $L^2$ initial data $\vu_0, \vw_0$ and by a classical mollification argument we can construct global solutions $\vu, \vw\in L^\infty_tL^2_x\cap L^2_t\dot{H}^1_x$ that satisfy some energy inequalities:
%%%%%%%%%%%%%%%%%%%%%%%%%%%%%%%%%%%%%%%%%%%%%%%%%%%
\begin{Definition}[{\bf Leray-type weak solutions}]
Let $\vu_0, \vw_0\in L^2(\R)$ with $\div(\vu_0)=0$. We will say that $(\vu,p,\vw)$ is a \emph{Leray-type weak solution} of the micropolar fluid equations  \eqref{MicropolarFluidsEquationsEqua1} and \eqref{MicropolarFluidsEquationsEqua2} with initial value  $\vu_0$ and $\vw_0$ if 
${\vu\in L^\infty(]0,+\infty[,L^2(\R))\cap L^2(]0,+\infty[,\dot H^1(\R) )}$, ${\vw\in L^\infty(]0,+\infty[,L^2(\R))\cap L^2(]0,+\infty[, H^1(\R) )}$ and if for every $t\in]0,+\infty[$ we have
the following energy inequality 
$$\|\vu(t,\cdot)\|_{L^2}^2+\|\vw(t,\cdot)\|_{L^2}^2+\int_0^t\|\vu(s,\cdot)\|_{\dot H^1}^2+2\|\vw(s,\cdot)\|_{\dot{H}^1}^2+\|\vw(s,\cdot)\|_{L^2}^2+2\|\div(\vw)(s,\cdot)\|_{L^2}^2ds\le\|\vu_{0}\|_{L^2}+\|\vw_{0}\|_{L^2}^2.$$
\end{Definition}
%%%%%%%%%%%%%%%%%%%%%%%%%%%%%%%%%%%%%%%%%%%%%%%%%%%
\noindent Leray-type weak solutions will constitute the main framework of this work, however, just as for the Navier-Stokes equations, the complete study of the properties of these solutions remains a challenging open problem for the micropolar fluid equations.\\

In this article, we first want to perform a separate study for each variable $\vu$ and $\vw$ in order to obtain, by considering a hypothesis on the single variable $\vu$, some regularity for the couple $(\vu, \vw)$. Then we will deduce a concentration phenomenon for the $L^3_x$ norm of $\vu$ when approaching a potential blow-up time. Let us stress here that we will avoid as much as possible any additional information over $\vw$ (except for the $L^\infty_tL^2_x\cap L^2_t\dot{H}^1_x$ framework). To do so, we need now to introduce some definitions that underline this separation of the information between the variables $\vu$ and $\vw$:
%%%%%%%%%%%%%%%%%%%%%%%%%%%%%%%%%%%%%%%%%%%%%%%%%%%
\begin{Definition}[{\bf Partial suitable solution}]\label{Def_PartialSuitable}
We will say that the triplet $\vu,\vw:]0, T[\times \R\longrightarrow \R$ and $p: ]0, T[\times\R\longrightarrow \mathbb{R}$ is a \emph{partial suitable solution} of the micropolar fluid equations \eqref{MicropolarFluidsEquationsEqua1} and \eqref{MicropolarFluidsEquationsEqua2} over a regular open set $\Omega\subset]0, T[\times\R$ with $0<T<+\infty$, if:
\begin{itemize}
\item[1)] we have $\vu,\vw \in L_t^{\infty}L_x^2(\Omega)\cap L_t^2\dot H_x^1(\Omega)$, $p\in L_{t,x}^{\frac{3}{2}}(\Omega)$ and the variables $(\vu,p,\vw)$ satisfy in the weak sense the equations \eqref{MicropolarFluidsEquationsEqua1} and \eqref{MicropolarFluidsEquationsEqua2} over $\Omega$, 
\item[2)] for all $\phi\in \mathcal{D}(\Omega)$ the following local energy inequality is satisfied
\end{itemize}
\begin{align}
\int_{\R} |\vu|^2\phi(t,\cdot) dx&+ 2 \int_{s<t} \int _{\R}|\grad \otimes \vu|^2\phi dx ds\le \int_{s<t}\int_{\R}(\partial_t\phi + \Delta \phi)|\vu|^2dyds+2\int_{s<t}\int_{\R}p (\vu\cdot \grad \phi ) dyds\nonumber\\
&+\int_{s<t}\int _{\R} |\vu|^2 (\vu \cdot \grad)\phi dx ds+\int_{s<t}\int_{\R}(\rot \vw)\cdot (\phi\vu)dyds.\label{LocalEnergyInequality}
\end{align}
\end{Definition}
%%%%%%%%%%%%%%%%%%%%%%%%%%%%%%%%%%%%%%%%%%%%%%%%%%%
\noindent Let us observe that the previous inequality is only related to the structure of the first equation \eqref{MicropolarFluidsEquationsEqua1} and it is not related to the evolution of $\vw$ given in \eqref{MicropolarFluidsEquationsEqua2}. This notion of \textit{partial suitable solution} was introduced in our previous work \cite{ChLl23}, where we studied the interdependence of the variables in the $\varepsilon $-regularity theory (based on the celebrated work of Caffarelli, Kohn and Nirenberg \cite{CKN} for the Navier-Stokes system). \\

\noindent Now we introduce the following definition of partial regular points.
%%%%%%%%%%%%%%%%%%%%%%%%%%%%%%%%%%%%%%%%%%%%%%%%%%%
\begin{Definition}[{\bf Partial regular point/Partial singular point}]\label{Def_PartialRegularPoints}
A point $(t_0,x_0)\in \Omega \subset ]0,T[\times \R$ is a \emph{partial regular point} for the micropolar fluid equations \eqref{MicropolarFluidsEquationsEqua1} and \eqref{MicropolarFluidsEquationsEqua2} if there exists $r>0$ small enough such that $]t_0-r^2,t_0[\times B_{x_0,r}\subset \Omega$ and such that $\vu\in L^\infty_{t,x}(]t_0-r^2,t_0[\times B_{x_0,r})$. On the other hand, we will say that a point $(t_0,x_0)$ is \emph{partially singular} if it is not \emph{partially regular}.
\end{Definition}
%%%%%%%%%%%%%%%%%%%%%%%%%%%%%%%%%%%%%%%%%%%%%%%%%%%
\noindent In the two previous definitions we do not impose any constraint in the variable $\vw$. However, as we shall see, it will be enough to impose some conditions to the velocity field $\vu$  to obtain a gain of information (regularity or integrability) for both variables $\vu$ and $\vw$. Concerning this last notion of partial regular points, the regularity of the variables $\vu$ and $\vw$ will be obtained from the local hypothesis $\vu\in L^\infty_tL^\infty_x$. Although this is a rather ``reasonable'' result, to the best of our knowledge it was not studied in detail before, so we give a proof in the Theorem \ref{Theo_SerrinMP} below.\label{ReferenceLocaleSerrin}\\

As we aim to study the behavior of the variable $\vu$ around some potential blow-up point, we need to establish some very specific regularity results that were not treated before. Thus, our first result explore a gain of integrability when assuming a local $L^\infty_tL^3_x$ hypothesis for the velocity field $\vu$:
%%%%%%%%%%%%%%%%%%%%%%%%%%%%%%%%%%%%%%%%%%%%%%%%%%%
\begin{Theorem}[{\bf Partial interior regularity}]\label{Theo_LocalPoints}
Consider $(\vu,p,\vw)$ a partial suitable solution over a regular set $\Omega\subset ]0,T[\times \R$ with $0<T<+\infty$ of the micropolar fluid equations \eqref{MicropolarFluidsEquationsEqua1} and \eqref{MicropolarFluidsEquationsEqua2} in the sense of the Definition \ref{Def_PartialSuitable} above.  Assume that for some point $(t_0,x_0)\in \Omega$ there exists $R>0$ such that we have $]t_0-R^2,t_0[\times B_{x_0,R}\subset \Omega$ and such that we have the information $\vu\in L^\infty(]t_0-R^2,t_0[,L^3(B_{x_0,R}))$. Then there exists $r>0$ with $0<r\le \frac{R}{2}$ such that ${\vu\in L^\infty_{t,x}(]t_0-r^2,t_0[\times B_{x_0,r}})$, \emph{i.e.} the point $(t_0, x_0)$ is partially regular in the sense of the Definition \ref{Def_PartialRegularPoints} above. 
\end{Theorem}
%%%%%%%%%%%%%%%%%%%%%%%%%%%%%%%%%%%%%%%%%%%%%%%%%%%
\noindent Some remarks are in order here. First note again that we only impose some additional information on $\vu$ and not on the variable $\vw$ (which is consistent with the general spirit of this article), however the conclusion applies only to $\vu$. Remark next that this additional control, namely the fact that $\vu\in L^\infty_tL_x^3$ (locally), is reminiscent of the endpoints of the Serrin criterion for the classical Navier-Stokes system where it is traditional to assume locally $\vu\in L^p_tL_x^q$ with $\frac{2}{p}+\frac{3}{q}< 1$ (see \cite{Serr62}, \cite{Serr63}). The case when $\frac{2}{p}+\frac{3}{q}=1$ with $q>3$ was obtained by \cite{Struwe88} and \cite{Takahashi90} while the endpoint $p=+\infty$ and $q=3$ (which is the case studied in the Theorem \ref{Theo_LocalPoints} above) was obtained for the Navier-Stokes equations in \cite{Eus03}. Note also that for the Navier-Stokes equations some of these results were generalized to the framework of parabolic Morrey spaces $\M^{p,q}_{t,x}$ in \cite{OLeary} (these spaces will constitute one of the main tools of this article, see the expressions (\ref{DefMorrey}) and (\ref{DefMorreyparabolico}) below for a precise definition of Morrey spaces). For the micropolar fluid equations see our recent works \cite{ChLl22} and \cite{ChLl23} where we assumed a local control of the velocity field in terms of the parabolic Morrey space $\M^{p,q}_{t,x}$ with $2<p\leq q$ and $5<q\leq 6$. Let us mention finally that the treatment of the endpoint $p=+\infty$ and $q=3$ as announced in Theorem \ref{Theo_LocalPoints} above seems to be new in the context of the micropolar fluid equations.\\

In our next result, assuming a global in space $L^3_x$ control, we will characterize the continuity in time information for the velocity field $\vu$ in terms of partial regular points. More precisely we have:
%%%%%%%%%%%%%%%%%%%%%%%%%%%%%%%%%%%%%%%%%%%%%%%%%%%
\begin{Theorem}\label{Theo_NormL3}
Let $(\vu,p,\vw)$ be a weak Leray-type solution over $]0,+\infty[\times \R$ of the micropolar system \eqref{MicropolarFluidsEquationsEqua1} and \eqref{MicropolarFluidsEquationsEqua2} with $\vu,\vw \in L^\infty_tL^2_x\cap L^2_t\dot{H}^1_x $ such that for some time $0<\delta<T<+\infty$ we have ${\vu\in L^\infty(]\delta,T[,L^3(\R))}$.  Then the velocity field $\vu$ satisfies $\vu\in \mathcal{C}(]\delta, T[,L^3(\R))$ if and only if each point $(t_0,x_0)\in ]\delta, T[$ is a partial regular point in the sense of Definition \ref{Def_PartialRegularPoints}. 
\end{Theorem}
%%%%%%%%%%%%%%%%%%%%%%%%%%%%%%%%%%%%%%%%%%%%%%%%%%%
\noindent One of the main differences between this result and the previous Theorem \ref{Theo_LocalPoints} lies in the fact that we no longer require here the partial suitability condition (\ref{LocalEnergyInequality}). Indeed, as we shall see later on, the global in space hypothesis $\vu\in L^\infty(]0,T[,L^3(\R))$ is strong enough to ensure an interesting global estimate. Again, the variable $\vw$ seems to play no particular role in the statement of the result, but must be studied very carefully in the computations.\\ 

To the best of our knowledge, Theorem \ref{Theo_LocalPoints} and Theorem \ref{Theo_NormL3} are new in the setting of the micropolar fluid equations \eqref{MicropolarFluidsEquationsEqua1} and \eqref{MicropolarFluidsEquationsEqua2}. These results, although interesting for their own sake, are however merely preliminary results: indeed, our first main theorem states a blow-up criterion for Leray-type weak solution of the micropolar fluid equations  \eqref{MicropolarFluidsEquationsEqua1} and \eqref{MicropolarFluidsEquationsEqua2}:
%%%%%%%%%%%%%%%%%%%%%%%%%%%%%%%%%%%%%%%%%%%%%%%%%%%
\begin{Theorem}[{\bf Blow-up}]\label{Theo_BlowUp1}
Let $(\vu,p,\vw)$ be a Leray-type weak solution of the micropolar fluid equations \eqref{MicropolarFluidsEquationsEqua1} and \eqref{MicropolarFluidsEquationsEqua2}. Let $0<\mathcal{T}\le +\infty$ be the maximal time so that we have the control $ \vu\in\mathcal{C}(]0, \mathcal{T}[,L^3(\R))$. If $\mathcal{T}<+\infty$, then 
\begin{equation*}
\sup_{0<t<\mathcal{T}}\|\vu(t,\cdot)\|_{L^3}=+\infty.
\end{equation*}
\end{Theorem}  
%%%%%%%%%%%%%%%%%%%%%%%%%%%%%%%%%%%%%%%%%%%%%%%%%%%
\noindent The proof of this theorem will heavily rely on the previous results stated above. With all these results at our disposal, we can now tackle our second main theorem which is related to a refinement of the blow-up criterion stated in Theorem \ref{Theo_BlowUp1} above: indeed, we want now to study the concentration of the $L^3_x$-norm of the velocity field $\vu$ on balls centered at a singular point $(\mathcal{T},0)$ whose radius shrinks to zero as $t$ tends to $\mathcal{T}$.
%%%%%%%%%%%%%%%%%%%%%%%%%%%%%%%%%%%%%%%%%%%%%%%%%%%
\begin{Theorem}[{\bf $L^3$ concentration effect}]\label{Theo_BlowUpConcentration}
Let $(\vu,p,\vw)$ be a Leray-type weak solution of the micropolar fluid equations \eqref{MicropolarFluidsEquationsEqua1} and \eqref{MicropolarFluidsEquationsEqua2}. Assume  that  $0<\mathcal{T}<+\infty$ is the maximal time such that we have ${\vu\in \mathcal{C}(]0,\mathcal{T}[,L^\infty(\R))}$. Assume  that the point $(\mathcal{T},0)$ is  a partial singular point in the sense of the Definition \ref{Def_PartialRegularPoints} and the time  $\mathcal{T}$ satisfies the following condition: for some $r_0>0$ such that ${0<\mathcal{T}-r_0^2},$ we have
\begin{equation}\label{Teo_hypLocalSmoothing}
\sup_{x_0\in \R}\sup _{r\in ]0,r_0]}\sup _{t\in]\mathcal{T}-r^2,\mathcal{T}]}\frac{1}{r}\int_{B_{x_0,r}}|\vu(t,x)|^2dx =\mathfrak{M} <+\infty.
\end{equation}
Then, there exists $\varepsilon>0$, $\mathfrak S=\mathfrak S(\mathfrak{M})>0$ and $0<\delta<\mathcal{T}$ such that for all $t\in ]\mathcal{T}-\delta,\mathcal{T}[$, we have
\begin{equation}\label{concentration}
\int_{B_{0,\sqrt{\frac{\mathcal{T}-t}{\mathfrak S}}}}|\vu(t,x)|^3dx\ge \varepsilon.
\end{equation}
\end{Theorem}
%%%%%%%%%%%%%%%%%%%%%%%%%%%%%%%%%%%%%%%%%%%%%%%%%%%
\noindent Of course, with the estimate (\ref{concentration}) above it is quite straightforward to observe the announced concentration phenomenon of the $L^3_x$ norm for the velocity field $\vu$ when  $t$ tends to the ``blow-up'' time $\mathcal{T}$. Let us remark now that the constraint given in the expression (\ref{Teo_hypLocalSmoothing}) is known in the literature of the Navier-Stokes equations as the \emph{type I condition} (see \cite{BarPran19}, \cite{BarPran20}, \cite{KanMiuTsai19} and the references therein) and it can be interpreted in terms of Morrey spaces. Indeed, if $\vu$ satisfies the condition (\ref{Teo_hypLocalSmoothing}), then we have $\vu\in L^\infty_t\M_x^{2,3}\subset \M_{t,x}^{2,5}$. As it might be expected, the fact that $\vu, \vw\in \M_{t,x}^{p,q}$ with $p=2$ and $q=5$ falls outside the scope of the Serrin regularity criterion stated in terms of Morrey spaces where we need to impose that $2<p\leq q$ and $5<q\leq 6$ (see \cite{ChLl22}). This suggest that the values $p=2$ and $q=5$ may constitute a threshold: above these values the additional parabolic Morrey information will provide enough ``integrability'' to deduce a gain of regularity, while at $p=2$ and $q=5$ (or below) the parabolic Morrey control will not produce a consequent gain of information.\\

We also note that, although it is not difficult to exhibit a ``domination'' of the variable $\vu$ over the variable $\vw$ when considering regularity results (in the sense that it is enough to impose some conditions on $\vu$ to obtain a gain for \emph{both} variables $\vu$ and $\vw$), the techniques developed in this article do not seem to provide any information about the behavior of $\vw$ close to a potential blow-up point. However, we can possibly conjecture that a blow-up for the variable $\vw$ will impact the behavior of the velocity field $\vu$, but the complete study of this problem would probably require some additional work which is beyond the scope of this article.\\

The plan of the paper is as follows: Theorem \ref{Theo_LocalPoints} will be studied in Section \ref{Sec_LocalPoints} and in Section \ref{Sec_RegularityNormL3} we prove Theorem \ref{Theo_NormL3}. Section \ref{Secc_BlowUp1} is devoted to the proof of Theorem \ref{Theo_BlowUp1} while the $L^3$-norm concentration effect stated in the Theorem \ref{Theo_BlowUpConcentration} is treated in the Seccion \ref{Sec_Concentration}. In the appendix  \ref{Sec_LocalregularityMP} and the appendix \ref{Sec_PartialRegularity} we present some regularity results related to the system \eqref{MicropolarFluidsEquationsEqua1}-\eqref{MicropolarFluidsEquationsEqua2} that were not explicitly proven before and that are needed here to perform some computations.
%%%%%%%%%%%%%%%%%%%%%%%%%%%%%%%%%%%%%%%%%%%%%%%%%%%
\subsection*{Notations}
Throughout this paper we fix the following notation for two different types of parabolic balls centered in a point $(t_0,x_0)\in ]0,+\infty[\times \R$: we define the sets ${\bf Q}_r(t_0,x_0)$ and $ Q_r(t_0,x_0)$ by
\begin{eqnarray}
{\bf Q}_r(t_0,x_0)&=&]t_0-r^2,t_0+r^2[\times B_{x_0,r},\label{Def_Bolas1}\\[3mm]
\text{and} \quad Q_r(t_0,x_0)&=&]t_0-r^2,t_0[\times B_{x_0,r},\label{Def_Bolas2}
\end{eqnarray}
for some $0<r^2<t_0$ and $B_{x_0,r}=B(x_0,r)$. When the context is clear we will write ${\bf Q}_r$ (or $Q_r$) instead of ${\bf Q}_r(t_0, x_0)$ (or $Q_r(t_0, x_0)$). Note that we clearly have $Q_r(t_0, x_0)\subset {\bf Q}_r(t_0, x_0)$.\\

\noindent Morrey spaces $\mathcal{M}_{x}^{p,q}(\R)$ with $1<p\leq q<+\infty$ are defined as the set
$$\mathcal{M}_{x}^{p,q}(\R)=\{\vf:\R\longrightarrow\R: \vf\in L^p_{loc}(\R), \|\vf\|_{\mathcal{M}_{x}^{p,q}}<+\infty\},$$
where
\begin{equation}\label{DefMorrey}
\|\vf\|_{\mathcal{M}_{x}^{p,q}}=\underset{x_{0}\in \R, r>0}{\sup}\left(\frac{1}{r^{3(1-\frac{p}{q})}}\int_{B_{x_0,r}}|\vf(x)|^{p}dx\right)^{\frac{1}{p}}.
\end{equation}
For $1< p\leq q<+\infty$, the \emph{parabolic} Morrey spaces $\mathcal{M}_{t,x}^{p,q}(\mathbb{R}\times \R)$ are defined as the set of measurable functions $\vf:\mathbb{R}\times\R\longrightarrow \R$ that belong to the space $(L^p_{t,x})_{loc}$ such that $\|\vf\|_{\M_{t,x}^{p,q}}<+\infty$ where
\begin{equation}\label{DefMorreyparabolico}
\|\vf\|_{\mathcal{M}_{t,x}^{p,q}}=\underset{x_{0}\in \R, t_{0}\in \mathbb{R}, r>0}{\sup}\left(\frac{1}{r^{5(1-\frac{p}{q})}}\int_{|t-t_{0}|<r^{2}}\int_{B_{x_0,r}}|\vf(t,x)|^{p}dxdt\right)^{\frac{1}{p}}.
\end{equation}
Although not explicitly present in the statement of our results, Morrey spaces will play a crucial role in our computations. Indeed, these functional spaces are a very useful tool when addressing problems related harmonic analysis or to the regularity of a large class of PDEs, see \emph{e.g.} \cite{Admas}, \cite{ChLl23}, \cite{Kukavica}, \cite{PGLR1}, \cite{Robinson} and the references therein for some interesting applications of these spaces.
%%%%%%%%%%%%%%%%%%%%%%%%%%%%%%%%%%%%%%%%%%%%%%%%%%%
%%%%%%%%%%%%%%%%%%%%%%%%%%%%%%%%%%%%%%%%%%%%%%%%%%%
\mysection{Proof of Theorem \ref{Theo_LocalPoints}}\label{Sec_LocalPoints}
Recall that we plan to prove that $\vu\in L^\infty_{t,x}(Q_{r}(t_0,x_0))$ for some $0<r\leq \frac{R}{2}$. For this purpose we begin by introducing some useful preliminary results and important properties satisfied by any partial suitable solution $(\vu,p,\vw)$ of the micropolar fluid equations \eqref{MicropolarFluidsEquationsEqua1} and \eqref{MicropolarFluidsEquationsEqua2} such that $\vu\in L_t^\infty L_x^3(Q_R(x_0,t_0))$ where $Q_R(x_0,t_0)=]t_0-R^2,t_0[\times B_{x_0,R}$. 
\begin{itemize} [leftmargin=*]
\item [$\bullet$] First, under the hypotheses given in Theorem \ref{Theo_LocalPoints} over $(\vu,p,\vw)$,  we can obtain the following information 
\begin{equation}\label{ContinuityInTimePartialSuitableSolution}
\vu\in\mathcal{C}\left(\left[t_0-\tfrac{R^2}{4},t_0\right],L^{\frac{5}{4}}(B_{x_0,\frac{R}{2}})\right).
\end{equation}
For proving this result, we will need the following lemma given in \cite[Lemma 9.6, pg 177]{Tsai19}. 
\begin{Lemma}\label{Lem_CoerciveEstimate}
Let $1<s<q<+\infty$. If $(\vv,p)$ is a weak solution of the time-dependent Stokes system, 
\begin{equation*}
\partial \vv -\Delta \vv +\grad p =\vf, \quad \div(\vv)=0,
\end{equation*}
such that for $R>0$, $\vv\in L_t^sL_x^1(Q_R)$, $p\in L^s_tL^1_x(Q_R)$ with an external force $\vf \in L^s_tL^q_x(Q_R) $. Then, for all $0<r<R$, we have
$$\|\partial_t \vv\|_{L^s_tL^q_x(Q_r)}+\|\Delta \vv\|_{L^s_tL^q_x(Q_r)}+\|\grad p\|_{L^s_tL^q_x(Q_r)}\le C(\|\vf\|_{L^s_tL^q_x(Q_R)}+\|\vv\|_{L^s_tL^1_t(Q_R)}+\|p\|_{L^s_tL^1_t(Q_R)}).$$
\end{Lemma}
This lemma is known in the literature as the \emph{coercive estimates for the Stokes system}, for further details about these estimates we refer to \cite[Theorem 5.4]{Robinson} or \cite[Proposition 6.7]{Seregin14}.\\

Now, we will see how to deduce \eqref{ContinuityInTimePartialSuitableSolution} by using the aforementioned lemma. Notice that by the Hölder inequality with $\frac{4}{5}=\frac{3}{10}+\frac{1}{2}$, and since $\vu\in  L^\infty_tL^2_x(Q_R)\cap L^2_t \dot H^1_x(Q_R)$  by hypothesis,  we obtain
\begin{equation*}
\|(\vu\cdot \grad)\vu\|_{L^{\frac{5}{4}}_{t}L^{\frac{5}{4}}_{x}(Q_R)}\le \|\vu\|_{L^{\frac{10}{3}}_{t}L^{\frac{10}{3}(Q_R)}_{x}}\|\grad\otimes\vu\|_{L^2_{t}L^2_x}\le\|\vu\|_{L^{2}_{t}L^6_x(Q_R)}^\frac{3}{5} \|\vu\|_{L^{\infty}_{t}L^2_x(Q_R)}^\frac{2}{5}\|\grad\otimes\vu\|_{L^2_{t}L^2_x(Q_R)}<+\infty.
\end{equation*}
Furthermore, since $\vw \in L^2_t\dot H^1_x(Q_R)$, $\vu\in L^\infty_tL^3_x(Q_R)$ and $p\in L^\frac{3}{2}_{t}L^\frac{3}{2}_{x}(Q_R)$ by hypotheses and since $Q_R$ is a bounded set, we obtain that $\rot \vw\in L^{2}_{t}L^2_x(Q_{R})\subset L^{\frac{5}{4}}_{t}L^{\frac{5}{4}}_{x}(Q_{R})$, $\vu\in L^\infty_tL^3_x(Q_R)\subset L^\frac54_{t} L^1_{x}(Q_R)$ and $p\in L^{\frac{3}{2}}_{t}L^{\frac{3}{2}}_{x}(Q_R)\subset L^\frac54_{t} L^1_{x}(Q_R)$. Thus, since $(\vu,p)$ satisfies the system
\begin{equation*}
\partial_t\vu-\Delta \vu +\grad p=(\vu\cdot \grad)\vu+\frac{1}{2}\rot \vw,\quad \div(\vu)=0,
\end{equation*} 
and we have deduced that $(\vu\cdot \grad )\vu, \rot \vw \in L^{\frac{5}{4}}_{t}L^{\frac{5}{4}}_{x}(Q_{R})$  and $\vu, p \in L^{\frac{5}{4}}_{t}L^{1}_{x}(Q_{R})$, then from the Lemma \ref{Lem_CoerciveEstimate} above, we obtain the following information over the time derivative of the velocity field: 
$$\partial_t \vu \in L_{t,x}^\frac{5}{4}(Q_\frac{R}{2}).$$ 
With this information at hand we obtain for almost all $t\in ]t_0-\frac{R^2}{4},t_0[$ that there exists a vector field $\vec U\in L^{\frac{5}{4}}(\R)$  such that we have the expression $\vu(t,\cdot)=\displaystyle{\int ^t_{t_0-\frac{R^2}{4}}\partial_t \vu (t,\cdot)dt+\vec U}$ from which we can deduce that $\vu\in\mathcal{C}\left(\left[t_0-\frac{R^2}{4},t_0\right],L^{\frac{5}{4}}(B_{x_0,\frac{R}{2}})\right)$ (see for instance \cite[Lemma 3.2]{Tsai19} and \cite[Corollary 1.4.36]{Caz98}).\\

It is worth noting that from \eqref{ContinuityInTimePartialSuitableSolution} we are able to study the behavior of the solution in the closed interval $[t_0-\frac{R^2}{4},t_0]$ even though some of the initial hypotheses are stated in the bigger (but open) interval $]t_0-R^2,t_0[$.\\

\item [$\bullet$] Secondly, observe that from the hypothesis $\vu\in L^\infty_tL^3_x(Q_R)$, we have that $\vu(t,\cdot)\in L^3(B_{x_0,R})$ for almost all $t\in ]t_0-R^2,t_0[$, however we will deduce, using \eqref{ContinuityInTimePartialSuitableSolution}, that for any $t\in [t_0-\frac{R^2}{4},t_0]$, we have ${\vu(t,\cdot)\in L^3(B_{x_0,\frac{R}{2}})}$ (and not only for almost all $t\in ]t_0-\frac{R^2}{4},t_0[$). Indeed, let $t\in [t_0-\tfrac{R^2}{4},t_0]$ and  $(t_k)_{k\in \mathbb{N}}$ be a sequence in $]t_0-\frac{R^2}{2},t_0[$ such that $t_k \underset{k\to +\infty}{\longrightarrow} t$. Since ${\|\vu(t_k,\cdot)\|_{L^3(B_{x_0,\frac{R}{2}})}\le \|\vu\|_{L^\infty_tL_x^3(Q_R)}}$, using the Banach-Alaoglu theorem, there exists a subsequence $(t_{k_j})_{j\in\mathbb{N}}$ such that $(\vu(t_{k_j},\cdot))_{j\in \mathbb{N}}$ converges weakly-$\ast$ to some $\vv(t,\cdot)$  in $L^3(B_{x_0,\frac{R}{2}})$.
On the other hand, by the continuity in the time variable given in   $\eqref{ContinuityInTimePartialSuitableSolution}$, we have $\vu(t_{k_j},\cdot)\underset{j\to +\infty}{\longrightarrow}\vu(t,\cdot)$ strongly in $L^\frac{5}{4}(B_{x_0,\frac{R}{2}})$. Hence by uniqueness of the limit, one has
${\vu(t,\cdot)=\vec v(t,\cdot)\in L^3(B_{x_0,\frac{R}{2}})\cap L^\frac{5}{4}(B_{x_0,\frac{R}{2}})}$, and then we have proved that
\begin{equation}\label{EveryTimeInL3}
\text{ for any}\quad t\in [t_0-\frac{R^2}{4},t_0], \quad \text{we have}\quad  \vu(t,\cdot)\in L^3(B_{x_0,\frac{R}{2}}).
\end{equation}

Similar to the previous point, we remark that we are able to deduce some information on the behavior  of $\vu$ in the closed interval $[t_0-\frac{R^2}{4},t_0]$.\\

\item [$\bullet$] We give now some remarks about the pressure. Notice that we can decompose the pressure $p$ into two parts
\begin{equation}\label{PresureDecomposition}
p=\mathfrak{p}+\Pi,
\end{equation}
where  ${\mathfrak{p}=\frac{1}{(-\Delta)}(\div( \div (\phi \vu\otimes \vu))}$  with  $\phi$ a positive test function supported in $B_{x_0,\rho}$ such that $\phi=1$ in $B_{x_0,\frac{\rho}{2}},$ for  $0<\rho\le R$, and $\Pi$ is defined by $\Pi=p-\mathfrak{p}$. From the definition of $\mathfrak{p}$ we have,
\begin{equation*}
\|\mathfrak{p}\|_{L^\infty_tL^\frac{3}{2}_{x}(Q_\rho)}\le\|\mathfrak{p}\|_{L^\infty(]t_0-\rho^2,t_0[,L^\frac{3}{2}(\R))}=\left\|\frac{1}{(-\Delta)}(\div( \div (\phi \vu\otimes \vu))\right\|_{L^\infty(]t_0-\rho^2,t_0[,L^\frac{3}{2}(\R))}.
\end{equation*}  
Using the fact that the Riesz transforms are bounded in $L^\frac{3}{2}(\R)$ and $\supp(\phi)\subset B_{x_0,\rho} $, we can write
\begin{equation*}
\|\mathfrak{p}\|_{L^\infty_tL^\frac{3}{2}_{x}(Q_\rho)}\le C \|\phi \vu\otimes \vu\|_{L^\infty(]t_0-\rho^2,t_0[,L^\frac{3}{2}(\R))}\le  C \|\phi\|_{L^\infty(\R)}\| \vu\otimes \vu\|_{L^\infty(]t_0-\rho^2,t_0[,L^\frac{3}{2}(B_{x_0,\rho}))}.
\end{equation*}
Thus, since $\vu\in L^\infty_tL^3_x(Q_R)$ by hypothesis, we obtain
\begin{equation}\label{PressureP1Estimate}
\|\mathfrak{p}\|_{L^\infty_tL^\frac{3}{2}_{x}(Q_\rho)}
\le C\| \vu\|_{L^\infty_tL^3_{x}(Q_\rho(t_0,x_0))}^2\le C\| \vu\|_{L^\infty_tL^3_{x}(Q_R(t_0,x_0))}^2<+\infty.
\end{equation}
Now, since $\Pi=p-\mathfrak{p}$ and $p$ satisfies the equation \eqref{Eq_Pression}, we have for all $t\in ]t_0-\frac{R^2}{4},t_0[$,
\begin{equation*}
\Delta \Pi(t,\cdot)=\Delta p(t,\cdot)-\Delta \mathfrak{p}(t,\cdot)
=-\div(\div(\vu\otimes\vu))-\Delta\frac{1}{(-\Delta)}(\div(\div(\phi \vu\otimes \vu))).
\end{equation*} 
Thus, since $\phi\equiv1 $ in $B_{x_0,\frac{\rho}{2}}$, we observe that for all $t\in ]t_0-\frac{\rho^2}{4},t_0[$ we have $\Delta \Pi(t,\cdot)=0$ over $B_{x_0,\frac{\rho}{2}}$. Now, by the local estimates for harmonic functions (see \cite[Theorem 7]{Evans}) we have for any $0<\rho\le R$ the estimate
$\|\Pi(t,\cdot)\|_{L^\infty(B_{x_0,\frac{\rho}{2}})}\le C\|\Pi(t,\cdot)\|_{L^\frac{3}{2}( B_{x_0,\rho})}$. Moreover, since  $\Pi=p-\mathfrak{p}$, we have
\begin{align}
\|\Pi(t,\cdot)\|_{L^\infty(B_{x_0,\frac{\rho}{2}})}&\le C\|\mathfrak{p}(t,\cdot)\|_{L^\frac{3}{2}(B_{x_0,\rho})}+C\|p(t,\cdot)\|_{L^\frac{3}{2}(B_{x_0,\rho})}. \label{PressureHarmonicFunction}
\end{align}
Finally, as  $p\in L^{\frac{3}{2}}_{t}L^{\frac{3}{2}}_{x}(Q_\rho)$ by hypothesis and since $\mathfrak{p}\in L^\infty_tL_x^\frac{3}{2}(Q_\rho)\subset L^{\frac{3}{2}}_{t}L^{\frac{3}{2}}_{x}(Q_\rho)$ by \eqref{PressureP1Estimate}, by taking the $L^\frac{3}{2}$-norm in the time interval $]t_0-\frac{\rho^2}{4},t_0[$ in the expression above, we obtain 
\begin{equation}\label{ControlPi}
\|\Pi\|_{L^\frac{3}{2}_tL^\infty_x(Q_{\frac{\rho}{2}})}\le  C\|\mathfrak{p}\|_{L^\frac{3}{2}_{t}L^\frac{3}{2}_{x}(Q_{\rho})}+C\|p\|_{L^\frac{3}{2}_{t}L^\frac{3}{2}_{x}(Q_{\rho})}< +\infty.
\end{equation}
The decomposition (\ref{PresureDecomposition}) of the pressure as well as the controls (\ref{PressureP1Estimate})-(\ref{ControlPi}) will be useful in the sequel. \\[3mm]
\end{itemize}
These three points finishes the preliminary results. Now, we will prove that  ${\vu\in L^\infty_{t,x}(Q_{r})}$ for some ${0<r\le \frac{R}{2}}$. For this, we proceed by contradiction, assuming that for all $0<r\le\frac{R}{2}$, we have ${\vu\notin L^\infty(Q_r(t_0,x_0))}$ \textit{i.e.}, the point $(t_0,x_0)$ is partially singular in the sense of the Definition \ref{Def_PartialRegularPoints}. The strategy will consist in applying a scaling argument around the point $(t_0,x_0)$ and to study the behavior of some limit functions in order to exhibit a contradiction. Thus,  let us consider, for any $k\in \mathbb{N}$ the sequence $$\lambda_k=\sqrt{t_0-t_k},$$ where $(t_k)_{k\in \mathbb{N}}$ is a sequence such that for all $k\in \mathbb{N}$, $0<t_k<t_0$ and $t_k \underset{k\to +\infty}{\longrightarrow}  t_0$. Notice that $(\lambda_{k})_{k\in\mathbb{N}}\underset{k\to +\infty}{\longrightarrow}   0$ and it is a bounded sequence.

We extend now the functions $(\vu,p,\vw)$ by $0$ outside $Q_{\frac{R}{2}}(t_0,x_0)$ and we denote them by $(\underline{\vu},\underline{p},\underline{\vw})$. For any $k\in \mathbb{N}$, consider now $\vu_k$, $\vw_k$ and $p_k$ the following scaled functions: for any ${(s,y)\in[0,1]\times \R},$
\begin{align}\label{DefUkWk}
\begin{aligned}
\vu_{k}(s,y)&=\lambda_k \underline{\vu}(t_k+\lambda_k^2s,x_0+\lambda_k y),\qquad p_{k}(s,y)=\lambda_k^2 \underline{p}(t_k+\lambda_k^2s,x_0+\lambda_ky)\\[2mm]
& \text{and}\qquad \vw_{k}=\lambda_k^2\underline{\vw}(t_k+\lambda_k^2s,x_0+\lambda_ky).
\end{aligned}
\end{align}
\begin{Remark}\label{Rem_SupportScaled}
The support of the functions $(\vu_k,p_k,\vw_k)$ is included in $Q_{\frac{R}{2\lambda_{k}}}(1,0)$.
In the following we will consider $k$ large enough  such that $(t_0-t_k)< \frac{R^2}{4}$ and therefore $1<\frac{R}{2\lambda_{k}}$ (recall that $(\lambda_{k})_{k\in\mathbb{N}}$ converges to $0$). Hence the values of the functions $\vu_k$, $p_k$ and $\vw_k$ in $]0,1[\times \R$ correspond to the ones of $(\vu,p,\vw)$ in $]t_k,t_0[\times B_{x_0,\frac{R}{2}}$.
\end{Remark}
It is worth noting that $(\vu_k,p_k,\vw_k)$ is not a solution of the micropolar fluid systems \eqref{MicropolarFluidsEquationsEqua1} and \eqref{MicropolarFluidsEquationsEqua2} due to the lack of scaling of these two equations, as it was pointed out in the page \pageref{Scaling} of the introduction. Nevertheless the triplet $(\vu_k,p_k,\vw_k)$ satisfies the equation
\begin{equation}\label{ForceNS}
\partial_t\vu_k = \Delta \vu_k-(\vu_k\cdot \vn_k)\vu_k-\grad p_k +\frac{1}{2}\rot \vw_k,
\end{equation}
which can be seen as the classical Navier-Stokes equations with an external force $\rot\vw_k$ which is ``given'' and belongs to the space $L^2_tL^2_x$. Now, we want to prove the following convergences
\begin{align*}
\rot\vw_{k}\xrightarrow[ k\to +\infty ]{} 0,\quad  p_{k}\xrightarrow[ k\to +\infty ]{} p_{\infty}\quad \text{and}\quad \vu_{k}\xrightarrow[ k\to +\infty ]{}\vu_{\infty},
\end{align*}
in order to deduce that $(\vu_\infty,p_{\infty})$ is a solution of the Navier-Stokes equations in $]0,1[\times \R$,
\begin{equation*}
\partial_t \vu_\infty=\Delta \vu_\infty -\div(\vu_\infty\otimes\vu_\infty)-\grad p_{\infty},
\end{equation*}
and then a careful study of the properties of the solution $(\vu_\infty,p_{\infty})$ will leads us to the wished contradiction.\\

\begin{itemize}[leftmargin=*]
\item [$\bullet$] First, we study the convergence of the sequence $(\rot\vw_k)_{k\in \mathbb{N}}$ in the domain $]0,1[\times \R$. By the definition of $\vw_k$ given in \eqref{DefUkWk}, since  ${\supp(\rot \vw_k)\subset Q_{\frac{R}{2\lambda_{k}}}(1,0)}=]1-\frac{R^2}{4\lambda_{k}^2},1[\times B_{0,\frac{R}{2\lambda_{k}}}$ and $1-\frac{R^2}{4\lambda_{k}^2}<0$  by Remark \ref{Rem_SupportScaled}, we have
\begin{align*}
\|\rot\vw_{k}\|_{L^2(]0,1[,L^2(\R))}^2&=\int_{0}^1\int_{\R}|\rot\vw_k|^2dyds=\int_{0}^1\int_{B_{0,\frac{R}{2\lambda_{k}}}}|\lambda_k^3(\rot\underline{\vw})(t_k+\lambda_{k}^2s,x_0+\lambda_{k}y)|^2dyds.
\end{align*}
Now, by a change of variable and since $t_0-\frac{R^2}{4}<t_k$ by Remark \ref{Rem_SupportScaled}, we can write 
\begin{align*}
\|\rot\vw_{k}\|_{L^2(]0,1[,L^2(\R))}^2 &= \lambda_k\int_{t_k}^{t_0}\int_{B_{x_0,\frac{R}{2}}}|\rot\underline{\vw}|^2dyds\le \lambda_k\int_{t_0-\frac{R^2}{4}}^{t_0}\int_{B_{x_0,\frac{R}{2}}}|\rot\underline{\vw}|^2dyds.
\end{align*}
Using the fact that $\rot \underline{ \vw}=\rot \vw|_{Q_\frac{R}{2}(t_0,x_0)}$ by construction, we obtain
\begin{align}
\|\rot\vw_{k}\|_{L^2(]0,1[,L^2(\R))}^2\le \lambda_k\int_{Q_{\frac{R}{2}}(t_0,x_0)}|\rot\vw|^2dyds
&\le \lambda_k\int_{Q_{R}(t_0,x_0)}|\rot\vw|^2dyds.\label{EstimateRotOmega}
\end{align}
Since $\vw \in L^2_t\dot H_x^1(Q_R)$ by hypothesis, and $(\lambda_{k})_{k\in \mathbb{N}}$ converges towards zero when $k\to +\infty$, we have 
$$\rot\vw_k\underset{k\to + \infty}{\longrightarrow} 0 \quad\text{strongly in}\quad L^2(]0,1[, L^2(\R)).$$
\item [$\bullet$] Now, we study the convergence of $(p_k)_{k\in\mathbb{N}}$. Recall that for any $k\in \mathbb{N}$, we have  ${p_k(s,y)=\lambda_{k}\underline{p}(t_k+\lambda_{k}^2s,x_0+\lambda_{k}y)}$, where $\underline{p}=p|_{Q_{\frac{R}{2}}(t_0,x_0)}$. Since we can split the pressure  ${p}={\mathfrak{p}}+{\Pi}$ by \eqref{PresureDecomposition}, we can write for any $k\in\mathbb{N}$, $p_k=\mathfrak{p}_{\;k}+\Pi_{\;k}$, where
\begin{equation*}
\mathfrak{p}_{\;k}(s,y)=\lambda_k^2\underline{\mathfrak{p}}(t_k+\lambda_k^2 s,x_0+\lambda_k y),\quad \quad  \Pi_{\;k}(s,y)=\lambda_k^2\underline{\Pi}(t_k+\lambda_k^2 s,x_0+\lambda_k y),
\end{equation*}
with  $\underline{\mathfrak{p}}=\mathfrak{p}|_{Q_{\frac{R}{2}}(t_0,x_0)}$ and  $\underline{\Pi}=\Pi|_{Q_{\frac{R}{2}}(t_0,x_0)}$. Thus, by homogeneity, one has
\begin{equation*}
\|\mathfrak{p}_k\|_{L^\infty(]0,1[,L^\frac{3}{2}(\R))}=\|\lambda_{k}^2\underline{\mathfrak{p}}(t_k+\lambda_{k}^2\cdot,x_0+\lambda_{k}\cdot)\|_{L^\infty(]0,1[,L^\frac{3}{2}(\R))}=\|\underline{\mathfrak{p}}\|_{L^\infty(]t_k,t_0[,L^\frac{3}{2}(\R))}.
\end{equation*}
Since we have $]t_k,t_0[\subset]t_0-\frac{R^2}{4},t_0[$ by Remark \ref{Rem_SupportScaled}, we write
\begin{align}
\|\mathfrak{p}_k\|_{L^\infty(]0,1[,L^\frac{3}{2}(\R)}\le\|\underline{\mathfrak{p}}\|_{L^\infty(]t_0-\frac{R^2}{4},t_0[,L^\frac{3}{2}(B_{x_0,\frac{R}{2}}))}=\|\underline{\mathfrak{p}}\|_{L^\infty_tL^\frac{3}{2}_{x}(Q_\frac{R}{2}(t_0,x_0))}.\notag
\end{align}
Using the fact that  $\underline{\mathfrak{p}}=\mathfrak{p}|_{Q_{\frac{R}{2}}(t_0,x_0)}$ and  $ \mathfrak{p}\in L^\infty_t L^\frac{3}{2}_x(Q_{\frac{R}{2}}(t_0,x_0))$ by \eqref{PressureP1Estimate}, we obtain the following uniform bound
\begin{align}
\|\mathfrak{p}_k\|_{L^\infty(]0,1[,L^\frac{3}{2}(\R)}\le \|\mathfrak{p}\|_{L^\infty_tL^\frac{3}{2}_{x}(Q_\frac{R}{2}(t_0,x_0))}\le C<+\infty.\label{pfrakUniformEstimate}
\end{align}
Hence, by the Banach-Alaoglu theorem there exists a subsequence  $(\mathfrak{p}_{{k}_j})_{j\in \mathbb{N}}$ and $p_\infty\in L^\infty(]0,1[, L^\frac{3}{2}(\R))$ such that 
\begin{equation}\label{Convergencepfrak}
\mathfrak{p}_{k_j}\xrightarrow[j\to +\infty]{*}p_\infty \quad \text{in}\quad L^\infty(]0,1[,L^\frac{3}{2}(\R)).
\end{equation}
Let us study now the sequence $(\Pi_k)_{k\in \mathbb{N}}$. Since $\supp(\Pi_k)\subset Q_{\frac{R}{2\lambda_{k}}}(1,0)$, we have
\begin{align*}
\|\Pi_k\|_{L^\frac{3}{2}(]0,1[,L^\infty(\R))}&=\|\lambda_k^2\underline{\Pi}(t_k+\lambda_k^2 \cdot,x_0+\lambda_k \cdot)\|_{L^\frac{3}{2}(]0,1[,L^\infty(B_{0,\frac{R}{2\lambda_{k}}}))},
\end{align*}
and by the homogeneity of the space $L^\frac{3}{2}_tL^\infty_x$, we obtain
\begin{align*}
\|\Pi_k\|_{L^\frac{3}{2}(]0,1[,L^\infty(\R))}=\lambda_k^{\frac{2}{3}}\|\underline{\Pi}\|_{L^\frac{3}{2}(]t_k,t_0[,L^\infty(B_{x_0,\frac{R}{2}}))}\le \lambda_k^{\frac{2}{3}}\|\underline{\Pi}\|_{L^\frac{3}{2}(]t_0-\frac{R^2}{4},t_0[,L^\infty(B_{x_0,\frac{R}{2}}))} .
\end{align*}
Now, since $\underline{\Pi}=\Pi|_{Q_{\frac{R}{2}}(t_0,x_0)}$ and $\Pi\in L^\frac{3}{2}_tL^\infty_x(Q_{\frac{R}{2}})$ by \eqref{ControlPi}, one gets
\begin{align}\label{ConvergencePi} \|\Pi_k\|_{L^\frac{3}{2}(]0,1[,L^\infty(\R))}\le C\lambda_k^{\frac{2}{3}}\|\Pi\|_{L^\frac{3}{2}_tL^\infty_x(Q_{\frac{R}{2}}(t_0,x_0))}\le  C\lambda_k^{\frac{2}{3}}.
\end{align}  Since $(\lambda_{k})_{k\in \mathbb{N}}$ tends to zero as $k\to +\infty$ we can deduce that $(\Pi_k)_{k\in \mathbb{N}}$ converges to zero strongly in $L^\frac{3}{2}(]0,1[,L^\infty(\R))$.\\ 

We have proved so far that  $(\mathfrak{p}_{k_j})_{j\in\mathbb{N}}$ converges to $p_\infty$ by \eqref{Convergencepfrak} and  $(\Pi_k)_{k\in \mathbb{N}}$ tends to zero. Now, since $p_k=\mathfrak{p}_k+\Pi_k$, we may obtain, up to a subsequence, the weak convergence of $(p_k)_{k\in \mathbb{N}}$ to $p_\infty$ in $L^\infty_tL^3_x$. However, for our purposes we need to study more in detail the convergence of  $(p_k)_{k\in \mathbb{N}}$. Indeed, let us prove that  $(p_k)_{k\in \mathbb{N}}$ is uniformly bounded in $(L^\frac{3}{2}(]0,1[,L^\frac{3}{2}(\R)))_{loc}$. For showing this claim, we will use again the decomposition ${p_k=\mathfrak{p}_k+\Pi_k}$ and the previous estimates obtained on $\mathfrak{p}_k$ and $\Pi_k$.  Thus, for any compact set ${\mathfrak{Q}=[a,b]\times B\subset ]0,1[\times \R}$, since ${p_k=\mathfrak{p}_k+\Pi_k}$,  we have
\begin{align*}
\int_\mathfrak{Q}|p_k|^\frac{3}{2}dyds &\le C \int_0^1\int _{B}|\mathfrak{p}_k|^\frac{3}{2}dyds+C\int_0^1\int _{B}|\Pi_k|^\frac{3}{2}dyds.
\end{align*}
Since $\supp (p_k)\subset Q_{\frac{R}{2\lambda_{k}}(1,0)}$ and  $\lambda_{k}\xrightarrow[k\to +\infty]{}0$, we can consider $k$ large enough such that $\mathfrak{Q}\subset Q_{\frac{R}{2\lambda_{k}}}(1,0)$. Now, by using the fact that  $(\mathfrak{p}_k)_{k\in \mathbb{N}}$ is uniformly bounded in $L^\infty(]0,1[,L^\frac{3}{2}(\R))$ by \eqref{pfrakUniformEstimate}, we have
\begin{align*}
\int_0^1\int _{B}|\mathfrak{p}_k|^\frac{3}{2}dyds= \int_0^1\|\mathfrak{p}_k(s,\cdot)\|_{L^\frac{3}{2}(B)}^\frac{3}{2}ds \le C\|\mathfrak{p}_k\|_{L^\infty(]0,1[,L^\frac{3}{2}(\R))}^{\frac{3}{2}}\le C.
\end{align*}
Moreover, since  $\|\Pi_k\|_{L^\frac{3}{2}(]0,1[,L^\infty(\R))}\le C\lambda_k^{\frac{2}{3}} $ by \eqref{ConvergencePi}, we obtain
\begin{align*}
\int_0^1\int _{B}|\Pi_k|^\frac{3}{2}dyds\le |B|\int_0^1\|\Pi_k(s,\cdot)\|_{L^\infty(B)}^\frac{3}{2}ds\le  C\|\Pi_k\|_{L^\frac{3}{2}(]0,1[,L^\infty(\R))}^\frac{3}{2}\le  C\lambda_k^{\frac{2}{3}}<C, 
\end{align*}
where we have used that $(\lambda_{k})_{k\in \mathbb{N}}$ is a bounded sequence. Then, from the previous two estimates we obtain  
\begin{align}\label{UniformBoundPressure}
\int_\mathfrak{Q}|p_k|^\frac{3}{2}dyds&\le C.
\end{align}
Thus, $(p_k)_{k\in\mathbb{N}}$ is uniformly bounded in $(L^\frac{3}{2}(]0,1[,L^\frac{3}{2}(\R)))_{loc}$ and by the Banach-Alaoglu theorem and the uniqueness of the limit, there exists a subsequence $(p_{k_j})_{j\in \mathbb{N}}$ such that
\begin{equation}\label{ConvergencePressure}
p_{k_j}\xrightarrow[j\to +\infty]{*}p_\infty\quad\text{in}\quad (L^\frac{3}{2}_{t}L^\frac{3}{2}_{x})_{loc}.
\end{equation}
Notice that we have obtained a refinement of the  weak-$\ast$ convergence of $(p_k)_{k\in \mathbb{N}}$ given in (\ref{Convergencepfrak}).\\

\item [$\bullet$] Now, let us study the convergence of $(\vu_{k})_{k\in \mathbb{N}}$. First, observe that $(\vu_{k})_{k\in \mathbb{N}}$ is uniformly bounded in $L^\infty (]0,1[,L^3(\R))$, indeed using the definition of $\vu_k$ given in \eqref{DefUkWk}, by homogeneity and since $t-\frac{R^2}{4}<t_k$ by Remark \ref{Rem_SupportScaled},  we have
\begin{align*}
\|\vu_k\|_{L^\infty (]0,1[,L^3(\R))}&=\|\lambda_{k}\underline{\vu}(t_k+\lambda_{k}^2\cdot, x_0+\lambda_{k}\cdot)\|_{L^\infty (]0,1[,L^3(\R))}\notag\\
&=\|\underline{\vu}\|_{L^\infty (]t_k,t_0[,L^3(\R))}\le\|\underline{\vu}\|_{L^\infty (]t_0-\frac{R^2}{4},t_0[,L^3(\R))}.
\end{align*}
Then, since $\underline{\vu}=\vu|_{Q_{\frac{R}{2}}(t_0,x_0)}$ and $\vu\in L^\infty_tL^3_x(Q_R)$ by hypothesis, we have
\begin{equation}\label{UniformBoundL3}
\|\vu_k\|_{L^\infty (]0,1[,L^3(\R))}\le
\|\vu\|_{L^\infty_tL^3_x(Q_\frac{R}{2}(t_0,x_0))}<+\infty,
\end{equation}
then from the Banach-Alaoglu theorem, there exists a subsequence $(\vu_{k_j})_{j\in \mathbb{N}}$ such that
\begin{equation}\label{uinfiL3norm}
\vu_{k_j} \xrightarrow[j\to +\infty]{*}  \vu_\infty\in L^\infty (]0,1[,L^3(\R)).
\end{equation} 
Nevertheless, this convergence it is not enough to conclude that $(\vu_\infty,p_\infty)$ satisfies the Navier-Stokes equations \eqref{LimitNavierStokes} and we need to use the local energy inequality in order to obtain stronger convergences. For this purpose, we fix $\varphi\in \mathcal{C}_0^\infty(]0,1[\times \R)$ such that $\supp(\varphi)\subset ]a,b[\times B$ with $0<a<b<1$ and $B$ a bounded set of $\R$. For any $k\in \mathbb{N}$, we define  ${\varphi_k(\cdot,\cdot)=\varphi(\frac{\cdot-t_k}{\lambda_k^2},\frac{\cdot-x_0}{\lambda_k})}$.  Since the extended functions $(\underline{\vu},\underline{p},\underline{\vw})$ satisfy the local energy inequality \eqref{LocalEnergyInequality}, we have for any $t\in ]t_0-R^2,t_0[$,
\begin{eqnarray*}
\int_{\R} |\underline{\vu}|^2\varphi_k(t,x) dx+ 2 \int_{-\infty}^{t} \int _{\R}|\grad \otimes \underline{\vu}|^2\varphi_k dx ds\le \int_{-\infty}^{t} \int_{\R}(\partial_t\varphi_k +\Delta \varphi_k)|\underline{\vu}|^2dyds\qquad\qquad\qquad\\
+2\int_{-\infty}^{t} \int_{\R}\underline{p} (\underline{\vu}\cdot \grad \varphi_k ) dyds+\int_{-\infty}^{t }\int _{\R} |\underline{\vu}|^2(\underline{\vu} \cdot \grad)\varphi_k dx ds
+\int_{-\infty}^{t} \int_{\R}(\rot \underline{\vw})\cdot (\varphi_k\underline{\vu})dyds.
\end{eqnarray*}
By applying the change of variable $\tau=\frac{s -t_k }{\lambda_k^2},\;\; y=\frac{x-x_0}{\lambda_k}$ and since $\supp(\varphi)\subset ]a,b[\times B$, we have
\begin{eqnarray}
\int_{\R} |\vu_{k}|^2\varphi(\tau,y) dy+ 2 \int_{a}^b \int _{\R}|\grad \otimes \vu_{k}|^2\varphi dy d\tau
\le \underbrace{\int_{a}^b\int_{\R}(\partial_t\varphi +\Delta \varphi)|\vu_{k}|^2dyd\tau}_{(I)}\qquad \qquad \qquad\nonumber\\+2\underbrace{\int_{a}^b \int_{\R}p_k (\vu_{k}\cdot \grad \varphi ) dyd\tau}_{(II)}+\underbrace{\int_{a}^b \int _{\R} |\vu_{k}|^2(\vu_{k} \cdot \grad)\varphi dy d\tau}_{(III)}+\underbrace{\int_{a}^b \int_{\R}(\rot \vw_k)\cdot (\varphi\vu_{k})dyd\tau}_{(IV)}.\label{TermsToBound}
\end{eqnarray}
Since $\supp(\vu_k)\subset Q_{\frac{R}{2\lambda_{k}}}(1,0)=]1-\frac{R^2}{4\lambda_{k}^2},1[\times B_{0,\frac{R}{2\lambda_{k}}}$, we consider $k$ large enough such that ${ B\subset B_{0,\frac{R}{2\lambda_{k}}}}$ (recall the sequence $(\lambda_{k})_{k\in \mathbb{N}}$ converges towards zero). Now, our aim consists in obtaining uniform estimates of  $(\vu_k)_{k\in\mathbb{N}}$, for which we shall control each term of the right-hand side of \eqref{TermsToBound}.
\begin{itemize}
\item[$*$] For $(I)$ by the Hölder inequality $(1=\frac{1}{3}+\frac{2}{3})$, since $\supp(\varphi)\subset ]a,b[\times B\subset ]0,1[\times B$ and  $(\vu_k)_{k\in \mathbb{N}}$ is uniformly bounded in $L^\infty(]0,1[, L^3(\R))$ by \eqref{UniformBoundL3}, we obtain
\end{itemize}
\begin{align}
\int_{a}^b \int_{\R}(\partial_t\varphi +\Delta \varphi)|\vu_{k}|^2dyd\tau\le \int_{0}^1\|\partial_t\varphi+\Delta \varphi\|_{L^3(B)}\|\vu_{k}(\tau,\cdot)\|^2_{L^3(B)}d\tau\le C\|\vu_{k}\|^2_{L^\infty_tL_x^3}\le C.\label{TermI}
\end{align}
\begin{itemize}
\item[$*$] For the term $(II)$ in \eqref{TermsToBound}, by the Hölder inequality $(1=\frac{2}{3}+\frac{1}{3}+\frac{1}{\infty})$, we have
\begin{align*}
\int_{a}^b \int_{\R}p_k(\vu_{k}\cdot \grad \varphi ) dy d\tau\le 
\int_a^b\|p_k(\tau,\cdot)\|_{L^{\frac{3}{2}}(B)}\|\vu_k(\tau,\cdot)\|_{L^3(B)}\|\grad\varphi(\tau,\cdot)\|_{L^\infty(B)}d\tau.
\end{align*}
Then, by the Hölder inequality in the time variable $(1=\frac{2}{3}+\frac{1}{\infty}+\frac{1}{3})$, and since $(\vu_k)_{k\in \mathbb{N}}$ is uniformly bounded in $L^\infty(]0,1[, L^3(\R))$  by \eqref{UniformBoundL3} and  $(p_k)_{k\in \mathbb{N}}$ is uniformly bounded in $(L^\frac{3}{2}(]0,1[,L^\frac{3}{2}(\R)))_{loc}$ by \eqref{UniformBoundPressure}, one has
\begin{align}
\int_{a}^b \int_{\R}p_k(\vu_{k}\cdot \grad \varphi ) dy d\tau\le C \|p_k\|_{L_{t}^\frac{3}{2}L_{x}^\frac{3}{2}(]a,b[\times B)}\|\vu_{k}\|_{L^\infty_tL^3_x}\|\grad\varphi\|_{L^3_tL_x^\infty}\le C.\label{TermII}
\end{align}

\item[$*$]  The term $(III)$ in \eqref{TermsToBound} follows immediately from  the fact that  $(\vu_k)_{k\in \mathbb{N}}$ is uniformly bounded in $L^\infty (]0,1[,L^3(\R))$, indeed we have
\begin{align}
\int_{a}^b \int _{\R} |\vu_{k}|^2(\vu_{k} \cdot \grad)\varphi dy d\tau\le \int_{0}^1 \int _{\R} |\vu_{k}|^3 |\grad\varphi| dy d\tau\le 
C \|\grad \varphi\|_{L^\infty_{t,x}}\|\vu_{k}\|_{L^\infty_tL^3_x}^3 \le C .\label{TermIII}
\end{align}
\item[$*$] For the last term $(IV)$ of \eqref{TermsToBound},  by the Hölder inequality $(1=\frac{1}{2}+\frac{1}{6}+\frac{1}{3})$ we have
\begin{align*}
\int_{a}^b \int_{\R}(\rot \vw_k)\cdot (\varphi\vu_{k})dy d\tau&\le C\int_{0}^1\|\rot \vw_k(\tau,\cdot)\|_{L^2(B)}\|\varphi(\tau,\cdot)\|_{L^6(B)}\|\vu_k(\tau,\cdot)\|_{L^3(B)}d\tau,
\end{align*}
now, by the Cauchy-Schwarz inequality in the time variable we obtain
\begin{align*}
\int_{a}^b \int_{\R}(\rot \vw_k)\cdot (\varphi\vu_{k})dy d\tau& \le C \|\rot \vw_k\|_{L_{t,x}^{2}}\|\varphi\|_{L^2_tL^6_x}\|\vu_k\|_{L^\infty_t L_x^3}.
\end{align*}
Since,  $\|\rot \vw_k\|_{L^2_{t}L^2_x}\le \lambda_{k} \|\rot \vw\|_{L^2_{t}L^2_x}$  by \eqref{EstimateRotOmega}  and $\|\vu_k\|_{L^\infty _tL^3_x}\le C$ by \eqref{UniformBoundL3},  we obtain
\begin{align}\label{TermIV}
\int_{a}^b \int_{\R}(\rot \vw_k)\cdot (\varphi\vu_{k})dy d\tau \le C \lambda_{k}^{\frac{1}{2}}\|\rot \vw\|_{L_{t,x}^{2}}\le C,
\end{align}
where we have used the fact that $(\lambda_k)_{k\in\mathbb{N}}$ is a bounded sequence.\\
\end{itemize}
Thus, from the estimates \eqref{TermI}-\eqref{TermIV}, we deduce that there exists a constant $C>0$ ( independent of $k$) such that the left-hand side of \eqref{TermsToBound} satisfies 
$\displaystyle{\int_{\R} |\vu_{k}|^2\varphi dy +\int_{a}^b \int _{\R}|\grad \otimes \vu_{k}|^2\varphi dy d\tau\le C}$. Therefore, we obtain that for any test function $\varphi \in \mathcal{D}(]0,1[\times \R)$  the sequence 
\begin{equation}\label{UniformEstimateL2L6}
(\varphi\vu_k)_{k\in\mathbb{N}}\quad \text{ remains uniformly bounded  in } \quad L_t^\infty L_x^2\cap L_t^2\dot H_x^1.\\[1mm]
\end{equation}
Moreover, since $(\vu_k,p_k)$ satisfies the Navier-Stokes equations \eqref{ForceNS}, it is possible to obtain that $(\varphi\partial_t \vu_{k})_{k\in \mathbb{N}}$ remains uniformly bounded in $L^{\frac{3}{2}}_t H_x^{-\frac{3}{2}}$, (see for instance Step 3 in the proof of the Theorem 14.1 of the book \cite{PGLR1}). Thus, by the Rellich-Lions theorem (see \cite[Theorem 12.1]{PGLR1}), we may find a subsequence $(\vu_{{k_j}})_{j\in \mathbb{N}}$ such that in the domain $]0,1[\times \R$ we have on the one hand that $\vu_{{k_j}}$ converges weak-$*$ to $\vu_\infty$ in $(L_t^\infty L_x^2\cap L_t^2\dot H_x^1)_{loc}$ and on the other hand that $\vu_{{k_j}}$ is strongly convergent to $\vu_\infty$ in $(L^2_tL^2_x)_{loc}$. Furthermore, we can deduce that we have
\begin{equation}\label{ConvergenceL3loc}
\vu_{k_j}\underset{j\to+\infty}{\longrightarrow}\vu_\infty\quad \text{in} \quad (L^3_{t}L^3_x)_{loc},
\end{equation}
indeed, since  $(\vu_{k_j})_{j\in \mathbb{N}}$ is uniformly bounded  in $(L^\infty_tL^2_x\cap L^2_tL^6_x)_{loc}$ by \eqref{UniformEstimateL2L6} and in $L^\infty_tL^3_x$ by \eqref{UniformBoundL3}, using an interpolation argument we obtain $
\|\vu_{k_j}\|_{L^4_{t}L^4_x}\le \|\vu_{k_j}\|_{L^2_{t}L^6_x}^{\frac12}\|\vu_{k_j}\|_{L^\infty_{t}L^3_x}^{\frac12}\le C<+\infty$. Thus, $(\vu_{k_j})_{j\in \mathbb{N}}$ is uniformly bounded in $(L^4_{t}L^4_x)_{loc}$. This fact together with the strong convergence in $(L^2_{t}L^2_x)_{loc}$, which is given by the Rellich-Lions theorem, imply the strong convergence in $(L^3_{t}L^3_x)_{loc}$.\\

We have now ended the study of the sequence $(\vu_{k})_{k\in \mathbb{N}}$.\\
\end{itemize}
Summarizing, we have obtained, up to a subsequence, that the triplet $(\vu_k,p_k,\rot\vw_k)$ converges to $(\vu_\infty,p_\infty,0)$ in some (strong or weak) sense, from which we can deduce that $(\vu_\infty,p_\infty)$ is a weak solution of the Navier-Stokes equations in $]0,1[\times \R$
\begin{equation}\label{LimitNavierStokes}
\partial_t \vu_\infty=\Delta \vu_\infty -\div(\vu_\infty\otimes\vu_\infty)-\grad p_{\infty}.
\end{equation}
Moreover, from the weak-$\ast$ convergence  of $(p_{k_j})_{j\in \mathbb{N}}$ given in \eqref{ConvergencePressure} and the strong convergence of $(\vu_{k_j})_{j\in \mathbb{N}}$ in $(L^3_tL^3_x)_{loc}$ given in \eqref{ConvergenceL3loc}, it is possible to deduce  that $(\vu_\infty,p_\infty)$ is a \emph{suitable solution} of the Navier-Stokes equation in $]0,1[\times \R$, in the sense of Definition 6.9 of \cite{PGLR1}. This fact can be seen with all details in the Step 4 in the proof of Theorem 14.1 of the book \cite{PGLR1}.\\

We will exploit this ``suitable'' property in the sequel but we need more properties on $(\vu_\infty, p_\infty)$ in order to obtain the desired contradiction. Indeed, we will show that on one hand, this solution is a nontrivial solution of the Navier-Stokes equations, and on the other hand, using the backward uniqueness theory developed in \cite{Eus03}, \cite{PGLR1} or \cite{Seregin14}, we will deduce that this solution must be identically zero, leading us to the wished contradiction.\\

Let us prove now that $(\vu_\infty, p_\infty)$ is nontrivial and for this we will study in a particular manner the limit that leaded us to  $(\vu_\infty, p_\infty)$. Consider first $0<\mathfrak{a}\ll \frac{1}{2}$ a small parameter.  Since $(\vu,p,\vw)$ is a partial suitable solution with $(t_0,x_0)$  a partial singular point, we can consider $k$ big enough such that $0<\mathfrak{a} \lambda_{k}\le \frac{R}{4}$ and we can use Proposition \ref{Propo_CaracSingular} in the appendix \ref{Sec_PartialRegularity} (which is valid for all radius $0<\mathfrak{a} \lambda_{k}<1$) to obtain the existence of a small parameter $\varepsilon>0$ such that 
\begin{equation*}
0<\varepsilon< \frac{1}{(\mathfrak{a}\lambda_k)^2}\int_{Q_{\mathfrak{a}\lambda_k}(t_0,x_0)}|\vu|^3+|p|^{\frac{3}{2}}dyds.
\end{equation*}
Now, observe that by a change of variable, we have
\begin{equation}\label{PositiviUk}
\varepsilon<\frac{1}{(\mathfrak{a}\lambda_k)^2}\int_{Q_{\mathfrak{a}\lambda_k}(t_0,x_0)}|\vu|^3+|p|^\frac{3}{2}dyds=\frac{1}{\mathfrak{a}^2}\int_{1-\mathfrak{a}^2}^{1}\int_{B_{0,\mathfrak{a}}}|\vu_k|^3+|p_k|^{\frac{3}{2}}dyds.
\end{equation}
For studying more in detail the previous expression we need to obtain some estimates about the pressure $(p_k)_{k\in \mathbb{N}}$.
Since $(\vu_k,p_k)$ is a solution of the system  \eqref{ForceNS} which can be seen as the Navier-Stokes equations with an external force $\rot\vw_k$, the pressure satisfies $-\Delta p_k=\div(\div(\vu_k\otimes\vu_k))$. Hence, following the same arguments as in \eqref{PresureDecomposition}, we can split the pressure  $p_k=\widetilde{\mathfrak{p}}_k+\widetilde \Pi_k $ where ${\widetilde{\mathfrak{p}}_k=\frac{1}{(-\Delta)}(\div( \div (\phi \vu_k\otimes \vu_k))}$  with  $\phi$ a positive test function supported in $B_{0,2\mathfrak{a}}$ such that $\phi=1$ in $B_{0,\mathfrak{a}},$ and $\widetilde{\Pi}_k$ is a harmonic function defined by $\widetilde{\Pi}_k={p}_k-\widetilde{\mathfrak{p}}_k$. Now, using the boundedness of the Riesz transform in $L^\frac32(\R)$, we have
\begin{align*}
\|\widetilde{\mathfrak{p}}_k(t,\cdot)\|_{L^\frac{3}{2}_{x}(B_{0,2\mathfrak{a}})}&\le\|\widetilde{\mathfrak{p}}_k(t,\cdot)\|_{L^\frac{3}{2}(\R)}=\left\|\frac{1}{(-\Delta)}(\div( \div (\phi \vu_k\otimes \vu_k))(t,\cdot)\right\|_{L^\frac{3}{2}(\R)}\notag\\
&\le C \|\phi (\vu_k\otimes \vu_k)(t,\cdot)\|^2_{L^\frac{3}{2}(\R)}.
\end{align*}  
Moreover since $\supp(\phi)\subset B_{0,2\mathfrak{a}}$, we have 
\begin{align}\label{PressurePkEstimate}
\|\widetilde{\mathfrak{p}}_k(t,\cdot)\|_{L^\frac{3}{2}_{x}(B_{0,2\mathfrak{a}})}\le C  \|\phi\|_{L^\infty(\R)}\|\vu_k(t,\cdot)\|^2_{L^3(B_{0,2\mathfrak{a}})}\le C\|\vu_k(t,\cdot)\|^2_{L^3(B_{0,2\mathfrak{a}})}.
\end{align}
Furthermore,  since $\widetilde{\Pi}_k$ is a harmonic function, by the same arguments as in \eqref{PressureHarmonicFunction}, {\textit{i.e.},} the local estimates for harmonic functions,  we obtain 
\begin{align*}
\|\widetilde{\Pi}_k(t,\cdot)\|_{L^\infty(B_{0,\mathfrak{a}})}\le C\|\widetilde{\Pi}_k(t,\cdot)\|_{L^\frac{3}{2}( B_{0,2\mathfrak{a}})}\le C\|\widetilde{\mathfrak{p}}_k(t,\cdot)\|_{L^\frac{3}{2}(B_{0,2\mathfrak{a}})}+C\|{p}_k(t,\cdot)\|_{L^\frac{3}{2}(B_{0,2\mathfrak{a}})},
\end{align*}
and using the estimate \eqref{PressurePkEstimate}, we have 
\begin{align}
\|\widetilde{\Pi}_k(t,\cdot)\|_{L^\infty(B_{0,\mathfrak{a}})}
\le C \|\vu_k(t,\cdot)\|^2_{L^3(B_{0,2\mathfrak{a}})}+C\|{p}_k(t,\cdot)\|_{L^\frac{3}{2}(B_{0,2\mathfrak{a}})} .\label{Pi_kEstimate}
\end{align} 
Now, coming back to \eqref{PositiviUk}, using that $p_k=\widetilde{ \mathfrak{p}}_k+\widetilde{\Pi}_k$, we have
\begin{align}
\varepsilon&<\frac{1}{\mathfrak{a}^2}\int_{1-\mathfrak{a}^2}^1\int_{B_{0,\mathfrak{a}}}|\vu_k|^3+|p_k|^{\frac{3}{2}}dyds\le\frac{1}{\mathfrak{a}^2}\int_{1-\mathfrak{a}^2}^1\int_{B_{0,\mathfrak{a}}}|\vu_k|^3+C|\widetilde {\mathfrak{p}}_k|^{\frac{3}{2}}+C|\widetilde \Pi_k|^{\frac{3}{2}}dyds \notag\\
\varepsilon&< \frac{1}{\mathfrak{a}^2}\int_{1-\mathfrak{a}^2}^1\int_{B_{0,\mathfrak{a}}}|\vu_k|^3dyds +\frac{C}{\mathfrak{a}^2}\int_{1-\mathfrak{a}^2}^1\int_{B_{0,\mathfrak{a}}}|\widetilde{ \mathfrak{p}}_k|^{\frac{3}{2}}dyds+\frac{C}{\mathfrak{a}^2}\int_{1-\mathfrak{a}^2}^1\int_{B_{0,\mathfrak{a}}}|\widetilde{\Pi}_k|^{\frac{3}{2}}dyds. \label{PositivityUk}
\end{align}
Let us study each term of the expression above. For the first one, it is easy to see that
\begin{equation}\label{FirstTermPositivUk}
\frac{1}{\mathfrak{a}^2}\int_{1-\mathfrak{a}^2}^1\int_{B_{0,\mathfrak{a}}}|\vu_k|^3dyds\le \frac{1}{\mathfrak{a}^2}\int_{1-(2\mathfrak{a})^2}^1\int_{B_{0,2\mathfrak{a}}}|\vu_k|^3dyds= \frac{1}{\mathfrak{a}^2}\|\vu_k\|_{L^3_tL^3_x(Q_{2\mathfrak{a}}(1,0))}^3.
\end{equation}
For  the second term in \eqref{PositivityUk}, from the estimate \eqref{PressurePkEstimate}, we have
\begin{align}
\frac{1}{\mathfrak{a}^2}\int_{1-\mathfrak{a}^2}^1\int_{B_{0,\mathfrak{a}}}|\widetilde{ \mathfrak{p}}_k|^{\frac{3}{2}}dyds&\le \frac{1}{\mathfrak{a}^2} \int_{1-(2\mathfrak{a})^2}^1\|\widetilde{\mathfrak{p}}_k(s,\cdot)\|_{L^\frac{3}{2}(B_{0,2\mathfrak{a}})}^{\frac{3}{2}}ds\le \frac{C}{\mathfrak{a}^2}\int_{1-(2\mathfrak{a})^2}^1 \left\|\vu_k(s,\cdot)\right\|_{L^3(B_{0,2\mathfrak{a}})}^{3}ds\notag\\
&\le\frac{C}{\mathfrak{a}^2}\|\vu_k\|_{L^3_tL^3_x(Q_{2\mathfrak{a}}(1,0))}^3.\label{SecondTermPositivityUk}
\end{align}
For the third term in \eqref{PositivityUk}, since $\|\widetilde \Pi_k(s,\cdot)\|_{L^\infty(B_{0,a})}\le C\|\vu_k(s,\cdot)\|^2_{L^3(B_{0,2\mathfrak{a}})}+C\|{p}_k(s,\cdot)\|_{L^\frac{3}{2}(B_{0,2\mathfrak{a}})}$ by \eqref{Pi_kEstimate},
we obtain
\begin{align*}
\frac{1}{\mathfrak{a}^2}\int_{1-\mathfrak{a}^2}^1\int_{B_{0,\mathfrak{a}}|}\widetilde{\Pi}_k|^{\frac{3}{2}}dyds&\le C\mathfrak{a}\int_{1-\mathfrak{a}^2}^1\|\widetilde{\Pi}_k(t,\cdot)\|_{L^\infty(B_{0,\mathfrak{a}})}^{\frac{3}{2}}ds\\
&\le C\mathfrak{a}\left(\int_{1-\mathfrak{a}^2}^1\|\vu_k(s,\cdot)\|^3_{L^3(B_{0,2\mathfrak{a}})}+\int_{1-\mathfrak{a}^2}^1\|{p}_k(s,\cdot)\|_{L^\frac{3}{2}(B_{0,2\mathfrak{a}})}^{\frac32}ds\right)\\
& \le C\mathfrak{a}(\|\vu_k\|^3_{L^\infty_tL_x^3}+\|{p}_k\|_{L^\frac{3}{2}_{t}L^\frac{3}{2}_{x}(Q_{2\mathfrak{a}}(1,0))}^{\frac{3}{2}}).
\end{align*}
Thus, since $(\vu_k)_{k\in \mathbb{N}}$ is uniformly bounded  in $L^\infty(]0,1[, L^3(\R))$ by \eqref{UniformBoundL3} and $(p_k)_{k\in\mathbb{N}}$ is uniformly bounded in $(L^\frac{3}{2}(]0,1[,L^\frac{3}{2}(\R)))_{loc}$ by \eqref{UniformBoundPressure}, we obtain that
\begin{equation}\label{ThirdTermPositivyUk}
\frac{1}{\mathfrak{a}^2}\int_{1-\mathfrak{a}^2}^1\int_{B_{0,\mathfrak{a}}}|\widetilde{\Pi}_k|^{\frac{3}{2}}dyds\le C\mathfrak{a}.
\end{equation}
Then, gathering all the estimates \eqref{FirstTermPositivUk}-\eqref{ThirdTermPositivyUk} in \eqref{PositivityUk}, we obtain
$$\varepsilon< \frac{C}{\mathfrak{a}^2}\|\vu_k\|_{L^3_tL^3_x(Q_{2\mathfrak{a}}(1,0))}^3+C'\mathfrak{a},$$
which we rewrite in the following manner $\mathfrak{a}^2\varepsilon-C'\mathfrak{a}^3< C \|\vu_k\|_{L^3_tL^3_x(Q_{2\mathfrak{a}}(1,0))}^3$. Now, by considering $\mathfrak{a}$ such that $\mathfrak{a}<\frac{\varepsilon}{C'}$, we can find  a constant $0<\varepsilon_\ast<\mathfrak{a}^2\varepsilon-C'\mathfrak{a}^3$, such that
\begin{equation*}
0<\varepsilon_\ast < C\|\vu_k\|^3_{L^3_{t,x}(Q_{2\mathfrak{a}}(1,0))}.
\end{equation*}
Thus, from the strong convergence in $(L^3_{t,x})_{loc}$ of $(\vu_{k_j})_{j\in\mathbb{N}}$ given in \eqref{ConvergenceL3loc}, we obtain
\begin{equation}\label{NonTriviality}
0<\varepsilon_\ast<\int_{Q_{2\mathfrak{a}}}|\vu_\infty|^3dyds.
\end{equation}
We have thus proven that $(\vu_\infty,p_\infty)$ is a nontrivial solution of the Navier-Stokes equations.\\

We will now exhibit a contradiction by showing that $\vu_\infty\equiv 0$. For this purpose,  we recall that the limit solution $(\vu_\infty,p_\infty)$ satisfies the Navier-Stokes equations \eqref{LimitNavierStokes} and therefore we may consider the backwards uniqueness and unique continuation theories developed in \cite{Eus03}, which can be summarized in the following proposition
%%%%%%%%%%%%%%%%%%%%%%%%%%%%%%%%%%%%%%%%%%%%%%%%%%%
\begin{Proposition}\label{Prop_Uniqueconti}
Let $(\vv,h)$ be a solution of the Navier-Stokes equations on $]0,1[\times \R$,  \textit{i.e.,} we have
\begin{equation*}
\partial \vv=\Delta \vv -(\vv\cdot\grad )\vv-\grad h,\quad \div(\vv)=0
\end{equation*} 
Assume moreover  that $\vv\in L_t^\infty L_x^2\cap L_t^2\dot H_x^1$ and for any $\phi \in \mathcal{C}_0^\infty(]0,1[\times \R)$, the pair $(\vv,h)$ satisfies the following local energy inequality 
\begin{align}\label{LocalEnergyNS}
\int_{\R} |\vv|^2\phi(t,\cdot) dy+ &2 \int_{s<t} \int _{\R}|\grad \otimes \vv|^2\phi dy ds\le \int_{s<t}\int_{\R}(\partial_t\phi + \Delta \phi)|\vv|^2dyds\\
&+2\int_{s<t}\int_{\R}h (\vv\cdot \grad \phi ) dyds +\int_{s<t}\int _{\R} |\vv|^2(\vv \cdot \grad)\phi dy ds.\notag
\end{align}
If $\vv \in L^\infty(]0,1[,L^3(\R))$ and $\vv(1,\cdot)=0$, then $\vv=0$ on $]0,1]\times \R$. 
\end{Proposition}
%%%%%%%%%%%%%%%%%%%%%%%%%%%%%%%%%%%%%%%%%%%%%%%%%%%
\noindent For a proof of this proposition we refer to the article \cite{Eus03} and the books \cite{Seregin14}, \cite{Tsai19}. \\

Let us now verify that the pair $(\vu_\infty,p_\infty)$ satisfies the hypotheses of the previous proposition. First, notice that $(\vu_\infty,p_\infty)$ satisfies the local energy inequality \eqref{LocalEnergyNS} since it is a suitable solution of the Navier-Stokes equations \eqref{LimitNavierStokes}. Moreover we also have that $\vu\in L^\infty(]0,1[,L^3(\R))$ by \eqref{uinfiL3norm}, thus it is enough to proof that $\vu_{\infty}(1,\cdot)=0$. For this purpose, remark that for any $j\in \mathbb{N}$, we have ${\vu_{k_j}\in L^\infty(]0,1[, L^2(\R))\subset L^1(]0,1[,H^{-\frac{3}{2}}(\R))}$ due to the spaces inclusions ${L^2(\R)\subset H^{-1}(\R)\subset H^{-\frac{3}{2}}(\R)}$ and that we have  ${\partial _t\vu_{k_j}\in L^\frac{3}{2}(]0,1[,H^{-\frac{3}{2}}(\R))\subset L^1(]0,1[,H^{-\frac{3}{2}}(\R))}.$ Therefore, by following the same lines that leaded us to deduce \eqref{ContinuityInTimePartialSuitableSolution}, (see as well \cite[pg 402]{PGLR1}), we can obtain that $${\vu_{k_j}\in \mathcal{C}([0,1],H^{-\frac{3}{2}}(\R))}.$$
Thus, if we consider   $\phi\in \mathcal{C}^\infty(\mathbb{R})$ such that $\phi=1$ on $]-\infty,\frac{3}{2}[$ and $\phi=0$ on $]2,+\infty[$,  by writing  for any $t\in [0,1]$, 
$$\vu_{k_j}(t,\cdot)=-\int_t^2 \partial_t(\phi\vu_{k_j})ds, $$
we can obtain that $\vu_\infty(t,\cdot)$ is the  weak-$*$ limit of $\vu_{k_j}(t,\cdot)$  in $H^{-\frac{3}{2}}(\R)$. It follows that for any $t\in [0,1]$, $\vu_{\infty}(t,\cdot)$ is well defined in a distributional sense.
In particular, $\vu_{\infty}(1,\cdot)$ is the weak-$\ast$ limit of $\vu_{k_j}(1,\cdot)$  in $H^{-\frac{3}{2}}(\R)$. Moreover, for any $\mathfrak{r}>0$, since $\vu_k(1,\cdot)=\lambda_{k}\vu(t_0,x_0+\lambda_k\cdot)$ and by the change of variable $z=x_0+\lambda_{k_j}y$, we have
\begin{align*}
\int_{B_{0,\mathfrak{r}}}|\vu_{k_j}(1,y)|^3dy=\int_{B_{0,\mathfrak{r}}}\lambda_{k_j}^3|\underline{\vu}(t_0,x_0+\lambda_{k_j}y)|^3dy=\int_{B_{x_0,\lambda_{k_j}\mathfrak{r}}}|\underline{\vu}(t_0,z)|^3dz.
\end{align*}
Notice that  $\vu(t_0,\cdot)\in L^3(B_{x_0,\frac{R}{2}})$ by \eqref{EveryTimeInL3} and since $\underline{\vu}=\vu|_{Q_{\frac{R}{2}}(t_0,x_0)}$, we have that $\underline{\vu}(t_0,\cdot)\in L^3(\R)$. Thus, since  $\lambda_{k_j}=\sqrt{t_0-t_{k_j}}\underset{j \to+\infty}{\longrightarrow}0$, the function $\mathds{1}_{B_{x_0,\lambda_{k_j}\mathfrak{r}}}$  converges pointwise  to $0$ as $j\to+\infty $. Hence by the  dominated convergence theorem, we have 
\begin{equation*}
\int_{B_{0,\mathfrak{r}}}|\vu_{k_j}(1,y)|^3dy=\int_{\R}\mathds{1}_{B_{x_0,\lambda_{k_j}\mathfrak{r}}}(z)|\underline{\vu}(t_0,z)|^3dz\underset{j \to+\infty}{\longrightarrow}0.
\end{equation*}
Therefore, the sequence $(\vu_{k_j}(1,\cdot))_{j\in \mathbb{N}}$ converges weakly-$\ast$ to $0$ in $L^3(\R)$ and then by the uniqueness of the limit, we obtain that $\vu_\infty(1,\cdot)=0$.\\

\noindent We have now all the hypotheses needed to apply  Proposition \ref{Prop_Uniqueconti} (\emph{i.e.} $(\vu_\infty,p_\infty)$ is a suitable solution, $\vu_\infty\in L^\infty(]0,1[,L^3(\R))$ by \eqref{uinfiL3norm} and $\vu_\infty(1,\cdot)=0$) so we obtain that $\vu_\infty=0$ on $]0,1[\times \R$. Thus, we have  $$\displaystyle {\int_{Q_{2\mathfrak{a}}}|\vu_\infty| dx ds=0,}$$ which is a contradiction to \eqref{NonTriviality}, and this allows to conclude that $\vu$ is actually bounded in $Q_{r}(t_0,x_0)$ for all $0<r\le \frac{R}{2}$. This ends the proof of the Theorem \ref{Theo_LocalPoints}. \hfill$\blacksquare $
%%%%%%%%%%%%%%%%%%%%%%%%%%%%%%%%%%%%%%%%%%%%%%%%%%%
%%%%%%%%%%%%%%%%%%%%%%%%%%%%%%%%%%%%%%%%%%%%%%%%%%%
\mysection{Proof of Theorem \ref{Theo_NormL3}}\label{Sec_RegularityNormL3}
We recall the setting of this theorem. Let $(\vu,p,\vw)$ be a Leray-type weak solution of the micropolar fluids equations \eqref{MicropolarFluidsEquationsEqua1} and \eqref{MicropolarFluidsEquationsEqua2} such that $\vu,\vw \in L^\infty(]0,+\infty[,L^2(\R))\cap L^2(]0,+\infty[,\dot{H}^1(\R))$ and we assume that for  ${0<\delta<T<+\infty},$ we have ${\vu\in L^\infty(]\delta,T[, L^3(\R))}$. Our aim consists in proving that under the previous assumptions the condition ${\vu \in \mathcal{C}(]\delta, T[,L^3(\R))}$ is equivalent to the fact that any point $(t_0,x_0)\in ]\delta,T[$ is partially regular in the sense of Definition \ref{Def_PartialRegularPoints}.\\

\noindent To do so, first we will establish some properties of the weak solution $(\vu,p,\vw)$ in this framework.
\begin{itemize}
\item [$\bullet$] We prove here that $\vu(t,\cdot)\in L^3(\R)$ for any $t\in[\delta,T]$. For showing this claim, in contrast to the proof of Theorem \ref{Theo_LocalPoints}, we take advantage of the properties of weak solutions of the Navier-Stokes equations. Indeed, since $(\vu,p)$ is a weak solution of the equation \eqref{MicropolarFluidsEquationsEqua1} which can be seen as the Navier-Stokes equations with an external force $\rot \vw\in L^2_{t,x}$, it is possible to deduce that  $\vu$ is $L^2(\R)$-weakly continuous in time \textit{i.e.,} for any  $t\in [\delta,T]$ the application 
\begin{equation}\label{WeaklyC}
t\longmapsto \int_{\R} \vu(t,x)\vphi (x)dx, 
\end{equation}
is continuous for every $\vphi\in L^{2}(\R)$.  See for instance Theorem 3.8 in the book \cite{Robinson} or Lemma 3.4 in \cite{Tsai19} for a proof of this fact and more details. Now, fix $t\in [\delta,T]$ and we consider a sequence $(t_k)_{k\in\mathbb{N}}$ in $]\delta,T[$ such that $t_k\underset{k\to+\infty}{\longrightarrow}t$. Since ${\|\vu(t_k,\cdot)\|_{L^3}\le \|\vu\|_{L^\infty_t L^3_x}}$, by the Banach-Alaoglu theorem, there exists a subsequence $(\vu(t_{k_j},\cdot))_{j\in \mathbb{N}}$ such that $\vu(t_{k_j},\cdot)\longrightarrow \vec v(t,\cdot)$ weakly-$\ast$ in $L^3(\R)$. 
On the other hand, since the application \eqref{WeaklyC} is continuous for every $\phi \in L^2(\R)$,  in particular it is continuous for any $\psi\in \mathcal{C}_0^{\infty}(\R)$, and hence we have
$\displaystyle{\int_{\R}\vu(t_{k_j},x)\psi(x)dx\underset{j\to+\infty}{\longrightarrow}\int_{\R}\vu(t,x)\psi(x)dx}$. Since $\mathcal{C}_0^{\infty}(\R)$ is dense in $L^\frac{3}{2}(\R)$, by the uniqueness of the limit, we obtain $\vu(t,\cdot)=\vec v(t,\cdot)\in L^3(\R)$. We have thus proved that
\begin{equation}\label{EverytimeinL3Global}
\text{for any} \quad  t\in [\delta, T], \quad \vu(t,\cdot)\in L^3(\R).
\end{equation} 

\item [$\bullet$] We prove now that, for any open set $B\subset\R$, the triplet $(\vu,p,\vw)$ is a partial suitable solution of the micropolar fluids equations in $]\delta,T[\times B$ in the sense of Definition \ref{Def_PartialSuitable}. Indeed, since  ${\vu,\vw\in L^\infty(]0,+\infty[L^2(\R))\cap L^2(]0,+\infty[,\dot H^1(\R))}$ we immediately have $${\vu,\vw\in L^\infty(]\delta,T[,L^2( B))\cap L^2(]\delta,T[,\dot H^1( B))}.$$ Thus, it is enough to show that $p\in L^\frac 32(]\delta,T[,L^\frac 32(B))$ and $(\vu,p,\vw)$ satisfies the following local energy inequality: for any $\phi \in \mathcal{C}_0^\infty(]\delta,T[\times B),$ 
\begin{align}
\int_{\R} |\vu|^2\phi(t,\cdot) dx
&+ 2 \int_{s<t} \int _{\R}|\grad \otimes \vu|^2\phi dx ds
\le \int_{s<t}\int_{\R}(\partial_t\phi + \Delta\phi)|\vu|^2dyds+2\int_{s<t}\int_{\R}p (\vu\cdot \grad \phi ) dyds\nonumber\\
&+\int_{s<t}\int _{\R} |\vu|^2(\vu \cdot \grad)\phi dx ds+
\int_{s<t}\int_{\R}(\rot \vw)\cdot (\phi\vu)dyds.\label{LocalEnergyInequalityTheorem2}
\end{align}
For proving that $p\in L^\frac 32(]\delta,T[,L^\frac 32(B))$, recall that the pressure satisfies the equation ${p=\frac{1}{(-\Delta)}\div\div(\vu\otimes\vu)}$ over $\R$. Hence, using the boundedness of the Riesz transforms in $L^\frac{3}{2}(\R)$, we have 
\begin{equation*}
\|p(t,\cdot)\|_{L^{\frac{3}{2}}(B)}\le\left\|\frac{1}{(-\Delta)}\div\div(\vu\otimes\vu)(t,\cdot)\right\|_{L^{\frac{3}{2}}(\R)}\le C\|\vu\otimes\vu(t,\cdot)\|_{L^{\frac{3}{2}}(\R)}\le C \|\vu(t,\cdot)\|_{L^3(\R)}^2.
\end{equation*} 
By considering the $L^\frac{3}{2}$-norm in the time interval $]\delta,T[$, in the expression above, and since ${\vu\in L^\infty(]\delta,T[,L^3(\R))}$ by hypothesis, one has 
\begin{equation}\label{PressureInL32}
\|p\|_{L_t^{\frac{3}{2}}L_x^{\frac{3}{2}}(]\delta,T[\times B)}\le C \|\vu\|_{L^{\frac{3}{2}}(]\delta,T[,L^3(\R))}^2\le C \|\vu\|_{L^\infty(]\delta,T[,L^3(\R))}^2<+\infty.
\end{equation} 

Now, let us prove that $(\vu,p,\vw)$ satisfies the local energy inequality \eqref{LocalEnergyInequalityTheorem2}. First notice that since  ${\vu\in L_t^2\dot H^1_x\subset L_t^2L_x^6}$ and $\vu \in L^\infty(]\delta,T[ ,L^3(\R))$ by hypothesis, by using an interpolation argument, we have
\begin{equation*}
\|\vu\|_{ L^4(]\delta,T[ ,L^4(B))}\le \|\vu\|_{ L^\infty(]\delta,T[ ,L^3(B))}^{\frac{1}{2}}\|\vu\|_{ L^2(]\delta,T[ ,L^6(B))}^\frac{1}{2}\le  \|\vu\|_{ L^\infty(]\delta,T[ ,L^3(\R))}^{\frac{1}{2}}\|\vu\|_{ L^2(]\delta,T[ ,L^6(\R))}^\frac{1}{2}<+\infty.
\end{equation*} 
Thus, since $(\vu,p,\vw)$ satisfies the first equation of the micropolar fluids equations $\eqref{MicropolarFluidsEquationsEqua1}$, and we have deduced that $\vu\in L^4(]\delta,T[,L^4(B))$ and $p\in L^\frac{3}{2}(]\delta,T[,L^\frac{3}{2}(B))$, it is then possible to see that each term in the local energy inequality \eqref{LocalEnergyInequalityTheorem2} is well defined. Therefore, since ${\vu,\vw\in L^\infty(]\delta,T[,L^2( B))\cap L^2(]\delta,T[,\dot H^1( B))}$, $p\in L^\frac{3}{2}(]\delta,T[,L^\frac{3}{2}(B)),$ and  the local energy inequality is satisfied, we obtain that for any open set $B\subset\R$,  the triplet $(\vu,p,\vw)$ is a partial suitable solution on $]\delta,T[\times B$.\\
\end{itemize}
Having proved the previous two points about the weak solution $(\vu,p,\vw)$ in the general framework considered in this section, we continue with the proof of Theorem \ref{Theo_NormL3}.\\ 

First let us show that if $\vu \in \mathcal{C}(]\delta,T[,L^3(\R))\cap L^\infty(]\delta, T[,L^3(\R))$ then any $(t_0,x_0)\in ]\delta, T[\times \R$ is partially regular in the sense of Definition \ref{Def_PartialRegularPoints}. Indeed, notice that since $\vu\in  L^\infty(]\delta, T[,L^3(\R)$, for any $(t_0,x_0)\in ]\delta,T[\times \R$, there exists $0<R<\sqrt{t_0-\delta}$ such that $\vu\in L^\infty_tL^3_x(Q_R(t_0,x_0))$. Moreover, we have seen that $(\vu,p,\vw)$ is a partial suitable solution on $ Q_R(t_0,x_0)$ in the sense of Definition \ref{Def_PartialRegularPoints}, then by using Theorem \ref{Theo_LocalPoints}, there exists a radius $0<r\le \frac{R}{2}$ such that $\vu\in L^\infty_{t,x}(Q_r(t_0,x_0))$ and therefore the point $(t_0,x_0)$ is partially regular in the sense of Definition \ref{Def_PartialRegularPoints}. Thus we have proved the first implication.\\

We turn now to the other direction: assume that any $(t_0,x_0)\in ]\delta, T[\times \R $ is a partial regular point in the sense of Definition \ref{Def_PartialRegularPoints} and we aim to prove that 
we have 
$$\vu\in \mathcal{C}(]\delta, T[, L^3(\R)).$$ 
To do so, first, we will deduce that the velocity $\vu$ satisfies that $\vu\in L^\infty_{t,x}(]\delta,T[\times \R)$ and   $\vu\in \mathcal{C}([\delta,T],L^2(\R))$.\\ 

Indeed, let us prove that $\vu$ is bounded on $]\delta, T[\times \R$. Fix $t_0\in ]\delta, T[$ and $0<R<\sqrt{t_0-\delta}$. Since ${\vu\in L^\infty (]\delta,T[,L^3(\R))}$ by hypothesis and $p\in L^\frac{3}{2}(]\delta,T[,L^\frac{3}{2}(\R))$ by \eqref{PressureInL32}, we have that    
$$\lim_{|x_0|\longrightarrow +\infty}\frac{1}{R^2}\int_{Q_R(t_0,x_0)}|\vu|^3+|p|^\frac32dyds=0.$$
Thus, for  $\varepsilon>0$ there exists $K>0$  such that for any $|x_0|>K$, we have $\displaystyle{\frac{1}{R^2}\int_{Q_R(t_0,x_0)}|\vu|^3+|p|^\frac32dyds\le \varepsilon}$. Therefore, using the $\varepsilon$- regularity theory developed in Theorem \ref{Theo_partialRegularity} in the appendix \ref{Sec_PartialRegularity}, there exists $0<r\le \frac{R}{2}$ such that ${\|\vu\|_{L^\infty_{t,x}(Q_r(t_0,x_0))}<C}$. Since this bound is valid for any $|x_0|>K$, we deduce that
\begin{equation}\label{BoundOutsideBk}
\vu\in L^{\infty}_{t,x}(]t_0-r^2,t_0[\times  B_{0,K}^c).
\end{equation}
Now, we will show that  $\vu$ is bounded on $]t_0-\rho^2,t_0[\times \bar B_{0,K}$, for some $\rho>0$ to be defined later. Notice that for any $y\in \bar B_{0,K}$, the point $(t_0,y)$ is partially regular by hypothesis and hence there exists $0<\mathfrak{r}_y<\sqrt{t_0-\delta}$ such that 
\begin{equation}\label{UboundedBalls}
\vu\in L^\infty_{t,x}(Q_{\mathfrak{r}_{y}}(t_0,y)),\quad \text{ where}\quad Q_{\mathfrak{r}_y}(t_0,y)=]t_0-\mathfrak{r}^2_y,t_0[\times B_{y,\mathfrak{r}_y}.
\end{equation}
Remark that the family $\{B_{y,\mathfrak{r}_y}: y\in \bar B_{0,K}\}$ forms a cover of $\bar B_{0,K}$. Thus, by the compactness of $\bar B_{0,K}$ and by the information given in \eqref{UboundedBalls}, there exists a finite sub-cover $\{B_{\mathfrak{r}_{y_i}}(t_0,y_i): i=1,\dots,n\}$ of $\bar B_{0,K}$ such that for all $1\le i\le n$, $\vu\in L^\infty_{t,x}(Q_{\mathfrak{r}_{y_i}}(t_0,y_i))$. Setting ${\rho=\min \{ \mathfrak{r}_{y_1},\dots,\mathfrak{r}_{y_n}\}}$, we have $\displaystyle{\vu\in  L^\infty_{t,x}(]t_0-{\rho}^2,t_0[\times \bigcup_{i=1}^nB_{y_i,\mathfrak{r}_{y_i}})}$. Now, since $\bar B_{0,K} \subset\displaystyle{\bigcup_{i=1}^n B_{y_i,\mathfrak{r}_{y_i}}}$, we have
${\vu\in L^\infty_{t,x}(]t_0-{\rho}^2,t_0[\times\bar B_{0,K})}$. Therefore, from the previous information and  \eqref{BoundOutsideBk}, we can easily deduce that  ${\vu\in L^\infty(]t_0-\min\{\rho,r\}^2,t_0[\times \R)}$. Moreover, since $t_0\in ]\delta,T[$ was arbitrary, one has 
\begin{equation}\label{BoundAllSpace}
\vu\in L^{\infty}_{t,x}(]\delta,T[\times \R).\\[4mm]
\end{equation} 
\noindent Now, let us prove that  $\vu\in \mathcal{C}([\delta,T],L^2(\R))$. For this, we remark that it is known that it is sufficient to verify that  ${\vu\in L^2(]\delta,T[, H^1(\R))}$ and ${\partial_t \vu \in L^2(]\delta,T[,H^{-1}(\R)}$ to obtain this fact (see for instance \cite[Theorem 1.33]{Robinson}). Thus, since ${\vu\in L^\infty([0,+\infty[,L^2(\R))\cap L^2(]0,+\infty[, \dot H^1(\R))}$ by hypothesis, we have
\begin{align*}
\|\vu\|_{L^2(]\delta,T[, H^1(\R))}^2&=\|\vu\|_{L^2(]\delta,T[, L^2(\R))}^2+\|\vu\|_{L^2(]\delta,T[,\dot H^1(\R))}^2\\
&\le C\|\vu\|_{L^\infty (]\delta,T[, L^2(\R))}^2+\|\vu\|_{L^2(]\delta,T[,\dot H^1(\R))}^2<+\infty.
\end{align*} 
and hence we obtain that $\vu\in L^2(]\delta,T[, H^1(\R))$. Now, for proving ${\partial_t \vu \in L^2(]\delta,T[,H^{-1}(\R)}$, recall that $\vu$ satisfies the equation $\partial_t \vu = \Delta \vu -\mathbb{P}(\div(\vu\otimes \vu))+\frac{1}{2}\rot\vw$ where $\mathbb{P}(\cdot)$ is the Leray projector. Thus, since $\dot H^{-1}(\R)\subset H^{-1}(\R)$ and $\mathbb{P}$ is a bounded operator in $\dot H^{-1}(\R)$, we have
\begin{eqnarray*}
\|\partial_t \vu(t,\cdot)\|_{ H^{-1}(\R)}&\le &\|\partial_t \vu(t,\cdot)\|_{ \dot H^{-1}(\R)}\\
&\le&\|\Delta \vu(t,\cdot) \|_{ \dot H^{-1}(\R)}+\|\mathbb{P}(\div(\vu\otimes \vu))(t,\cdot)\|_{ \dot H^{-1}(\R)}+\|\frac{1}{2}\rot\vw(t,\cdot)\|_{ \dot H^{-1}(\R)}\\
&\le&\|\Delta \vu(t,\cdot) \|_{ \dot H^{-1}(\R)}+C\|\div(\vu\otimes \vu)(t,\cdot)\|_{ \dot H^{-1}(\R)}+\|\frac{1}{2}\rot\vw(t,\cdot)\|_{ \dot H^{-1}(\R)},
\end{eqnarray*}
and we can write
\begin{align*}
\|\partial_t \vu(t,\cdot)\|_{ H^{-1}(\R)}&\le \| \vu(t,\cdot) \|_{ \dot H^{1}(\R)}+C\|(\vu\otimes \vu)(t,\cdot)\|_{L^{2}(\R)}+C\|\vw(t,\cdot)\|_{L ^{2}(\R)}.
\end{align*}
By considering the $L^2$-norm in the time interval $]\delta,T[$ in the expression above, since ${\vu,\vw\in  L^\infty(]\delta,T[,L^{2}(\R))\cap L^2(]\delta,T[,\dot H^{1}(\R))}$ by hypothesis and since $\vu\in L^\infty(]\delta, T[,L^\infty(\R)) $ by \eqref{BoundAllSpace}, we have
\begin{eqnarray*}
\|\partial_t \vu\|_{ L^2(]\delta,T[,H^{-1}(\R))}\le \| \vu \|_{ L^2(]\delta,T[,\dot H^{1}(\R))}+\|\vu\otimes \vu\|_{ L^2(]\delta,T[,L^{2}(\R))}+\frac{1}{2}\|\vw\|_{ L^2(]\delta,T[,L^2(\R))}\qquad \qquad\\
\le  \| \vu \|_{ L^2(]\delta,T[,\dot H^{1}(\R))}+C\|\vu\|_{L^\infty(]\delta, T[,L^\infty(\R))}\|\vu\|_{ L^\infty(]\delta,T[, L^{2}(\R))}+C\|\vw\|_{ L^\infty(]\delta,T[,L^{2}(\R))}<+\infty.
\end{eqnarray*}
Thus, since we have proved that ${\vu\in L^2(]\delta,T[, H^1(\R))}$ and ${\partial_t \vu \in L^2(]\delta,T[,H^{-1}(\R))}$, it is possible to deduce that ${\vu\in \mathcal{C}([\delta,T],L^2(\R))}$.\\[4mm]

Having established that $\vu\in L^\infty_{t,x}(]\delta,T[\times \R)$ and  $\vu\in \mathcal{C}([\delta,T],L^2(\R))$, we will now prove that ${\vu\in \mathcal{C}(]\delta, T[, L^3(\R))}$ \textit{i.e.} we will study the continuity of the function
\begin{align}\label{ContinuityFunL3}
\begin{aligned}
]\delta,T[&\longrightarrow L^3(\R)\\
t&\longmapsto \vu(t,\cdot)
\end{aligned}
\end{align}
Remark that the previous function is well-defined since for any $t\in [\delta,T]$, we have $\vu(t,\cdot)\in L^3(\R)$ by  \eqref{EverytimeinL3Global}. Now, let $\varepsilon>0$ and  $t_1,t_2\in ]\delta,T[$. Since $\vu\in L^\infty(]\delta,T[,L^\infty(\R)) $ by \eqref{BoundAllSpace}, we have
\begin{align*}
\|\vu(t_1,\cdot)-\vu(t_2,\cdot)\|_{L^3(\R)}^3=\int_{\R}|\vu(t_1,x)-\vu(t_2,x)|^3dx&=\int_{\R}|\vu(t_1,x)-\vu(t_2,x)|^2|\vu(t_1,x)-\vu(t_2,x)|dx\\
&\le  2\|\vu\|_{L^\infty(]\delta,T[,L^\infty(\R))}\int_{\R}|\vu(t_1,x)-\vu(t_2,x)|^2dx.
\end{align*} 
On the other hand since $\vu\in \mathcal{C}([\delta,T],L^2(\R))$, there exists $\delta=\delta(\varepsilon)$ such that if $|t_1-t_2|\le \delta$, we have
\begin{equation*}
\int_{\R}|\vu(t_1,x)-\vu(t_2,x)|^2dx\le \frac{\varepsilon}{2\|\vu\|_{L^\infty(]\delta,T[,L^\infty(\R))}}.
\end{equation*}
Hence, there exists $\delta=\delta(\varepsilon)$ such that if $|t_1-t_2|\le \delta$, we have
\begin{equation}\label{ConclusionContinuitiLe}
\|\vu(t_1,\cdot)-\vu(t_2,\cdot)\|_{L^3(\R)}^3< \varepsilon.
\end{equation}
Thus, the function \eqref{ContinuityFunL3} is continuous and therefore we conclude that ${\vu\in \mathcal{C}(]\delta,T[,L^3(\R))}$, which finishes the proof of Theorem \ref{Theo_NormL3}. \hfill  $\blacksquare$
%%%%%%%%%%%%%%%%%%%%%%%%%%%%%%%%%%%%%%%%%%%%%%%%%%%
%%%%%%%%%%%%%%%%%%%%%%%%%%%%%%%%%%%%%%%%%%%%%%%%%%%
\section{Proof of Theorem \ref{Theo_BlowUp1}}\label{Secc_BlowUp1}
Let us recall the framework: we consider $(\vu,p,\vw)$ a Leray-type weak solution of the micropolar equations \eqref{MicropolarFluidsEquationsEqua1} and \eqref{MicropolarFluidsEquationsEqua2} and  let $0<\mathcal{T}\le +\infty$ be the maximal time such that  $ \vu\in\mathcal{C}(]0, \mathcal{T}[,L^3(\R))$. We thus want  to prove that if $\mathcal{T}<+\infty$, then $\underset{0<t<\mathcal{T}}{\sup}\|\vu(t,\cdot)\|_{L^3}=+\infty$.\\

\noindent To this end, we will need the following proposition.
%%%%%%%%%%%%%%%%%%%%%%%%%%%%%%%%%%%%%%%%%%%%%%%%%%%
\begin{Proposition}\label{Propo_ExtendeBeyondL3}
Let $(\vu,p,\vw)$ be a Leray-type weak solution of the micropolar equations \eqref{MicropolarFluidsEquationsEqua1} and \eqref{MicropolarFluidsEquationsEqua2} such that  for some $0<T_1<+\infty$ we have $\vu\in \mathcal{C}(]0, T_1[,L^3(\R))$. Then, the following assertions are equivalent:
\begin{enumerate}[leftmargin=*]
\item [1)]  For some $0<T_2<+\infty$ such that $T_1<T_2$, the velocity $\vu$ may be extended to the time interval $]0, T_2[$ such that we have the control $\vu\in \mathcal{C}(]0, T_2[,L^3(\R))$.
\item [2)] For every $x_0\in \R$, any point $(T_1,x_0)$ is partially regular in the sense of Definition \ref{Def_PartialRegularPoints}.
\end{enumerate} 
\end{Proposition}
%%%%%%%%%%%%%%%%%%%%%%%%%%%%%%%%%%%%%%%%%%%%%%%%%%%
\noindent\textbf{Proof.} Let us prove that $1)$ implies $2)$. Assume that for some $0<T_1<T_2<+\infty,$ the velocity $\vu$ may be extended to $]0, T_2[$  such that $\vu\in \mathcal{C}(]0, T_2[,L^3(\R))$. Notice that since $T_2<+\infty$ we have ${\vu\in L^\infty(]0,T_2[,L^3(\R))}$ and that we also have $\vu,\vw \in L^\infty(]0,+\infty[,L^2(\R))\cap L^2(]0,+\infty[,\dot{H}^1(\R))$ since it is a Leray-type weak solution. We can thus apply Theorem \ref{Theo_NormL3} and we obtain that any point $(t,x)\in ]0,T_2[\times \R$ is partially regular in the sense of Definition \ref{Def_PartialRegularPoints}. Since $0<T_1<T_2$ by hypothesis, it follows that for any $x_0\in \R$, the point $(T_1,x_0)$ is partially regular,  and this completes the proof of the first implication.\\

Now, we show  the converse \textit{i.e.,} we will prove that $2)$ implies $1)$ and we assume that for every $x_0\in \R$, the point $(T_1,x_0)$ is partially regular in the sense of Definition \ref{Def_PartialRegularPoints}. First, we remark that any point $(t,x)\in ]0, T_1]\times \R$ is also partially regular in the sense of Definition \ref{Def_PartialRegularPoints}: indeed since $\vu\in \mathcal{C}(]0, T_1[,L^3(\R))$ by hypothesis, and $T_1<+\infty$, we have $\vu\in L^\infty(]0,T_1[,L^3(\R))$. Since $\vu,\vw$ is a Leray-type solution we have $\vu,\vw\in L^\infty(]0,+\infty[,L^2(\R))\cap L^2(]0,+\infty[,\dot{H}^1(\R)) $. Thus, by Theorem \ref{Theo_NormL3}, any point $(t,x)\in]\delta,T_1[\times \R$ is partially regular. Recall that the case $ t=T_1$ follows from the assumption $2)$. \\

By using the same arguments as the ones used to deduce \eqref{BoundAllSpace}, we have that $ L^{\infty}_{t,x}(]0,T_1]\times \R)$ and similarly we can deduce that $\vu \in \mathcal{C}([0,T_1],L^2(\R))$. Therefore, following the same lines as in \eqref{ConclusionContinuitiLe}, we have
$\vu\in \mathcal{C}(]0, T_1],L^3(\R))$. It is worth noting that we are considering now the interval $]0,T_1]$.\\ 

\noindent To continue and in order to extend the solution beyond $t=T_1<+\infty$, we will use the following useful result.
%%%%%%%%%%%%%%%%%%%%%%%%%%%%%%%%%%%%%%%%%%%%%%%%%%%
\begin{Lemma}
Let $\vec f:]0,T[\times \R\longrightarrow \R$ be an exterior force with $\div(\vf)=0$, such that ${\vf \in L^p(]0,T[,L^q(\R))}$ with  $\frac{2}{p}+\frac{3}{q}=3$ and $\frac{3}{2}<q< 3$. Consider $\vv_0$ be a divergence-free initial data in $L^3(\R)$. Then, there exists $0<T_0<T$ and an unique solution $(\vv,h)$ of the forced Navier-Stokes equation 
\begin{equation*}
\begin{cases}
\partial \vv=\Delta \vv -(\vv\cdot\grad )\vv-\grad h+\vec f,\quad \div(\vv)=0,\\
\vv(0,\cdot)=\vv_0,
\end{cases}
\end{equation*} 
such that $\vv \in\mathcal{C}([0,T_0[, L^3(\R))\cap L^4(]0,T_0[, L^6(\R)).$ 
\end{Lemma}
%%%%%%%%%%%%%%%%%%%%%%%%%%%%%%%%%%%%%%%%%%%%%%%%%%%
\noindent For a proof of this result, we refer to Theorem 15.5 in \cite{PGLR1}. Remark that, since $\vw\in L^2([0,+\infty[,\dot H^1(\R))$ by hypothesis, we have for any $1\ll \kappa<+\infty $ that ${\rot \vw\in L^{2}(]0,\kappa T_1[,L^2(\R))}\subset  L^\frac{4}{3}(]0,\kappa T_1[,L^2(\R))$. Therefore, by considering the previous proposition with $\vu(T_1,\cdot)\in L^3(\R)$ as initial data and $\rot\vw$ as external force in $L^\frac{4}{3}_t(]T_1,\kappa T_1[,L^2(\R))$, there exists $0<T_1<T_2<\kappa T_1$ and a solution $(\vv,h)$ of the forced Navier-Stokes equations such that $\vv\in \mathcal{C}([T_1, T_2[,L^3(\R))$. Since, $(\vu,p)$ can be seen as a Leray-type weak solution of the same equation satisfied by $(\vv,h)$ (starting from the same initial data and the same external force), by a weak-strong uniqueness argument  we have that  $\vu=\vv\in \mathcal{C}([T_1, T_2[,L^3(\R))$ and hence the solution can be extended beyond $t=T_1$ such that $\vu\in \mathcal{C}(]0,T_2[,L^3(\R))$. This completes the proof of the second implication and this proves Proposition \ref{Propo_ExtendeBeyondL3}. \hfill $\blacksquare$\\

%%%%%%%%%%%%%%%%%%%%%%%%%%%%%%%%%%%%%%%%%%%%%%%%%%%
\noindent\textbf{Proof of Theorem \ref{Theo_BlowUp1}}. Let $0<\mathcal{T}\le +\infty$  be the maximal time such that  $ \vu\in\mathcal{C}(]0, \mathcal{T}[,L^3(\R))$. Recall that we want to prove that if $\mathcal{T}<+\infty$, then $\underset{0<t<\mathcal{T}}{\sup}\|\vu(t,\cdot)\|_{L^3}=+\infty$.
Assume the contrary, \textit{i.e.,} we have $\vu\in L^\infty(]0,\mathcal{T}[,L^3(\R))$. 
Since  $0<\mathcal{T}<+\infty$ is the maximal time such that $ \vu\in\mathcal{C}(]0, \mathcal{T}[,L^3(\R))$, from Proposition \ref{Propo_ExtendeBeyondL3}, there exists a point $x_0\in \R$ such that $(\mathcal{T},x_0)$ has to be a partial singular point in the sense of Definition \ref{Def_PartialRegularPoints}. 

On the other hand, for the same point $(\mathcal{T},x_0)$ since  $\vu\in L^\infty(]0,\mathcal{T}[,L^3(\R))$ by assumption, we can find  $0<R<\sqrt{\mathcal{T}}$ such that ${\vu\in L^\infty_tL^3_x(Q_R(\mathcal{T},x_0))}$. Moreover, since $(\vu,p,\vw)$ is a Leray-type solution we have $\vu,\vw\in L^\infty_tL^2_x(Q_R(\mathcal{T},x_0))\cap L^2_t\dot{H}^1_x(Q_R(\mathcal{T},x_0))$. 
Then, since we have moreover that $(\vu,p,\vw)$ is a partial suitable solution in $Q_R(\mathcal{T},x_0)$ and $\vu\in L^\infty_tL^3_x(Q_R(\mathcal{T},x_0))$, we can apply Theorem \ref{Theo_LocalPoints} and it follows that $(\mathcal{T},x_0)$ is actually a partial regular point in the sense of Definition \ref{Def_PartialRegularPoints}. This is a  contradiction since we have seen that $(\mathcal{T},x_0)$ is partially singular. Thus, the quantity $L^\infty_tL^3_x$ should explode and this  finishes the proof of the Theorem \ref{Theo_BlowUp1}. \hfill $\blacksquare$
%%%%%%%%%%%%%%%%%%%%%%%%%%%%%%%%%%%%%%%%%%%%%%%%%%%%%%%%%%%%%%%%%%%%%%%%%%%%%%%%%%%%%%%%%%%%%%%%%%%%%%
\mysection{The $L^3$-norm concentration effect}\label{Sec_Concentration}
In this section we will prove Theorem \ref{Theo_BlowUpConcentration}. More precisely, we will deduce  the concentration effect of the $L^3$-norm of the velocity $\vu$ around a partial singular point $(\mathcal{T},0)$ when $\vu\in \mathcal{C}([0,\mathcal{T}[,L^\infty(\R))$. Thus, if we assume that for $r_0>0$ such that ${0<\mathcal{T}-r_0^2}$ we have
\begin{equation*}
\sup_{x_0\in \R}\sup _{r\in ]0,r_0]}\sup _{t\in]\mathcal{T}-r^2,\mathcal{T}]} \frac{1}{r}\int_{B_{x_0,r}}|\vu(t,x)|^2dx =\mathfrak{M} <+\infty,
\end{equation*}
we will deduce that there exists $\varepsilon>0$, $\mathfrak{S}=\mathfrak{S}(\mathfrak{M})$ and $0<\delta<\mathcal{T}$ such that for all $t\in ]\mathcal{T}-\delta,\mathcal{T}[,$ we have
\begin{equation*}
\int_{B_{0,\sqrt{\frac{\mathcal{T}-t}{\mathfrak{S}}}}}|\vu(t,x)|^3dx\ge \varepsilon.
\end{equation*}
Before beginning the proof of Theorem \ref{Theo_BlowUpConcentration},  we need to introduce the following notion and some propositions.
%%%%%%%%%%%%%%%%%%%%%%%%%%%%%%%%%%%%%%%%%%%%%%%%%%%
\begin{Definition}[{\bf Partial local Leray solution}]\label{Def_partiaLocalLeray}
We will say that $\vu,\vw:]0, T[\times \R\longrightarrow \R$ and ${p: ]0, T[\times\R\longrightarrow \mathbb{R}}$ is a \emph{partial local Leray solution} of the micropolar fluids equations \eqref{MicropolarFluidsEquationsEqua1} and \eqref{MicropolarFluidsEquationsEqua2} with initial data ${\vu_0,\vw_0\in L^2(\R)}$ if:
\begin{enumerate}
\item[1)]we have $\displaystyle{\sup_{0<t<T}\sup_{x\in \R}\int_{B_{x,1}}|\vu|^2dy+\sup_{x\in \R}\int_{0}^T\int_{B_{x,1}}|\grad\otimes\vu|^2dydt<+\infty}$,
\item[2)]the triplet $(\vu,p,\vw)$ is a partial suitable solution on $]0,T[\times \R$ in the sense of Definition \ref{Def_PartialSuitable},
\item[3)]for every compact subset $K$ of $\R$, we have
\begin{equation}\label{LocalLerayInitialData}
{\underset{t\longrightarrow 0^+}{\lim}\displaystyle{\int_K}|\vu(t,y)-\vu_0(y)|^2dy=0},
\end{equation} 
\item[4)]for any $R>0$ we have $\displaystyle{\lim_{|x_0|\to +\infty}\int_0^{\min\{R^2,T\}} \int_{B{x_0,R}}|\vu|^2dyds=0}$.
\end{enumerate}
\end{Definition}
%%%%%%%%%%%%%%%%%%%%%%%%%%%%%%%%%%%%%%%%%%%%%%%%%%%
\noindent This notion of \emph{local Leray} solution is borrowed from the theory of the Navier-Stokes problem. See in particular \cite[Definition 14.1]{PGLR1}, and \cite[Appendix B]{Seregin14} where the global setting considered there (the $L^2_{uloc}$ space) is slightly more general than the one considered here. However, for our purposes the $L^2$ setting stated above is enough. Note again that in the previous definition, we are not imposing any particular hypothesis over the variable $\vw$, leading us to the previous ``\emph{partial}'' notion of local Leray solutions.
\begin{Remark}\label{Rem_PartialLocalEqualLocal}
It is worth noting that if $(\vu,p,\vw)$ is a partial local Leray solution of the micropolar fluids equations  \eqref{MicropolarFluidsEquationsEqua1} and \eqref{MicropolarFluidsEquationsEqua2} in the sense of Definition \ref{Def_partiaLocalLeray}, we may say  that $(\vu,p)$ is a local Leray solution of the forced Navier-Stokes equations in the sense of the Definition 14.1 of the book \cite{PGLR1} where the quantity $\frac{1}2\rot\vw\in L^2_{t,x}$ can be considered as an external force, \textit{i.e.,}
$\partial_t \vu = \Delta \vu -(\vu \cdot \vn)\vu-\vn p +\frac{1}{2}\rot\vw$, $\div(\vu)=0$.
\end{Remark}
\noindent We present now some lemmas to highlight some properties of the partial local Leray solutions introduced above. First, we remark that the pressure can be studied in the same way as in the classical Navier-Stokes equations since the variable  $\vw$ is not present in the equation \eqref{Eq_Pression}. Thus, we have the following local decomposition
%%%%%%%%%%%%%%%%%%%%%%%%%%%%%%%%%%%%%%%%%%%%%%%%%%%
\begin{Lemma}\label{Propo_LocalDescompositionPressure}
Let $(\vu,p,\vw)$ be a partial local Leray solution in the sense of Definition \ref{Def_partiaLocalLeray} of the micropolar fluids equations \eqref{MicropolarFluidsEquationsEqua1} and \eqref{MicropolarFluidsEquationsEqua2}. Then, the pressure $p$ can be decomposed as follows: for all $x_0\in \R$ and $r>0$, there exists $\mathfrak{h}\in L^\frac{3}{2}(]0,T[)$ such that 
\begin{align}
p(x,t)-\mathfrak{h}(t)&=\frac{1}{(-\Delta)}(\div(\div(\mathds{1}_{B_{x_0,3r}} \vu\otimes \vu)))+\int_{|y-x_0|>3r}(\mathbb{K}(x-y)-\mathbb{K}(-y))\vu\otimes\vu(t,y))dy\notag\\
&=\mathfrak{p}_1+\mathfrak{p}_2,\label{LocalDescompopressu}
\end{align} 
where $\mathbb{K}$ is the kernel of the singular integral operator $\frac{1}{(-\Delta)}(\div(\div))$.
\end{Lemma}
%%%%%%%%%%%%%%%%%%%%%%%%%%%%%%%%%%%%%%%%%%%%%%%%%%%
\noindent For a proof of this lemma in the setting of the Navier-Stokes equations, we refer to \cite[Lemma 3.4]{KanMiuTsai20}, see also \cite[Theorem 4]{BarPran19} and the article \cite{Seregin14}. \\

Now, we observe that since $\vu$ satisfies the condition \eqref{LocalLerayInitialData}, it is possible to rewrite the local  energy inequality \eqref{LocalEnergyInequality} in terms of the initial data as follows.
%%%%%%%%%%%%%%%%%%%%%%%%%%%%%%%%%%%%%%%%%%%%%%%%%%%
\begin{Lemma}
Let $(\vu,p,\vw)$ be a partial local Leray solution in the sense of Definition \ref{Def_partiaLocalLeray} on $]0,T[\times \R$. For all $\Phi \in\mathcal{C}_0^{\infty}(\R)$ and for  all $t\in ]0,T[$ we have
\begin{align}
\int_{\R} |\vu|^2\Phi dy&+ \int_{0}^t \int _{\R}|\grad \otimes \vu|^2\Phi dy ds
\le \int_{\R}|\vu_0|^2 \Phi dy+ \int_{0}^t \int_{\R}|\vu|^2\Delta \Phi dyds+
\int_{0}^t \int_{\R}\rot \vw \cdot \vu\Phi dyds \nonumber\\
&+\int_{0}^t \int _{\R} (|\vu|^2+2[\mathfrak{p}_1+\mathfrak{p}_2])\vu\cdot \grad\Phi dyds.\label{LocalEnergyInequalityInitalData}
\end{align}
\end{Lemma}
%%%%%%%%%%%%%%%%%%%%%%%%%%%%%%%%%%%%%%%%%%%%%%%%%%%
\noindent It is worth noting that we are able to take here test functions constant in time in the local energy inequality. We refer to \cite[Remark 1.2]{Mayea} for a proof of this result (the term $\rot \vw$ is considered here as an ``external force'').\\

Having announced these previous results, we  present now the main tool to prove the concentration effect of the $L^3-$norm stated in Theorem \ref{Theo_BlowUpConcentration}.  
%%%%%%%%%%%%%%%%%%%%%%%%%%%%%%%%%%%%%%%%%%%%%%%%%%%
\begin{Lemma}\label{Lem_localSmoothing2}
Let $(\vu,p,\vw)$ be a partial local Leray solution of the micropolar fluids equations 
on $]0,1[\times \R$  associated to the initial data $\vu_0,\vw_0\in L^2(\R)$ in the sense of Definition \ref{Def_partiaLocalLeray}, such that there exists $M>0$ with
\begin{equation}\label{Hypho_EstimatesAPriori}
\sup_{0<t<1}\sup_{x\in \R}\int_{B_{x,1}}|\vu(t,x)|^2dy+\sup_{x\in \R}\int_{0}^{1}\int_{B_{x,1}}|\grad\otimes\vu|^2dyds\le M.
\end{equation}
Assume moreover that for some $0<R<\frac{1}{12}$ and $S>0$ we have
\begin{equation}\label{Hypho_EstimatesOmega}
\|\vw\|_{L^\infty_tL^2_x(Q_1(1,0))}^2<CR.
\end{equation}
and
\begin{equation}\label{Hypo_LocalSmoothing2}
{\displaystyle{\sup_{R<\mathfrak{r}\le 1}\frac{1}{\mathfrak{r}}\int_{B_{0,\mathfrak{r}}}|\vu_0|^2dy}}=S<+\infty.
\end{equation}
Then,  there exists
$T_1=T_1(M,S)<1$  and a universal constant $\mathfrak{\mathfrak{c}}>1$ such that for all $\mathfrak{r}>0$ with $R\le \mathfrak{r}\le\frac{1}{3}$ and for $t>0$ such that  $0<t\le T^\ast=\min\{ T_1,\mathfrak{c}\lambda \mathfrak{r}^2\}$ with $\lambda=\frac{1}{1+S^2}$, we have the control
\begin{equation}\label{ConclusionLocalSmothiing}
E_\mathfrak{r}(t)=\sup_{0<s<t}\frac{1}{\mathfrak{r}}\int _{B_{0,\mathfrak{r}}}|\vu|^2dy+\frac{1}{\mathfrak{r}}\int_0^t \int_{B_{0,\mathfrak{r}}}|\grad\otimes\vu|^2dyds
+\frac{1}{\mathfrak{r}^2}\int_0^t\int_{B_{0,\mathfrak{r}}}|p-\mathfrak{h}|^{\frac{3}{2}}dyds
<C S.
\end{equation}
\end{Lemma}
%%%%%%%%%%%%%%%%%%%%%%%%%%%%%%%%%%%%%%%%%%%%%%%%%%%
\noindent This result was originally established within the framework of the classical Navier-Stokes equations without force in \cite[Theorem 3.1]{KanMiuTsai19}. In our case, since we are dealing with the micropolar fluids equations \eqref{MicropolarFluidsEquationsEqua1} and \eqref{MicropolarFluidsEquationsEqua2}, we need to take into account  the term $\rot \vw$ in the  equation related to the evolution of $\vu$ and this lead us to the condition (\ref{Hypho_EstimatesOmega}). 
%%%%%%%%%%%%%%%%%%%%%%%%%%%%%%%%%%%%%%%%%%%%%%%%%%%
\begin{Remark}\label{Rem_SmallNessOmega} This result will be applied later on to a suitable re-scaled system. In particular, the smallness hypothesis (\ref{Hypho_EstimatesOmega}) will be a consequence of this re-scaling. See formula (\ref{InforotvW}) below for more details. 
\end{Remark}

\noindent\textbf {Proof of Lemma \ref{Lem_localSmoothing2}.} Let $r>0$ such that $R\le r\le\frac{1}{3}$ be a fixed radius and $0<t<T^*<1$ for some $T^*$ to be fixed later. First, notice that by the local decomposition of the pressure given in the Lemma \ref{Propo_LocalDescompositionPressure}, we have $p-\mathfrak{h}=\mathfrak{p}_1+\mathfrak{p}_2$ where $\mathfrak{p}_1$ and $\mathfrak{p}_2$ are given in \eqref{LocalDescompopressu}. Then $E_r(t)$ can be written in the following manner
\begin{align*}
E_r(t)
&=\sup_{0<s<t}\frac{1}{r}\int _{B_{0,r}}|\vu|^2dy+\frac{1}{r}\int_0^t \int_{B_{0,r}}|\grad\otimes\vu|^2dyds
+\frac{1}{r^2}\int_0^t\int_{B_{0,r}}|\mathfrak{p}_1+\mathfrak{p}_2|^{\frac{3}{2}}dyds.\notag
\end{align*}
Moreover it is easy to see that
\begin{align}
E_r(t)
&\le\left(\sup_{0<s<t}\frac{1}{r}\int _{B_{0,r}}|\vu|^2dy+\frac{1}{r}\int_0^t \int_{B_{0,r}}|\grad\otimes\vu|^2dyds\right)+\frac{C}{r^2}\int_{0}^t \int_{B_{0,2r}}|\mathfrak{p}_1|^\frac{3}{2}+|\mathfrak{p}_2|^\frac{3}{2}dyds\label{FirstInequalityEr}.
\end{align}
In order to control the expression above in terms of the initial data \eqref{Hypo_LocalSmoothing2}, we will study more in detail the  terms inside the parentheses above. For this, we may use the local energy inequality \eqref{LocalEnergyInequalityInitalData} by considering a well-chosen test function. Indeed, let $\phi\in \mathcal{C}_0^{\infty}(\R)$ be a positive function such that $ \phi = 1$ in $B_{0,r}$, $\supp(\phi)\subset B_{0,2r}$ and for all $k\in \mathbb{N}$ and all multi-index $\alpha\in \mathbb{N}^3$, such that $|\alpha|\le k$ we have $\|D^\alpha\phi \|_{L^\infty}\le C_k r^{-k}$. Now, with this auxiliary function in the local energy inequality \eqref{LocalEnergyInequalityInitalData}, one has
\begin{align*}
\int_{\R} |\vu|^2\phi dy&+ \int_{0}^t \int _{\R}|\grad \otimes \vu|^2\phi dy ds
\le \int_{\R}|\vu_0|^2 \phi dy+ \int_{0}^t \int_{\R}|\vu|^2\Delta \phi dyds+
\int_{0}^t \int_{\R}\rot \vw \cdot (\vu\phi )dyds \nonumber\\
&+\int_{0}^t \int _{\R} (|\vu|^2+2[\mathfrak{p}_1+\mathfrak{p}_2])\vu\cdot \grad\phi dyds.
\end{align*}
Since by integration by parts we have
$$\int_{0}^t \int_{\R}\rot \vw \cdot (\vu\phi )dyds =\int_{0}^t \int_{\R}\rot \vu \cdot (\vw\phi )dyds +\int_{0}^t \int_{\R} (\vw\cdot \vu )\wedge \grad\phi dyds, $$ 
we obtain that the terms inside  parentheses in the left-hand side of the expression \eqref{FirstInequalityEr} can be bounded as follows:
\begin{align*}
\sup_{0<s<t}\int_{\R} |\vu|^2\phi dy&+ \int_{0}^t \int _{\R}|\grad \otimes \vu|^2\phi dy ds
\le \int_{\R}|\vu_0|^2 \phi dy+ \int_{0}^t \int_{\R}|\vu|^2\Delta \phi dyds\int_{0}^t \int_{\R}\rot \vu \cdot (\vw\phi )dyds  \nonumber\\
&+\int_{0}^t \int_{\R} (\vw\cdot \vu )\wedge \grad\phi dyds+\int_{0}^t \int _{\R} |\vu|^2\vu\cdot \grad\phi dyds+\int_{0}^t \int _{\R}2[\mathfrak{p}_1+\mathfrak{p}_2]\vu\cdot \grad\phi dyds.
\end{align*}
Now, by the properties of the test function $\phi$ we obtain
\begin{eqnarray*}
\sup_{0<s<t}\int_{B_{0,r}} |\vu|^2 dy+ \int_{0}^t \int _{B_{0,r}}|\grad \otimes \vu|^2 dy ds&\le& 
C\int_{B_{0,2r}}|\vu_0|^2 dy+ \frac{C}{r^2}\int_{0}^t \int_{B_{0,2r}}|\vu|^2 dyds\\
&&+C\int_{0}^t \int_{B_{0,2r}}|\rot \vu ||\vw| dyds+\frac{C}{r}\int_{0}^t \int_{\R} |\vw||\vu| dyds\\
&&+\frac{C}{r}\int_{0}^t \int _{\R} |\vu|^3 dyds+\frac{C}{r}\int_{0}^t \int _{\R}|\mathfrak{p}_1||\vu|+|\mathfrak{p}_2||\vu|dyds.
\end{eqnarray*}
Note that from the H\"older inequality $(1=\frac{2}{3}+\frac{1}{3})$ and the Young inequality, one has for the last term above the estimate $\displaystyle{\int_{0}^t \int _{\R}|\mathfrak{p}_1||\vu|+|\mathfrak{p}_2||\vu|dyds\le C \int_{0}^t \int _{\R}|\mathfrak{p}_1|^\frac{3}{2}+|\mathfrak{p}_2|^\frac{3}{2}+|\vu|^3dyds}$, hence, by applying the previous estimate in the inequality above, we have
\begin{eqnarray*}
\sup_{0<s<t}\int_{B_{0,r}} |\vu|^2 dy+ \int_{0}^t \int _{B_{0,r}}|\grad \otimes \vu|^2 dy ds\le 
C\int_{B_{0,2r}}|\vu_0|^2 dy+ \frac{C}{r^2}\int_{0}^t \int_{B_{0,2r}}|\vu|^2 dyds\qquad \\
+C\int_{0}^t \int_{B_{0,2r}}|\rot \vu ||\vw| dyds+\frac{C}{r}\int_{0}^t \int_{\R} |\vw||\vu| dyds\\
+\frac{C}{r}\int_{0}^t \int _{\R} |\vu|^3 dyds+{\frac{C}{r^2}\int_{0}^t \int_{B_{0,2r}}|\mathfrak{p}_1|^\frac{3}{2}dyds}+{\frac{C}{r^2}\int_{0}^t \int_{B_{0,2r}}|\mathfrak{p}_2|^\frac{3}{2}dyds}.
\end{eqnarray*}
Now, by multiplying  by $\frac{1}{r}$ in the expression above, we obtain that the following bound for the terms in parentheses in \eqref{FirstInequalityEr}
\begin{align}
&\sup_{0<s<t}\frac{1}{r}\int_{B_{0,r}} |\vu|^2dy
+ \frac{1}{r}\int_{0}^t \int _{B_{0,r}}|\grad \otimes \vu|^2 dyds
\le \frac{C}{r} \int_{B_{0,2r}}|\vu_0|^2 dy+{ \frac{C}{r^3}\int_{0}^t \int_{B_{0,2r}}|\vu|^2 dyds}\notag\\\nonumber 
&+{\frac{C}{r} \int_{0}^t \int_{B_{0,2r}}|\rot \vu| |\vw| dyds}
+{\frac{C}{r^2}\int_{0}^t \int_{B_{0,2r}}| \vw| |\vu| dyds}+ {\frac{C}{r^2}\int_0^t\int_{B_{0,2r}}|\vu|^3dy ds}+{\frac{C}{r^2}\int_{0}^t \int_{B_{0,2r}}|\mathfrak{p}_1|^\frac{3}{2}dyds}\\
&+{\frac{C}{r^2}\int_{0}^t \int_{B_{0,2r}}|\mathfrak{p}_2|^\frac{3}{2}dyds}.\nonumber
\end{align}
Thus, by replacing the previous estimate in \eqref{FirstInequalityEr}, one has
\begin{align}
E_r(t)&\le \biggl(\frac{C}{r} \int_{B_{0,2r}}|\vu_0|^2 dy+{ \frac{C}{r^3}\int_{0}^t \int_{B_{0,2r}}|\vu|^2 dyds}+{\frac{C}{r} \int_{0}^t \int_{B_{0,2r}}|\rot \vu| |\vw| dyds}\notag\\
&+{\frac{C}{r^2}\int_{0}^t \int_{B_{0,2r}}| \vw| |\vu| dyds}+ {\frac{C}{r^2}\int_0^t\int_{B_{0,2r}}|\vu|^3dy ds}+{\frac{C}{r^2}\int_{0}^t \int_{B_{0,2r}}|\mathfrak{p}_1|^\frac{3}{2}dyds}\notag\\
&+{\frac{C}{r^2}\int_{0}^t \int_{B_{0,2r}}|\mathfrak{p}_2|^\frac{3}{2}dyds}\biggl)+\frac{C}{r^2}\int_{0}^t \int_{B_{0,2r}}|\mathfrak{p}_1|^\frac{3}{2}+|\mathfrak{p}_2|^\frac{3}{2}dyds\notag\\[3mm]
&\le\frac{C}{r} \int_{B_{0,2r}}|\vu_0|^2 dy+\underbrace{ \frac{C}{r^3}\int_{0}^t \int_{B_{0,2r}}|\vu|^2 dyds}_{I_1}+\underbrace{\frac{C}{r} \int_{0}^t \int_{B_{0,2r}}|\rot \vu| |\vw| dyds}_{I_2}+\underbrace{\frac{C}{r^2}\int_{0}^t \int_{B_{0,2r}}| \vw| |\vu| dyds}_{I_3}\notag\\
&+\underbrace{\frac{C}{r^2}\int_0^t\int_{B_{0,2r}}|\vu|^3dy ds}_{I_4}+\underbrace{\frac{C}{r^2}\int_{0}^t \int_{B_{0,2r}}|\mathfrak{p}_1|^\frac{3}{2}dyds}_{I_5}+\underbrace{\frac{C}{r^2}\int_{0}^t \int_{B_{0,2r}}|\mathfrak{p}_2|^\frac{3}{2}dyds}_{I_6}.\label{Teo7_termsToBound}
\end{align}
Now, in order to obtain the wished estimate (\ref{ConclusionLocalSmothiing}) and for $R\leq \mathfrak{r}\leq \frac{1}{3}$, we will study the following expression  
$$\mathcal{E}_\mathfrak{r}(t)=\displaystyle{\sup_{\mathfrak{r}\le r\le \frac{1}{3}}E_r(t)}.$$ 
Remark that we have by construction $E_\mathfrak{r}(t)\leq \mathcal{E}_\mathfrak{r}(t)$, and to study the term $\mathcal{E}_\mathfrak{r}(t)$ we split the previous supremum into two parts:
\begin{equation}\label{SplitER}
\mathcal{E}_\mathfrak{r}(t)\le\sup_{\mathfrak{r}\le r\le \frac1{12}}E_r(t)+\sup_{\frac1{12}<r\le \frac13}E_r(t).
\end{equation}
In the following, we study each one of the terms above separately. 
\begin{itemize}
\item [$\bullet$] \underline{Assume $\mathfrak{r}\le r\le\frac{1}{12}$:} Note that from \eqref{Teo7_termsToBound}, we have 
\begin{equation}\label{Teo7_termsToBound1}
\sup_{\mathfrak{r}\le r\le \frac1{12}}E_r(t)\le \sup_{\mathfrak{r}\le r\le \frac1{12}} \frac{C}{r} \int_{B_{0,2r}}|\vu_0|^2 dy+\sup_{\mathfrak{r}\le r\le \frac1{12}}\sum_{j=1}^{6}I_j
\end{equation} 
For the term $I_1$ of \eqref{Teo7_termsToBound1}, by the definition of $E_r(t)$ and  since $\mathfrak{r}<2r<\frac{1}{6}<\frac{1}{3}$,  we have
\begin{align}
\sup_{\mathfrak{r}\le r\le \frac1{12}}I_{1}&=\sup_{\mathfrak{r}\le r\le \frac1{12}}\frac{C}{r^3}\int_{0}^t \int_{B_{0,2r}}|\vu|^2 dyds=\sup_{\mathfrak{r}\le r\le \frac1{12}}\frac{2C}{r^2}\int_{0}^t \frac{1}{2r}\int_{B_{0,2r}}|\vu|^2 dyds\notag\\
&\le \sup_{\mathfrak{r}\le r\le \frac1{12}} \frac{2C}{r^2}\int_{0}^t E_{2r}(s) ds \le \frac{C}{R^2}\int_0^t\mathcal{E}_\mathfrak{r}(s)ds.\label{Case1I_1}
\end{align}
For the term $I_2$ in \eqref{Teo7_termsToBound1}, by the Cauchy-Schwarz and Young inequalities we can write
\begin{eqnarray*}
\sup_{\mathfrak{r}\le r\le \frac1{12}}I_{2}&=&\sup_{\mathfrak{r}\le r\le \frac1{12}}\frac{C}{r} \int_{0}^t \int_{B_{0,2r}}|\rot \vu| |\vw| dyds\\
&\le &\sup_{\mathfrak{r}\le r\le \frac1{12}} \frac{C}{r}\int_0^t \left(\int_{B_{0,2r}} |\rot \vu| ^2dy\right)^{\frac{1}{2}}
\left(\int_{B_{0,2r}} |\vw| ^2dy\right)^{\frac{1}{2}}ds\\
&\le &\sup_{\mathfrak{r}\le r\le \frac1{12}}\frac{1}{16r} \int_0^t\int_{B_{0,2r}}|\rot\vu| ^2dyds+\sup_{\mathfrak{r}\le r\le \frac1{12}}\frac{C}{r}\int_0^t\int_{B_{0,2r}} |\vw| ^2dyds.
\end{eqnarray*}
Since $|\rot\vu|^2\le 2|\grad\otimes\vu|^2$ and $B_{0,2r}\subset B_{0,1}$ due to the fact that $2r<1$ ,  we have
\begin{align*}
\sup_{\mathfrak{r}\le r\le \frac1{12}}I_{2}&\le \sup_{\mathfrak{r}\le r\le \frac1{12}}  \frac{1}{8r} \int_0^t\int_{B_{0,2r}}|\grad\otimes\vu| ^2dyds+\sup_{\mathfrak{r}\le r\le \frac1{12}}\frac{C}{r}\int_0^t\int_{B_{0,2r}} |\vw| ^2dyds\\
&\le\sup_{\mathfrak{r}\le r\le \frac1{12}}\frac{1}{8r} \int_0^t\int_{B_{0,2r}}|\grad\otimes\vu| ^2dyds+\sup_{\mathfrak{r}\le r\le \frac1{12}}\frac{C}{r}\int_0^t\|\vw(s,\cdot)\|_{L^2(B_{0,1})} ^2ds.
\end{align*}
Now, by using the definition of $E_r(t)$,  we obtain
\begin{align*}
\sup_{\mathfrak{r}\le r\le \frac1{12}}I_2 &\le \sup_{\mathfrak{r}\le r\le \frac1{12}}\frac{1}{4} E_{2r}(t)+\sup_{\mathfrak{r}\le r\le \frac1{12}}\frac{Ct}{r}\|\vw\|_{L^\infty_t L^2_x(Q_1(1,0))}^2.
\end{align*}
Then, since  $\frac{1}{r}\le\frac{1}{\mathfrak{r}}$ and by the definition of $\mathcal{E}_\mathfrak{r}$, we have
\begin{equation*}
\sup_{\mathfrak{r}\le r\le \frac1{12}}I_2\le\frac{1}{4} \mathcal{E}_\mathfrak{r}+\frac{Ct}{\mathfrak{r}} \|\vw\|_{L_t^\infty L^2_x(Q_1(1,0))}^2.
\end{equation*}
Furthermore, since  $\|\vw\|_{L_t^\infty L^2_x(Q_1(1,0))}^2<CR\leq C\mathfrak{r}$ by the hypothesis \eqref{Hypho_EstimatesOmega} (and since we are assuming that $R\leq \mathfrak{r}\leq \frac{1}{3}$), we finally obtain
\begin{equation}
\sup_{\mathfrak{r}\le r\le \frac1{12}}I_2\le\frac{1}{4} \mathcal{E}_\mathfrak{r}+Ct.\label{Case1I_2}
\end{equation}
For the term $I_3$ in \eqref{Teo7_termsToBound1}, notice that it can be rewritten as follows
\begin{align*}
\sup_{\mathfrak{r}\le r\le \frac1{12}}I_{3}=\sup_{\mathfrak{r}\le r\le \frac1{12}}\frac{C}{r^2}\int_{0}^t \int_{B_{0,2r}}| \vw| |\vu| dyds=\sup_{\mathfrak{r}\le r\le \frac1{12}}C\int_{0}^t \int_{B_{0,2r}}\frac{1}{r^\frac12}| \vw| \frac{1}{r^\frac32}|\vu| dyds.
\end{align*}
Hence, by the Cauchy-Schwarz and Young inequalities, we have 
\begin{align*}
\sup_{\mathfrak{r}\le r\le \frac1{12}}I_{3}&\le\sup_{\mathfrak{r}\le r\le \frac1{12}} \frac{C}{r}\int_0^t \int_{B_{0,2r}} |\vw| ^2dyds
+ \sup_{\mathfrak{r}\le r\le \frac1{12}}\frac{1}{2r^3}\int_0^t \int_{B_{0,2r}} |\vu| ^2dyds\notag\\
&\le \sup_{\mathfrak{r}\le r\le \frac1{12}}\frac{C}{r}\int_0^t \int_{B_{0,2r}} |\vw| ^2dyds
+ \sup_{\mathfrak{r}\le r\le \frac1{12}}\frac{1}{r^2}\int_0^t E_{2r}(s)ds.\notag
\end{align*}
Again, since $B_{0,2r}\subset B_{0,1}$  we write
\begin{align*}
\sup_{\mathfrak{r}\le r\le \frac1{12}}I_{3}& \le \sup_{R\le r\le \frac1{12}} \frac{C}{r}\int_0^t \int_{B_{0,1}} |\vw| ^2dyds
+ \sup_{\mathfrak{r}\le r\le \frac1{12}}\frac{1}{r^2}\int_0^t E_{2r}(s)ds\notag\\
& \le \sup_{\mathfrak{r}\le r\le \frac1{12}} \frac{Ct}{r}\|\vw\|^2_{L_t^\infty L^2_x(Q_1(1,0))}
+ \sup_{\mathfrak{r}\le r\le \frac1{12}} \frac{1}{r^2}\int_0^t E_{2r}(s)ds.
\end{align*}
Therefore,  since $\|\vw\|^2_{L_t^\infty L^2_x(Q_1(1,0))}<CR\leq C\mathfrak{r}$ by \eqref{Hypho_EstimatesOmega} and by the definition of $\mathcal{E}_\mathfrak{r}$, we have
\begin{equation}
\sup_{\mathfrak{r}\le r\le \frac1{12}}I_{3}\le Ct+ \frac{1}{\mathfrak{r}^2}\int_0^t \mathcal{E}_{\mathfrak{r}}(s)ds.\label{Case1I_3}
\end{equation}
For the  term $I_4$ in \eqref{Teo7_termsToBound1}, first notice that  by the classical Gagliardo-Niremberg inequality  with $\frac{1}{3}=\theta(\frac{1}{2}-\frac{1}{3})+\frac{1-\theta}{2}$ (see \cite{Brezis}),  we obtain
$$\|\vu(t,\cdot)\|_{L^3(B_{0,2r})}\le \|\grad \otimes\vu(t,\cdot)\|_{L^2(B_{0,2r})}^{\frac{1}{2}}\|\vu(t,\cdot)\|^{\frac{1}{2}}_{L^2(B_{0,2r})}+\|\vu(t,\cdot)\|_{L^2(B_{0,2r})}.$$
Moreover, by the inequality above, one has
\begin{eqnarray}
\frac{C}{r^2}\int_0^t\int_{B_{0,2r}}|\vu|^3dy ds=\frac{C}{r^2} \int_{0}^t\|\vu(s,\cdot)\|_{L^3(B_{0,2r})}^3ds\hspace{4cm}\notag\\
\le \frac{C}{r^2} \int_{0}^t\|\grad \otimes\vu(s,\cdot)\|_{L^2(B_{0,2r})}^{\frac{3}{2}}\|\vu(s,\cdot)\|^{\frac{3}{2}}_{L^2(B_{0,2r})}+\|\vu(s,\cdot)\|_{L^2(B_{0,2r})}^{3}ds.\label{ControlNormL3}
\end{eqnarray}
Hence by the previous estimate and by the Young inequality for the sum (with $1=\frac{3}{4}+\frac{1}{4}$)  we obtain 
\begin{align*}
&\sup_{\mathfrak{r}\le r\le \frac1{12}}I_4=\sup_{\mathfrak{r}\le r\le \frac1{12}}\frac{C}{r^2}\int_0^t\int_{B_{0,2r}}|\vu|^3dy ds\\
&\le \sup_{\mathfrak{r}\le r\le \frac1{12}}\frac{C}{r^2} \int_{0}^t\|\grad \otimes\vu(s,\cdot)\|_{L^2(B_{0,2r})}^{\frac{3}{2}}\|\vu(s,\cdot)\|^{\frac{3}{2}}_{L^2(B_{0,2r})}+\|\vu(s,\cdot)\|_{L^2(B_{0,2r})}^{3}ds\nonumber\\
&\leq\sup_{\mathfrak{r}\le r\le \frac1{12}}\int_{0}^t \frac{1}{r^\frac{3}{4}}\|\grad \otimes\vu(s,\cdot)\|_{L^2(B_{0,2r})}^{\frac{3}{2}}\frac{C}{r^\frac{5}{4}}\|\vu(s,\cdot)\|^{\frac{3}{2}}_{L^2(B_{0,2r})}ds+\sup_{\mathfrak{r}\le r\le \frac1{12}}\frac{C}{r^2}\int_{0}^t \|\vu(s,\cdot)\|_{L^2(B_{0,2r})}^{3}ds\\
&\le\sup_{\mathfrak{r}\le r\le \frac1{12}}\frac{1}{8r}\int_{0}^t\|\grad \otimes\vu(s,\cdot)\|_{L^2(B_{0,2r})}^2ds+ \frac{C}{r^5}\int_0^t \|\vu(s,\cdot)\|_{L^2(B_{0,2r})}^6ds +\frac{C}{r^2}\int_0^t \|\vu(s,\cdot)\|_{L^2(B_{0,2r})}^{3}ds,
\end{align*}
and we can write
\begin{align*}
\sup_{\mathfrak{r}\le r\le \frac1{12}}I_4&\le\sup_{\mathfrak{r}\le r\le \frac1{12}}\frac{1}{8r}\int_{0}^t\|\grad \otimes\vu(s,\cdot)\|_{L^2(B_{0,2r})}^2ds+\sup_{\mathfrak{r}\le r\le \frac1{12}} \frac{8C}{r^2}\int_0^t \left(\frac{1}{2r}\int_{B_{0,2r}}|\vu|^2dy\right)^3ds \\&\quad+\sup_{\mathfrak{r}\le r\le \frac1{12}}\frac{2^\frac{3}{2}C}{r^\frac{1}{2}}\int_0^t \left(\frac{1}{2r}\int_{B_{0,2r}}|\vu|^2dy\right)^\frac{3}{2}ds.
\end{align*}
Now, by the definition of $E_r(t)$, we have
\begin{align*}
\sup_{\mathfrak{r}\le r\le \frac1{12}}I_4 &\le \sup_{\mathfrak{r}\le r\le \frac1{12}} \frac{1}{4} E_{2r}(t)+\sup_{\mathfrak{r}\le r\le \frac1{12}}\frac{C}{r^2}\int_0^t E^3_{2r}(s)ds +\sup_{\mathfrak{r}\le r\le \frac1{12}}\frac{C}{r^\frac{1}{2}}\int_0^tE^\frac{3}{2}_{2r}(s)ds,
\end{align*}
but since $r<1$, it follows that $\frac{1}{r^\frac{1}{2}}<\frac1{r^2}$ and we obtain
\begin{align*}
\sup_{\mathfrak{r}\le r\le \frac1{12}}I_4&\le \sup_{\mathfrak{r}\le r\le \frac1{12}}\frac{1}{4} E_{2r}(t)+\sup_{\mathfrak{r}\le r\le \frac1{12}}\frac{C}{r^2}\int_0^t E^3_{2r}(s)ds+\sup_{\mathfrak{r}\le r\le \frac1{12}}\frac{C}{r^2}\int_0^tE^\frac{3}{2}_{2r}(s)ds,
\end{align*}
using again the fact that $\frac{1}{r}\le\frac{1}{\mathfrak{r}}$ and by the definition of $\mathcal{E}_\mathfrak{r}(t)$, we finally have
\begin{align}
\sup_{\mathfrak{r}\le r\le \frac1{12}}I_4 \le\frac{1}{4} \mathcal{E}_\mathfrak{r}(t)+\frac{C}{\mathfrak{r}^2}\int_0^t\mathcal{E}^3_\mathfrak{r}(s)ds+\frac{C}{\mathfrak{r}^2}\int_0^t\mathcal{E}_\mathfrak{r}^\frac{3}{2}(s)ds.\label{Case1I_4}
\end{align} 
Now we study  the  term $I_5$ of \eqref{Teo7_termsToBound1}. First, notice that by the definition of $\mathfrak{p}_1$ given in \eqref{LocalDescompopressu} and since the kernel $\mathbb{K}(\cdot)=\frac{1}{(-\Delta)}\div\div(\cdot)$ is bounded on $L^\frac{3}{2}(\R)$ (since the Riesz transforms are bounded in such spaces), one has 
\begin{align*}
\|\mathfrak{p}_1\|^{\frac{3}{2}}_{L^{\frac{3}{2}}(B_{0,2r})}\leq \left\|\frac{1}{(-\Delta)}(\div(\div(\mathds{1}_{B_{0,3r}} \vu\otimes \vu)))\right\|_{L^\frac{3}{2}(\R)}^{\frac32}\le C \|\mathds{1}_{B_{0,3r}} |\vu|^2\|^{\frac{3}{2}}_{L^{\frac{3}{2}}(\R)} = C\int_{B_{0,3r}}|\vu|^3dy.
\end{align*} 
Thus, we obtain
\begin{equation}
\sup_{\mathfrak{r}\le r\le \frac1{12}}I_5=\sup_{\mathfrak{r}\le r\le \frac1{12}}\frac{C}{r^2}\int_0^t\int_{B_{0,2r}}|\mathfrak{p}_1|^{\frac{3}{2}}ds\le\sup_{\mathfrak{r}\le r\le \frac1{12}}\frac{C}{r^2}\int_{B_{0,3r}}|\vu|^3dy .\label{termI5}
\end{equation}
Furthermore  since $3r< \frac{1}{4}<\frac{1}{3}$ (recall $\mathfrak{r}\leq r\leq \frac 1 {12}$), we can use the same arguments as for the term $I_4$ and we have
\begin{equation}
\sup_{\mathfrak{r}\le r\le \frac1{12}}I_5\le\frac{1}{4} \mathcal{E}_\mathfrak{r}(t)+\frac{C}{\mathfrak{r}^2}\int_0^t\mathcal{E}^3_\mathfrak{r}(s)ds+\frac{C}{\mathfrak{r}^2}\int_0^t\mathcal{E}_\mathfrak{r}^\frac{3}{2}(s)ds.\label{Case1I_5}
\end{equation}
Now, we study the term $I_6$ of \eqref{Teo7_termsToBound1}. This term is the most technical one, and we will follow the same lines given in the proof of  \cite[Theorem 3.1]{KanMiuTsai19}. Thus, first recall the following estimate for the kernel $\mathbb{K}$: for all $x\in B_{0,2r}$ and $y\in \R\setminus B_{0,3r}$, we have
\begin{equation*}
|\mathbb{K}(x-y)-\mathbb{K}(y)|\le \frac{|x|}{|y|^4}.
\end{equation*}
By the definition of $\mathfrak{p}_2$ given in \eqref{LocalDescompopressu} and since $x\in B_{0,2r}$, we have $|x|<2r$ and 
\begin{align}
|\mathfrak{p}_2(t,x)|=\left| \int_{|y|>3r}(\mathbb{K}(x-y)-\mathbb{K}(-y))(\vu\otimes\vu)(t,y)dy \right|&\le \int_{\R\setminus B_{0,3r}}\frac{|x|}{|y|^4}|\vu|^2dy\notag\\
&\le 2r \int_{\R\setminus B_{0,3r}}\frac{1}{|y|^4}|\vu|^2dy.\label{EstimateP2}
\end{align}
Now, in order to estimate more in detail the  expression above, we need to study the integration domain of the previous integral. Remark that since $3r<\frac{1}{4}$, it is possible to deduce that there exists $\mathcal{N}=\mathcal{N}(r)\in \mathbb{N}$ such that
$$\R\setminus B_{0,3r}\subset \bigcup_{k=1}^\mathcal{N} A_k(r)\cup ( \R\setminus B_{0,\frac{1}{4}}),$$
where $A_k(r)=B_{0,2^kr}\setminus B_{0, 2^{k-1}r}$ and such that for any $1\le k\le \mathcal{N}$ we have $2^{k}r\le \frac{1}{3}$.
Thus, from \eqref{EstimateP2} one has
\begin{align*}
|\mathfrak{p}_2(s,x)|\le  2r \int_{ \bigcup_{k=1}^\mathcal{N} A_k(r)\cup(\R\setminus B_{0,\frac{1}{4}})}\frac{1}{|y|^4}|\vu|^2dy\le 2r\sum_{k=1}^{\mathcal{N}}\int_{A_k(r)}\frac{|\vu|^2}{|y|^4}dy+2r \int_{\R\setminus B_{0,\frac{1}{4}}}\frac{1}{|y|^4}|\vu|^2dy.
\end{align*}
Since  $A_k(r)=B_{0,2^kr}\setminus B_{0, 2^{k-1}r}$, and $2r<\frac{1}{6}<1$ we have 
\begin{align}
|\mathfrak{p}_2(s,x)|&\le\sum_{k=1}^{\mathcal{N}}2r\int_{B_{0,2^kr}\setminus B_{0, 2^{k-1}r}}\frac{|\vu|^2}{|y|^4}dy+  2r \int_{\R\setminus B_{0,\frac{1}{4}}}\frac{1}{|y|^4}|\vu|^2dy\notag\\
&\le \sum_{k=1}^{\mathcal{N}}\frac{2r}{2^{4(k-1)}r^4}\int_{B_{0,2^kr}\setminus B_{0, 2^{k-1}r}}|\vu|^2dy+ \int_{\R\setminus B_{\frac{1}{4}}}\frac{1}{|y|^4}|\vu|^2dy\notag\\
&\le
\sum_{k=1}^{\mathcal{N}}\frac{C}{2^{4(k-1)}r^3}\int_{B_{0,2^kr}}|\vu|^2dy+ \int_{\R\setminus B_{\frac{1}{4}}}\frac{1}{|y|^4}|\vu|^2dy.\label{EstimaP2I2}
\end{align}
Now, we study each term of the expression above separately. For the first one, notice that
by the definition of $E_r(t)$, one has 
\begin{align*}
\frac{1}{2^{4(k-1)}r^3}\int_{B_{0,2^kr}}|\vu|^2dy=\frac{C}{2^{3k}r^2}\left(\frac{1}{2^kr}\int_{B_{0,2^kr}}|\vu|^2dy\right)\le
\frac{1}{2^{3k}r^2}E_{2^kr}(s).
\end{align*}
Thus, since for any $1\le k\le \mathcal{N}$ we have $2^kr\le \frac{1}{3}$, and $\frac{1}{r}<\frac{1}{R}$ it follows that 
\begin{equation*}
\sum_{k=1}^{\mathcal{N}}\frac{C}{2^{4(k-1)}r^3}\int_{B_{0,2^kr}}|\vu|^2dy\le\sum_{k=1}^{\mathcal{N}}\frac{2C}{2^{3k}r^2}E_{2^kr}(s)\le \frac{C}{r^2}\mathcal{E}_R(s)\sum_{k=1}^{\mathcal{N}}\frac{1}{2^{3k}}\le\frac{C}{r^2}\mathcal{E}_R(s).
\end{equation*} 
Now, let us study the second term of \eqref{EstimaP2I2}. Since ${\R\setminus B_{0,\frac{1}{4}}\subset\displaystyle{\bigcup_{k=0}^{+\infty}}B_{0,2^{k-1}}\setminus B_{0,2^{k-2}}}$ we have
\begin{align}
\int_{\R\setminus B_{0,\frac{1}{4}}}\frac{1}{|y|^4}|\vu|^2dy\le \sum_{k=0}^{+\infty}\int_{B_{0,2^{k-1}}\setminus B_{0,2^{k-2}}}\frac{1}{|y|^4}|\vu|^2dy\le\sum_{k=0}^{+\infty}\frac{1}{(2^{k-2})^4}\int_{B_{0,2^{k}}}|\vu|^2dy.\label{Term1EstimateP2I2}
\end{align}
Notice that by a change of variable and since $\vu$ satisfies the estimate \eqref{Hypho_EstimatesAPriori}, we obtain 
\begin{equation*}
\int_{B_{0,2^{k}}}|\vu(t,y)|^2dy=(2^k)^3\int_{B_{0,1}}|\vu(t,\frac{z}{2^k})|^2dz\le 2^{3k}\sup_{0<t<1}\sup_{x\in \R} \int_{B_{x,1}}|\vu(t,y)|^2dy\le  2^{3k}M.
\end{equation*}
By considering the previous equality in \eqref{Term1EstimateP2I2}, we have
\begin{equation}\label{EstimateL2uloc}
\int_{\R\setminus B_{0,\frac{1}{4}}}\frac{1}{|y|^4}|\vu|^2dy\le\sum_{k=0}^{+\infty}\frac{(2^{k})^3}{(2^{k-2})^4}M\le CM\sum_{k=0}^{+\infty}\frac{1}{2^k}\le CM.
\end{equation}
Therefore, by applying \eqref{Term1EstimateP2I2} and the expression above in \eqref{EstimaP2I2}, for any $\mathfrak{r}\le r\le \frac{1}{12}$ we obtain
\begin{align*}
|\mathfrak{p}_2(s,x)|&\le \frac{C}{r^2}\mathcal{E}_\mathfrak{r}(s)+ C M.
\end{align*}
From the previous estimate it follows that for the term $I_6$ of \eqref{Teo7_termsToBound1} we have
\begin{align*}
\sup_{\mathfrak{r}\le r\le \frac1{12}}I_{6}&=\sup_{\mathfrak{r}\le r\le \frac1{12}}\frac{1}{r^2}\int_0^t\int_{B_{0,2r}}|\mathfrak{p}_2|^{\frac{3}{2}}dyds\\
&\le \sup_{\mathfrak{r}\le r\le \frac1{12}} \frac{1}{r^2}\int_0^t\int_{B_{0,2r}}\frac{C}{r^3}\mathcal{E}^\frac{3}{2}_\mathfrak{r}(s)dyds+\sup_{\mathfrak{r}\le r\le \frac1{12}}\frac{1}{r^2}\int_0^t\int_{B_{0,2r}} C M^\frac{3}{2}dyds.\notag
\end{align*}
Moreover, since $|B_{0,2r}|=Cr^3$ and   $r<\frac{1}{1
2}<1$ we have
\begin{equation}\label{Case1I_6}
\sup_{\mathfrak{r}\le r\le \frac1{12}}I_{6}\le \sup_{\mathfrak{r}\le r\le \frac1{12}}\frac{C}{r^2}\int_0^t\mathcal{E}^\frac{3}{2}_\mathfrak{r}(s)ds+\sup_{\mathfrak{r}\le r\le \frac1{12}}rC M^\frac{3}{2} t  \le\frac{C}{\mathfrak{r}^2}\int_0^t \mathcal{E}^{\frac{3}{2}}_\mathfrak{r}(s)ds+ C  M^\frac{3}{2}t.\\[4mm]
\end{equation}

We have finished the study of each term $I_j$ with $1\leq j\leq 6$ given in (\ref{Teo7_termsToBound1}). Thus, gathering the estimates \eqref{Case1I_1}-\eqref{Case1I_3}, \eqref{Case1I_4}, \eqref{Case1I_5} and \eqref{Case1I_6} in \eqref{Teo7_termsToBound1}, and since $\displaystyle{\sup_{\mathfrak{r}\le r\le \frac1{12}}\frac{C}{r} \int_{B_{0,2r}}|\vu_0|^2dy<S}$ by \eqref{Hypo_LocalSmoothing2} (recall also that we have $R\leq \mathfrak{r}$), we have proved that for all $\mathfrak{r}\le r<\frac{1}{12}$ we have
\begin{align*}
\sup_{\mathfrak{r}\le r\le \frac1{12}}{E}_r(t)&\le S+\frac{C}{\mathfrak{r}^2}\int_0^t\mathcal{E}_\mathfrak{r}(s)ds+\frac{1}{4} \mathcal{E}_\mathfrak{r}+Ct\\
&+Ct+ \frac{1}{\mathfrak{r}^2}\int_0^t \mathcal{E}_{\mathfrak{r}}(s)ds+\frac{1}{4} \mathcal{E}_\mathfrak{r}(t)+\frac{C}{\mathfrak{r}^2}\int_0^t\mathcal{E}^3_\mathfrak{r}(s)ds+\frac{C}{\mathfrak{r}^2}\int_0^t\mathcal{E}_\mathfrak{r}^\frac{3}{2}(s)ds.\\
&+\frac{1}{4} \mathcal{E}_\mathfrak{r}(t)+\frac{C}{\mathfrak{r}^2}\int_0^t\mathcal{E}^3_\mathfrak{r}(s)ds+\frac{C}{\mathfrak{r}^2}\int_0^t\mathcal{E}_\mathfrak{r}^\frac{3}{2}(s)ds+ CM^\frac{3}{2}t +\frac{C}{R^2}\int_0^t \mathcal{E}^3_\mathfrak{r}(s)ds\\[2mm]
&\le S+Ct+CM^\frac32 t+\frac{3}{4}\mathcal{E}_\mathfrak{r}(t)+\frac{C}{\mathfrak{r}^2}\int_0^t\mathcal{E}_\mathfrak{r}(s)+\mathcal{E}^3_\mathfrak{r}(s)+\mathcal{E}^\frac32_\mathfrak{r}(s)ds,
\end{align*}
and we rewrite the previous estimate as follows:
\begin{equation}\label{firstHalf_LocalSmoothing}
\sup_{\mathfrak{r}\le r\le \frac1{12}}{E}_r(t)\le S +\mathfrak{C}_1(M)t+ \frac{3}{4}\mathcal{E}_\mathfrak{r}(t)+\frac{C}{\mathfrak{r}^2}\int_0^t\mathcal{E}_\mathfrak{r}(s)+\mathcal{E}^3_\mathfrak{r}(s)+\mathcal{E}^\frac32_\mathfrak{r}(s)ds.
\end{equation}
Now, let us study the second term of the right-hand side of \eqref{SplitER}.
\item[$\bullet$]\underline{ Assume $\frac{1}{12}< r\le\frac{1}{3}$}: similar to the previous case, notice that from \eqref{Teo7_termsToBound}, we have 
\begin{equation}\label{Teo7_termsToBound2}
\sup_{\frac1{12}<r\le \frac13}E_r(t)\le \sup_{\frac1{12}<r\le \frac13} \frac{C}{r} \int_{B_{0,2r}}|\vu_0|^2 dy+\sup_{\frac1{12}<r\le \frac13}\sum_{j=1}^{6}I_j
\end{equation} 
For the  term $I_1$ of \eqref{Teo7_termsToBound2} as we have  $B_{0,2r}\subset B_{0,1}$ since $2r\le\frac{2}{3}<1$, since $\frac{1}{r}<12$ and since $\|\vu\|^2_{L^\infty_tL_x^2(Q_1(1,0))}\le M$ by the hypothesis  \eqref{Hypho_EstimatesAPriori}, it follows that 
\begin{equation}\label{Case2I_1}
\sup_{\frac1{12}<r\le \frac13}I_{1}=\sup_{\frac1{12}<r\le \frac13}\frac{C}{r^3}\int_{0}^t \int_{B_{0,2r}}|\vu|^2 dyds\le C\int_{0}^t \int_{B_{0,1}}|\vu|^2 dyds\le C M t.
\end{equation}
For the term $I_2$ in \eqref{Teo7_termsToBound2}, by the Cauchy-Schwarz inequality  and since $\frac{1}{r}\le 12$ and  $B_{0,2r}\subset B_{0,1}$, one has 
\begin{align*}
\sup_{\frac1{12}<r\le \frac13}I_{2}=\frac{C}{r} \int_{0}^t \int_{B_{0,2r}}|\rot \vu| |\vw| dyds&\le C\int_0^t \left(\int_{B_{0,2r}} |\rot \vu| ^2dy\right)^{\frac{1}{2}}
\left(\int_{B_{0,2r}} |\vw| ^2dy\right)^{\frac{1}{2}}ds\nonumber\\
&\le C\|\vw\|_{L^\infty_tL^2_x(Q_1(1,0))} \int_0^t\|\rot \vu(s,\cdot)\|_{L^2(B_{0,1})} ds.
\end{align*}
Moreover,  since $|\rot\vw|^2\le 2|\grad\otimes\vu|^2$ and by the Cauchy-Schwarz inequality in the time variable we obtain
\begin{align*}
\sup_{\frac1{12}<r\le \frac13}I_{2}\le C\|\vw\|_{L^\infty_tL^2_x(Q_1(1,0))}\int_0^t \|\grad\otimes \vu(s,\cdot)\|_{L^2(B_{0,1})}ds
&\le  C \|\vw\|_{L^\infty_tL^2_x(Q_1(1,0))}  \|\grad\otimes \vu\|_{L^2_tL^2_x(Q_1(1,0))}t^\frac{1}{2}.
\end{align*}
Therefore, since $\|\grad\otimes \vu\|^2_{L^2_tL_x^2(Q_1(1,0))}<M$ by  \eqref{Hypho_EstimatesAPriori},  $ \|\vw\|_{L^\infty_tL^2_x(Q_1(1,0))}^2<CR<C$ by \eqref{Hypho_EstimatesOmega} (recall $R<1$),  we have
\begin{equation}
\sup_{\frac1{12}<r\le \frac13} I_2\le C M ^\frac12 t^\frac12. \label{Case2I_2}
\end{equation}
For the term $I_3$ of \eqref{Teo7_termsToBound2}, since  $\frac{1}{r}\le 12$ and by the H\"older inequality we have
\begin{align*}
\sup_{\frac1{12}<r\le \frac13}I_3=\sup_{\frac1{12}<r\le \frac13}\frac{C}{r^2}\int_{0}^t \int_{B_{0,2r}}| \vw| |\vu| dyds&\le C\sup_{\frac1{12}<r\le \frac13}\int_0^t \left(\int_{B_{0,2r}} |\vw| ^2dy\right)^{\frac{1}{2}}
\left(\int_{B_{0,2r}} |\vu| ^2dy\right)^{\frac{1}{2}}ds.
\end{align*}
Then, since  $B_{0,2r}\subset B_{0,1}$, one has
\begin{equation*}
\sup_{\frac1{12}<r\le \frac13}I_3 \le  C\int_0^t  \|\vw(s,\cdot)\|_{L^2(B_{0,1})} \|\vu(s,\cdot)\|_{L^2(B_{0,1})}ds.
\end{equation*}
Thus, using the fact that $\|\vu\|^2_{L^\infty_tL_x^2(Q_1(1,0))}<M$ by  \eqref{Hypho_EstimatesAPriori} and  $ \|\vw\|_{L^\infty_tL_x^2(Q_1(1,0))}<C$, one has
\begin{align}
\sup_{\frac1{12}<r\le \frac13}I_3\le  C M^\frac12 t.\label{Case2I_3}
\end{align}
For the term $I_{4}$ in \eqref{Teo7_termsToBound2}, using the same arguments as in \eqref{ControlNormL3} we have. 
\begin{align*}
\sup_{\frac1{12}<r\le \frac13}I_{4}&=\sup_{\frac1{12}<r\le \frac13}\frac{C}{r^2}\int_0^t\int_{B_{0,2r}}|\vu|^3dy ds\le \sup_{\frac1{12}<r\le \frac13}\frac{C}{r^2} \int_{0}^t\|\grad \otimes\vu(s,\cdot)\|_{L^2(B_{0,2r})}^{\frac{3}{2}}\|\vu(s,\cdot)\|^{\frac{3}{2}}_{L^2(B_{0,2r})}ds\\
&\;\;+\sup_{\frac1{12}<r\le \frac13} \frac{C}{r^2} \int_{0}^t\|\vu(s,\cdot)\|_{L^2(B_{0,2r})}^{3}ds.\notag
\end{align*}
Using the fact $\frac{1}{r^2}\le C$ and  $B_{0,2r}\subset B_{0,1}$, we obtain
\begin{align*}
\sup_{\frac1{12}<r\le \frac13}I_{4} &\le C\|\vu\|^{\frac{3}{2}}_{L^\infty_t L_x^2(Q_1)}\int_{0}^t \|\grad \otimes\vu(s,\cdot)\|_{L^2(B_{0,1})}^\frac32ds+t\|\vu\|^{3}_{L^\infty_t L_x^2(Q_1)}.
\end{align*}
Thus,  by the Hölder inequality $(1=\frac{3}{4}+\frac{1}{4})$ in the time variable, and since $2r<1$, and $t<1$, one has
\begin{align}
\sup_{\frac1{12}<r\le \frac13}I_{4}&\le C t^\frac{1}{4} \|\vu\|^{\frac{3}{2}}_{L^\infty_t L_x^2(Q_1)}\left(\int_{0}^t \|\grad \otimes\vu(s,\cdot)\|_{L^2(B_{0,2r})}^2 ds\right)^{\frac{3}{4}}+t\|\vu\|^{3}_{L^\infty_t L_x^2(Q_1)}\notag\\
&\le C t^\frac{1}{4} \|\vu\|^{\frac{3}{2}}_{L^\infty_t L_x^2(Q_1)}\|\grad \otimes\vu\|_{L^2_tL_x^2(Q_{1})}^\frac{3}{2}+t\|\vu\|^{3}_{L^\infty_t L_x^2(Q_1)}.\notag
\end{align}
Since $\|\vu\|^2_{L^\infty_tL_x^2(Q_1(1,0))}\le M$  and $\|\grad\otimes\vu\|^2_{L^2_tL_x^2(Q_1(1,0))}\le M$ by  \eqref{Hypho_EstimatesAPriori}, we have
\begin{align}
\sup_{\frac1{12}<r\le \frac13}I_4\le & C  M^\frac{3}{2} t^\frac{1}{4}+  M^\frac{3}{2}t\label{Case2I_4}.
\end{align}
For the term $I_5$ of \eqref{Teo7_termsToBound2}, recall that from \eqref{termI5}, we obtain
\begin{equation*}
\sup_{\frac1{12}<r\le \frac13}I_5=\frac{C}{r^2}\int_0^t\int_{B_{0,2r}}|\mathfrak{p}_1|^{\frac{3}{2}}dyds\le\frac{C}{r^2}\int_0^t\int_{B_{0,3r}}|\vu|^3dy ds.
\end{equation*} 
Since $3r\le1$, we can apply the same arguments as in \eqref{Case2I_4}, and  we have
\begin{equation}
\sup_{\frac1{12}<r\le \frac13}I_5\le M^\frac{3}{2} t^\frac{1}{4} + M^\frac{3}{2}t.
\end{equation}
For the last term $I_6$ of \eqref{Teo7_termsToBound2}, first notice that by the estimate \eqref{EstimateP2} and since $\frac{1}{4}\le 3r$, we get
\begin{align*}
|\mathfrak{p}_2(s,x)|\le C \int_{\R\setminus B_{0,3r}}\frac{1}{|y|^4}|\vu|^2dy
\le  C \int_{\R\setminus B_{0,\frac{1}{4}}}\frac{1}{|y|^4}|\vu|^2dy.
\end{align*}
Therefore, using the same arguments as in \eqref{EstimateL2uloc}, one has $|\mathfrak{p}_2(s,x)|\le C M.$ Hence, since $r\le\frac{1}{3}<1$, it follows that
\begin{align}
\sup_{\frac1{12}<r\le \frac13}I_{6}=\frac{1}{r^2}\int_0^t\int_{B_{0,2r}}|\mathfrak{p}_2|^{\frac{3}{2}}dyds\le \sup_{\frac1{12}<r\le \frac13}r C M^\frac{3}{2}t\le C M^\frac{3}{2}t .
\label{Case2I_6}
\end{align}
Thus, gathering the estimates \eqref{Case2I_1}-\eqref{Case2I_6} in \eqref{Teo7_termsToBound}, since $\displaystyle{\frac{C}{r} \int_{B_{0,2r}}|\vu_0|^2dy<S}$ by \eqref{Hypo_LocalSmoothing2}, and $t<t^\frac{1}{2}<t^\frac{1}{4}$ due to $t<1$ it follows that  for all $\frac{1}{12}\le r\le \frac{1}{3}$, we have
\begin{align}
\sup_{\frac1{12}<r\le \frac13}{E}_r(t)&\le S+CMt+CM^\frac12 t^\frac12+ CM^\frac12 t+CM^\frac32t^\frac14+CM^\frac32t+ C M^\frac{3}{2}t\notag\\
&\le S +\mathfrak{C}_2(M)t^\frac{1}{4},\label{secondHalf_LocalSmoothing}
\end{align}
and this finishes the study of the previous quantity in the case when $\frac{1}{12}< r\le\frac{1}{3}$.\\
\end{itemize}
Thus, applying the estimates \eqref{firstHalf_LocalSmoothing} and \eqref{secondHalf_LocalSmoothing} in \eqref{SplitER}, we have proved that
\begin{align*}
\mathcal{E}_\mathfrak{r}(t)&\le\sup_{\mathfrak{r}\le r\le \frac1{12}}E_r(t)+\sup_{\frac1{12}<r\le \frac13}E_r(t)\\& \le S +\mathfrak{C}_1(M)t+ \frac{3}{4}\mathcal{E}_\mathfrak{r}(t)+\frac{C}{R^2}\int_0^t\mathcal{E}_\mathfrak{r}(s)+\mathcal{E}^3_\mathfrak{r}(s)+\mathcal{E}^\frac32_\mathfrak{r}(s)ds+ S +\mathfrak{C}_2(M)t^\frac{1}{4}\\[2mm]
&\le2S+\mathfrak{C}_1(M)t +\mathfrak{C}_2(M)t^\frac{1}{4}+ \frac{3}{4}\mathcal{E}_\mathfrak{r}(t)+\frac{C}{\mathfrak{r}^2}\int_0^t\mathcal{E}_\mathfrak{r}(s)+\mathcal{E}^3_\mathfrak{r}(s)+\mathcal{E}^\frac32_\mathfrak{r}(s)ds.
\end{align*}
Fix now the time $0<\mathcal{T}_1\leq 1$ such that
\begin{equation}\label{Def_T1}
\mathcal{T}_1=\min\{1,\frac{S}{2\mathfrak{C}_1},\frac{S^4}{2\mathfrak{C}_2^4}\}.
\end{equation}
Notice that for all $t< \mathcal{T}_1$, we have $\mathfrak{C}_1t\le \frac S2$ and $\mathfrak{C}_2t^\frac 14\le \frac S2$, hence it follows that
\begin{equation*}
\frac14\mathcal{E}_\mathfrak{r}(t)\le 3S+\frac{C}{\mathfrak{r}^2}\int_0^t\mathcal{E}_\mathfrak{r}(s)+\mathcal{E}^3_\mathfrak{r}(s)+\mathcal{E}^\frac32_\mathfrak{r}(s)ds.
\end{equation*}
Observing that if $\mathcal{E}_\mathfrak{r}>1$ we have $\mathcal{E}^\frac32_\mathfrak{r}<\mathcal{E}^3_\mathfrak{r}$ and if $\mathcal{E}_\mathfrak{r}\le 1$ we have $\mathcal{E}^\frac32_\mathfrak{r}<\mathcal{E}_\mathfrak{r}$, it is then enough to study for any $t< \mathcal{T}_1$ the expression
\begin{align}\label{EstimateBeforeGronwall}
\mathcal{E}_\mathfrak{r}(t)\le 12S +\frac{C}{\mathfrak{r}^2}\int_0^t\mathcal{E}_\mathfrak{r}(s)+\mathcal{E}_\mathfrak{r}^3(s)ds.
\end{align}
In order to estimate more in detail the expression above,  we can use the following Gronwall-type inequality
%%%%%%%%%%%%%%%%%%%%%%%%%%%%%%%%%%%%%%%%%%%%%%%%%%%
\begin{Lemma}\label{Lem_Gronwall}
Let $f\in L_{loc}^{\infty}([0,T_1[)$ be a function such that for all $t\in ]0,T_1[$, for some $ a,b>0$ and $m\ge1 $, we have
$$f(t)\le a+b \int _0^t(f(s)+f^m(s))ds.$$
Then, there exists a universal constant $\mathfrak{\mathfrak{c}}>1$ such that for all $t\in ]0,T]$ with $T= \min \left\{T_1, \frac{\mathfrak{c}}{b(1+a^{m-1})}\right\}$, we have $f(t)\le 2a $ .
\end{Lemma}
%%%%%%%%%%%%%%%%%%%%%%%%%%%%%%%%%%%%%%%%%%%%%%%%%%%
\noindent For a proof of this result we refer to \cite[Lemma 2.2]{BradTsai}. Now, by applying the previous lemma to the expression \eqref{EstimateBeforeGronwall} with  $a=12S$, $b=\frac C{\mathfrak{r}^2}$, $m=3$ and $T_1=\mathcal{T}_1$ given in (\ref{Def_T1}), there exists a universal constant  $\mathfrak{\mathfrak{c}}>1$ such that for
\begin{equation*}
T^\ast=\min\{\mathcal{T}_1, \mathfrak{\mathfrak{c}} \lambda  \mathfrak{r}^2\},\quad  \text{where}\quad \lambda=\frac{1}{(1+S^2)},
\end{equation*}   
we have for all $0<t<T^\ast$ the estimate
$$\mathcal{E}_\mathfrak{r}(t)\le 24S.$$
Since we have $E_{\mathfrak{r}}(t)\le \mathcal{E}_\mathfrak{r}(t)=\underset{\mathfrak{r}\leq r\leq \frac{1}{3}}{\sup}E_r(t)$, we finally obtain
\begin{equation}\label{ConclusionLemmaLocalSmoothing}
E_{\mathfrak{r}}(t)\le CS,
\end{equation} which finishes the proof of the Lemma \ref{Lem_localSmoothing2}. \hfill $\blacksquare$

%%%%%%%%%%%%%%%%%%%%%%%%%%%%%%%%%%%%%%%%%%%%%%%%%%%
\begin{Corollary}\label{Coro_localSmothing}
Under the hypothesis of Lemma \ref{Lem_localSmoothing2}, for any $R\le r\leq \frac 13$,  such that $\sqrt{\lambda} r\le \sqrt{\frac{\mathcal{T}_1}{\mathfrak{c}}}$ where $\mathcal{T}_1$ is given in \eqref{Def_T1}, $\lambda=\frac{1}{1+S^2}$ and  $\mathfrak{c}>1$, we have
\begin{equation*}
\frac{1}{r^2}\int_0^{\lambda r^2}\int_{B_{0,r}}|\vu|^3+|p-\mathfrak{h}|^\frac{3}{2}dyds<C(S^\frac{3}{2}+S).
\end{equation*}
\end{Corollary}
%%%%%%%%%%%%%%%%%%%%%%%%%%%%%%%%%%%%%%%%%%%%%%%%%%%
\noindent \textbf{Proof.} First note that we have $0<R<\frac{1}{12}$ is arbitrary and it can be chosen arbitrarily small, so the conditions $R\le r<\frac 13$ and $\sqrt{\lambda} r\le \sqrt{\frac{\mathcal{T}_1}{\mathfrak{c}}}$ are compatible. Now, remark that from the estimate  \eqref{ControlNormL3}, we have 
\begin{align*}
\frac{1}{r^2}\int_0^{\lambda r^2}\int_{B_{0,r}}|\vu|^3dyds&\le\frac{1}{r^2}\int_{0}^{\lambda r^2}\|\grad \otimes\vu(s,\cdot)\|_{L^2(B_{0,r})}^{\frac{3}{2}}\|\vu(s,\cdot)\|^{\frac{3}{2}}_{L^2(B_{0,r})}ds+\frac{1}{r^2}\int_{0}^{\lambda r^2} \|\vu(s,\cdot)\|_{L^2(B_{0,r})}^{3}ds\\
&\le \frac{1}{r^2}\sup_{0<s<\lambda r^2}\|\vu(s,\cdot)\|^{\frac{3}{2}}_{L^2(B_{0,r})}\int_{0}^{\lambda r^2}\|\grad \otimes\vu(s,\cdot)\|_{L^2(B_{0,r})}^{\frac{3}{2}}ds\\
&\quad+C\lambda\sup_{0<s<\lambda r^2}\|\vu(s,\cdot)\|^{3}_{L^2(B_{0,r})}.
\end{align*}
Using the Hölder inequality in the time variable $(1=\frac{3}{4}+\frac{1}{4})$, one has
\begin{align}
\frac{1}{r^2}\int_0^{\lambda r^2}\int_{B_{0,r}}|\vu|^3dyds&\le \frac{C\lambda^\frac{1}{4}}{r^\frac{3}{2}}\sup_{0<s<\lambda r^2}\|\vu(s,\cdot)\|^{\frac{3}{2}}_{L^2(B_{0,r})}\left( \int_{0}^{\lambda r^2}\|\grad \otimes\vu\|_{L^2(B_{0,r})}^{2}ds\right)^\frac{3}{4} \notag\\
&\quad+C\lambda\sup_{0<s<\lambda r^2}\|\vu(s,\cdot)\|^{3}_{L^2(B_{0,r})},\notag
\end{align}
which can be rewritten as follows
\begin{align*}
\frac{1}{r^2}\int_0^{\lambda r^2}\int_{B_{0,r}}|\vu|^3dyds&\le C\lambda^\frac{1}{4}\sup_{0<s<\lambda r^2}\left(\frac{1}{r}\int _{B_{0,r}}|\vu|^2dy\right)^\frac{3}{4}\left( \frac{1}{r}\int_{0}^{\lambda r^2}\int_{B_{0,r}}|\grad\otimes\vu|^2dyds\right)^\frac{3}{4} \\
&\quad+C\lambda r^\frac32\sup_{0<s<\lambda r^2}\left(\frac{1}{r}\int_{B_{0,r}}|\vu(s,y)|^2dy\right)^\frac{3}{2}.
\end{align*}
Since $E_r(t)=\displaystyle{\sup_{0<s<t}\frac{1}{r}\int _{B_{0,r}}|\vu|^2dy+\frac{1}{r}\int_0^t \int_{B_{0,r}}|\grad\otimes\vu|^2dyds +\frac{1}{r^2}\int_0^t\int_{B_{0,r}}|p-\mathfrak{h}(t)|^{\frac{3}{2}}dyds}$, and if we set $t=\lambda r^2$, we thus can write
\begin{align*}
\frac{1}{r^2}\int_0^{\lambda r^2}\int_{B_{0,r}}|\vu|^3dyds&\le C\lambda^{\frac14} E_r^{\frac{3}{2}}(\lambda r^2)+C\lambda r^\frac32 E_r^{\frac{3}{2}}(\lambda r^2)\le C(\lambda^{\frac14} +r^\frac32\lambda)E_r^{\frac{3}{2}}(\lambda r^2).
\end{align*}
Moroever, since  $r<1$ and $\lambda<1$, we have
\begin{align*}
\frac{1}{r^2}\int_0^{\lambda r^2}\int_{B_{0,r}}|\vu|^3dyds\le C E_r^{\frac{3}{2}}(\lambda r^2).
\end{align*}
Since $E_r(t)$ is an increasing fonction in $t$, and  $\mathfrak{c}>1$ we have $ E_r^{\frac{3}{2}}(\lambda r^2)< E_r^{\frac{3}{2}}(\mathfrak{c}\lambda r^2)$. Thus, since we have considered $0<r<1$ such that $\mathfrak{c}\lambda r^2\le \mathcal{T}_1$, and  it follows from \eqref{ConclusionLemmaLocalSmoothing} that $E_r(\mathfrak{c}\lambda r^2)<CS$ and then 
\begin{equation*}
\frac{1}{r^2}\int_0^{\lambda r^2}\int_{B_{0,r}}|\vu|^3dyds\le  C E_r^{\frac{3}{2}}(\mathfrak{c}\lambda r^2)<CS^\frac32.
\end{equation*}
On the other hand, from the definition of $E_{r}$ we immediately have
\begin{equation*}
\frac{1}{r^2}\int_{0}^{\lambda r^2}\int_{B_{0,r}}|p-\mathfrak{h}|^{\frac{3}{2}}dyds\le E_r(\lambda r^2)\le E_r(\mathfrak{c}\lambda r^2)\le C S.
\end{equation*}
Gathering the previous two estimates we find that for  $\lambda=\frac{1}{1+S^2}$ and for any  $R\leq r\leq \frac13$ such that $\sqrt{\lambda} r\le \mathcal{T}_1$, one has
\begin{equation*}
\frac{1}{r^2}\int_0^{\lambda r^2}\int_{B_{0,r}}|\vu|^3+|p-\mathfrak{h}|^\frac{3}{2}dyds<C(S^\frac{3}{2}+S),
\end{equation*}
which finishes the proof of Corollary \ref{Coro_localSmothing}. \hfill $\blacksquare$
\begin{Remark}
It is worth noting that throughout the proofs of the Lemma \ref{Lem_localSmoothing2} and Corollary \ref{Coro_localSmothing}, we treated $\vw$ as an external force and we can simply consider that $(\vu,p)$ is a \emph{local Leray} solution as in Remark \ref{Rem_PartialLocalEqualLocal} such that the hypotheses \eqref{Hypho_EstimatesOmega} over $\vw$ are satisfied.
\end{Remark}
%%%%%%%%%%%%%%%%%%%%%%%%%%%%%%%%%%%%%%%%%%%%%%%%%%%
\noindent \textbf{Proof of the Theorem \ref{Theo_BlowUpConcentration}.} 
Let $(\vu,p,\vw)$ be a Leray-type weak solution of the micropolar fluids equations \eqref{MicropolarFluidsEquationsEqua1} and \eqref{MicropolarFluidsEquationsEqua2}. Let ${\mathcal{T}}>0$ be the maximal time such that $\vu\in \mathcal{C}(]0,\mathcal{T}[, L^\infty(\R))$  and  the point $(\mathcal{T},0)$ is  a partial singular point in the sense of the Definition \ref{Def_PartialRegularPoints}. Assume that for  a fixed $r_0>0$ with ${0<\mathcal{T}-r_0^2},$ we have
\begin{equation}\label{HypoTpyI}
\sup_{x_0\in \R}\sup _{r\in ]0,r_0]}\sup _{t\in]\mathcal{T}-r^2,\mathcal{T}]}\frac{1}{r}\int_{B_{x_0,r}}|\vu(t,x)|^2dx =\mathfrak{M} <+\infty.
\end{equation}
Our aim consists in proving that there exits $\varepsilon>0$, $\mathfrak{S}=\mathfrak{S}(\mathfrak{M})$ and $\delta>0$ such that for all $t\in ]\mathcal{T}-\delta,\mathcal{T}[$ we have
\begin{equation}\label{NormConcentration}
\int_{B_{0,\sqrt{\frac{\mathcal{T}-t}{\mathfrak{S}}}}}|\vu(t,x)|^3dx\ge \varepsilon.
\end{equation}
First, notice that it is enough to show that there exits $\varepsilon_\ast>0$, $\mathfrak{S}=\mathfrak{S}(\mathfrak{M})$ and $\delta>0$ such that for all $t\in ]\mathcal{T}-\delta,\mathcal{T}[$, we have
\begin{equation}\label{NormL2Explotion}
\frac{1}{\sqrt{{\mathcal{T}}-t}}\int_{B_{0,\sqrt{\frac{{\mathcal{T}}-t}{\mathfrak{S}}}}}|\vu(t,x)|^2dx\ge \varepsilon_\ast,
\end{equation}
indeed, if \eqref{NormL2Explotion} holds, by the Hölder inequality $(1=\frac{2}{3}+\frac{1}{3})$ we have
\begin{equation*}
\varepsilon_\ast\le\frac{1}{\sqrt{{\mathcal{T}}-t}}\int_{B_{0,\sqrt{\frac{{\mathcal{T}}-t}{\mathfrak{S}}}}}|\vu(t,x)|^2dx\le \frac{C}{\sqrt{\mathfrak{S}}}\left( \int_{B_{0,\sqrt{\frac{{\mathcal{T}}-t}{\mathfrak{S}}}}}|\vu(t,x)|^3dx\right)^{\frac{2}{3}},
\end{equation*}
which in turn implies \eqref{NormConcentration} with $\varepsilon=\frac{{(\varepsilon_\ast)}^\frac{3}{2}\mathfrak{S}^\frac{3}{4}}{C}$.\\

Now, for proving \eqref{NormL2Explotion} we will use a contradiction argument. Thus, assume that for all $\varepsilon_\ast>0$, for all $0<\mathfrak{S}<1$ and for $\delta=\min\{\sqrt{\frac{\mathcal{T}}{2}}, \frac{r_0\mathfrak{S}^\frac32}{2}\}$ there exists $ {\mathcal{T}}-\delta^2\le t_0\le \mathcal{T}$ such that
\begin{equation}\label{controlSingularPoint}
\frac{1}{\sqrt{\mathcal{T}-t_0}}\int_{B_{0,\sqrt{\frac{\mathcal{T}-t_0}{\mathfrak{S}}}}}|\vu(t_0,x)|^2dx<\varepsilon_\ast.
\end{equation}
The strategy consists in  applying a particular scaling limit procedure to  the solution $(\vu,p,\vw)$ in order to obtain that  the point $(\mathcal{T},0)$ is partially regular in the sense of Definition \ref{Def_PartialRegularPoints} leading us to the wished contradiction.\\

Thus, let $0<\mathfrak{S}<1$ to be fixed later and consider $\gamma=\sqrt{\frac{\mathcal{T}-t_0}{\mathfrak{S}}}$. Notice that since $ {\mathcal{T}}-\delta^2< t_0< {\mathcal{T}}$ and $\delta^2\le r_0^2\mathfrak{S}^3\le r_0^2\mathfrak{S}$ (recall $0<\mathfrak{S}<1$), we have
\begin{equation}\label{IntervalGamma}
\gamma=\sqrt{\frac{\mathcal{T}-t_0}{\mathfrak{S}}}< \sqrt{\frac{\delta^2}{\mathfrak{S}}}< r_0.
\end{equation} 
Now, we scale the functions $\vu$, $p$ and $\vw$ as follows: for all $(s,y)\in [0,\mathfrak{S}[\times \R$, we consider  
$$\vec{u}^\gamma(s,y)=\gamma \vu(t_0+\gamma^2 s,\gamma y),\;\; p^{\gamma}(s,y)=\gamma^2 p(t_0+\gamma^2 s,\gamma y)\;\; \text{and}\;\; 
\vw^{\gamma}(s,y)=\gamma^2 \vw(t_0+\gamma^2 s,\gamma y),$$ 
recall that the first equation of the micropolar fluids equations \eqref{MicropolarFluidsEquationsEqua1} is invariant under the previous scaling. Remark also that since $\vu\in \mathcal{C}(]0,\mathcal{T}[,L^\infty(\R))$, we have that $(\vu^\gamma,p^\gamma)$ is a strong solution of the system above.\\

We want now to apply Lemma \ref{Lem_localSmoothing2} and Corollary \ref{Coro_localSmothing} in order to obtain that there exists $\tilde \varepsilon>0$ and $\rho>0$ such that
\begin{equation*}
\frac{1}{\rho^2} \int_{\mathfrak{S}-\rho^2}^{\mathfrak{S}}\int_{B_{0,\rho}}|\vec{u}^\gamma|^3+|p^{\gamma}-\mathfrak{h}^\gamma|^\frac{3}{2}dyds \le\tilde \varepsilon,
\end{equation*}
where $\mathfrak{h}^{\gamma}=\gamma^2\mathfrak{h}(\gamma^2\cdot)$, and $\mathfrak{h}$ is given by the local decomposition of the pressure (see Lemma \ref{Propo_LocalDescompositionPressure}). Then by re-scaling back to the variables $(\vu,p)$ and using the $\varepsilon$-regularity theory developed in the appendix \ref{Sec_PartialRegularity}, we will be able to deduce that the point $(\mathcal{T},0)$ is partially regular which is a contradiction.\\

Since we want to apply Lemma  \ref{Lem_localSmoothing2}, we need some information on the initial data and for this, since $\vu \in \mathcal{C}(]0,\mathcal{T}[,L^\infty(\R))$, we can consider $\vec{u}^\gamma_0(\cdot)=\gamma\vu(t_0,\gamma\cdot)$ as initial data such that $(\vec{u}^\gamma,p^{\gamma})$ is a solution of the forced Navier-Stokes equations
\begin{equation}\label{NSrotvW}
\begin{cases}
\partial_t \vec{u}^\gamma =\Delta \vec{u}^\gamma-(\vec{u}^\gamma\cdot\grad)\vec{u}^\gamma-\grad p^{\gamma} +\frac{1}{2}\rot\vw^{\gamma},\\[2mm]
\vec{u}^\gamma(0,\cdot)=\vec{u}^\gamma_0.
\end{cases}
\end{equation}
Furthermore, to deduce that $(\vec{u}^\gamma,p^{\gamma})$ is a local Leray solution of the system above, it will convenient to write $\rot \vw^{\gamma} = \div(\mathbb{W}^\gamma)$ where we have 
\begin{equation}\label{ExternalForceW}
{\mathbb{W}^\gamma=\left(\begin{matrix}
0&\omega^\gamma_3& -\omega^\gamma_2\\
\omega^\gamma_1 & 0 & -\omega^\gamma_3\\
-\omega^\gamma_2 & \omega^\gamma_1 & 0
\end{matrix}\right)},
\end{equation}
and we thus obtain that the pair $(\vec{u}^\gamma,p^\gamma)$ is a solution of the following system
\begin{equation*}
\begin{cases}
\partial_t \vec{u}^\gamma =\Delta \vec{u}^\gamma-(\vec{u}^\gamma\cdot\grad)\vec{u}^\gamma-\grad p^{\gamma} +\frac{1}{2}\div(\mathbb{W}^\gamma),\\[2mm]
\vec{u}^\gamma(0,\cdot)=\vec{u}^\gamma_0.
\end{cases}
\end{equation*}
At this point, we can apply the theory of local Leray solutions of the Navier-Stokes equations which is given in the following result:
%%%%%%%%%%%%%%%%%%%%%%%%%%%%%%%%%%%%%%%%%%%%%%%%%%%
\begin{Lemma}\label{Theo_LocalLeray}
Let $\vv_0$ be an initial data and $\mathbb{F}$ be a tensor field such that
\begin{equation*}
\sup_{x_0\in \R}\int_{B_{x_0,1}}|\vv_0(y)|^2dy<M^\ast\quad\text{and}\quad\sup_{x_0\in \R} \int_0^{1}\int _{B_{x_0,1}}|\mathbb{F}|^2dyds<M^\ast.
\end{equation*}
Then, there exists a local Leray solution $(\vv,q)$ in the sense of the Definition 14.1 of the book \cite{PGLR1}  of the forced Navier-Stokes equations
\begin{equation*}
\begin{cases}
\partial_t\vv =\Delta \vv-(\vv\cdot\grad)\vv-\grad q +\div(\mathbb{F}),\\[2mm]
\vv(0,\cdot)=\vv_0,
\end{cases}
\end{equation*}
on $]0,T[\times \R$, such that $T=\min\{1,\frac{1}{C(1+M^\ast)^4}\}$ and 
\begin{equation*}
\sup_{0<t<T}\sup_{x\in \R}\int_{B_{x,1}}|\vv(t,y)|^2dy+\sup_{x\in \R}\int_{0}^{T}\int_{B_{x,1}}|\grad\otimes\vv|^2dy\le C(M^\ast).
\end{equation*}
\end{Lemma}
%%%%%%%%%%%%%%%%%%%%%%%%%%%%%%%%%%%%%%%%%%%%%%%%%%%
\noindent For a proof of the previous lemma we refer to \cite[Theorem 14.1, pg 455 ]{PGLR1}.
As we can see, with this lemma at hand we can construct a local Leray solution as long as we have some mild decay on the initial data $\vec{u}^\gamma_0$ and on the external force $\mathbb{W}^\gamma$ and for this we only need to verify the following uniform controls 
\begin{equation}\label{EstimateInitialData}
\sup_{x_0\in \R}\int_{B_{x_0,1}}|\vec{u}^\gamma_0(y)|^2dy< +\infty,
\end{equation} 
and 
\begin{equation}\label{EstimateForce}
\sup_{x_0\in \R} \int_0^{1}\int _{B_{x_0,1}}|\mathbb{W}^\gamma|^2dyds<+\infty.
\end{equation}
Let us study the initial data. Since  $\vec{u}^\gamma_0(\cdot)=\gamma\vu(t_0,\gamma\cdot)$, by a change of variable, we obtain 
\begin{equation*}
\sup_{x_0\in \R}\int_{B_{x_0,1}}|\vec{u}^\gamma_0(y)|^2dy=\sup_{x_0\in \R}\frac{1}{\gamma}\int_{B_{\gamma x_0,\gamma}}|\vu(t_0,y)|^2dy=\sup_{z\in \R}\frac{1}{\gamma}\int_{B_{z,\gamma}}|\vu(t_0,y)|^2dy.
\end{equation*}
On the other hand, recall that by \eqref{IntervalGamma} we have $\gamma<r_0$ and $\mathcal{T}-\gamma^2<t_0<\mathcal{T}$. Thus, by the hypothesis \eqref{HypoTpyI}, we have
\begin{equation*}
\sup_{z\in \R}\frac{1}{\gamma}\int_{B_{z,\gamma}}|\vu(t_0,y)|^2dy\le \sup_{x_0\in \R}\sup _{r\in ]0,r_0]}\sup _{\mathcal{T}-r^2\le t\le\mathcal{T}}\frac{1}{r}\int_{B_{x_0,r}}|\vu(t,x)|^2dx\le \mathfrak{M}<+\infty,
\end{equation*}
and we obtain the following control on the initial data $\vec{u}^\gamma_0$
\begin{equation*}
\sup_{x_0\in \R}\int_{B_{x_0,1}}|\vec{u}^\gamma_0(y)|^2dy\le \mathfrak{M}<+\infty.
\end{equation*}
Let us study now the external force $\mathbb{W}^\gamma$ defined in (\ref{ExternalForceW}). Since $\omega_i^\gamma(\cdot,\cdot)=\gamma^2\omega_i(t_0+\gamma^2\cdot,\gamma\cdot)$, one has
\begin{align*}
\sup_{x_0\in \R} \int_0^{1}\int _{B_{x_0,1}}|\mathbb{W}^\gamma|^2dyds&=\sup_{x_0\in \R} \int_0^{1}\int _{B_{x_0,1}}\sum_{i=1}^32|\omega_i^\gamma|^2dyds\\&=\sup_{x_0\in \R} \int_0^{1}\int _{B_{x_0,1}}\sum_{i=1}^32|\gamma^2\omega_i(t_0+\gamma^2s,\gamma y)|^2dyds.
\end{align*}
Moreover, by a change of variable, we have
\begin{align*}
\sup_{x_0\in \R} \int_0^{1}\int _{B_{x_0,1}}|\mathbb{W}^\gamma|^2dyds&=\sup_{x_0\in \R} \frac{1}{\gamma}\int_{t_0}^{t_0+\gamma^2}\int _{B_{\gamma x_0,\gamma}}\sum_{i=1}^32|\omega(s,y)|^2dyds,
\end{align*}
and since $\gamma<r_0$ by \eqref{IntervalGamma},  one has
\begin{align*}
\sup_{x_0\in \R} \int_0^{1}\int _{B_{x_0,1}}|\mathbb{W}^\gamma|^2dyds&\le  C\sup_{z\in \R}\sup_{t_0\le t \le t_0+\gamma^2}\frac{1}{\gamma}\gamma^2\int _{B_{z,\gamma}}|\vw(t,y)|^2dy \\
&\le   Cr_0\sup_{z\in \R}\sup_{t_0\le t\le t_0+\gamma^2}\int _{B_{z,\gamma}}|\vw(t,y)|^2dy .
\end{align*}
Then, since $\vw\in L^\infty_tL^2_x$, we obtain
\begin{equation*}
\sup_{x_0\in \R} \int_0^{1}\int _{B_{x_0,1}}|\mathbb{W}^\gamma|^2dyds\le C\|\vw\|_{L^\infty_tL^2_x}^2<+\infty,
\end{equation*}
and we obtain the uniform control \eqref{EstimateForce} on the external force.\\

Thus, from the  estimate \eqref{EstimateInitialData} and the previous control, we can apply the Lemma  \ref{Theo_LocalLeray}, and therefore for $T=\min\{1,\frac{1}{C(1+M)^4}\}$ with $M=M(\mathfrak{M},\|\vw\|^2_{L^\infty_tL^2_x})>0$,  there exists  a local Leray solution $(\vv,q)$ of the  system \eqref{NSrotvW} (recall $\div({\mathbb{W}^\gamma})=\rot\vw^\gamma$) on $]0,T[\times \R$, such that for some constant $M_1>0$,
\begin{equation}\label{Hypo2}
\sup_{0<t<T}\sup_{x\in \R}\int_{B_{x,1}}|\vv(t,x)|^2dy+\sup_{x\in \R}\int_{0}^{T}\int_{B_{x,1}}|\grad\otimes\vv|^2dy\le M_1.
\end{equation}
It is worth noting that since $(\vec{u}^\gamma,p^\gamma)$ is  a strong solution of \eqref{NSrotvW}  by a weak-strong uniqueness argument (see \cite[Theorem 14.7]{PGLR1}), it follows that $\vec{u}^\gamma=\vv$ and $p^\gamma=q$  on $]0,\min\{\mathfrak{S},T\}[\times \R$. Thus, instead of studying $(\vec{u}^\gamma,p^\gamma)$ on $]0,\mathfrak{S}[\times \R$, we will apply Lemma \ref{Lem_localSmoothing2} and Corollary \ref{Coro_localSmothing} to the pair $(\vv,q)$ on $]0,T[\times\R$ and later by fixing $\mathfrak{S}\ll1$ small enough (in order to obtain the uniqueness on the intervall $]0, \mathfrak{S}[$)  we can come back to  the variables $(\vec{u}^\gamma,p^\gamma)$.\\

Let us verify that the triplet $(\vv,q,\vw^\gamma)$ satisfies the hypotheses of Lemma \ref{Lem_localSmoothing2}. Notice that since $(\vv,q)$ is already a local Leray solution by construction and we have  the control \eqref{Hypo2}, we only  need to verify  the following points:
\begin{itemize}
\item We have
\begin{equation}\label{InforotvW}
\|\vw^{\gamma}\|_{L^\infty(]0,T[,L^2( B_{0,1}))}^2< C\mathfrak{S}.
\end{equation}
\item For all $\mathfrak{r}>0$ such that $\sqrt{\mathfrak{S}}<\mathfrak{r}<1$, we have 
\begin{equation}\label{Hypo3}
{\displaystyle{\sup_{\sqrt{\mathfrak{S}}<\mathfrak{r}\le 1}\frac{1}{\mathfrak{r}}\int_{B_{0,\mathfrak{r}}}|\vec{u}^\gamma_0|^2dy}}< \varepsilon_\ast.
\end{equation} 

\end{itemize}
For obtaining \eqref{InforotvW}, since $T<1$ by construction and ${\vw^{\gamma}(\cdot,\cdot)=\gamma^2 \vw(t_0+\gamma^2 \cdot,\gamma \cdot)}$, by a change of variable we have
\begin{eqnarray*}
\|\vw^{\gamma}\|_{L^\infty(]0,T[,L^2( B_{0,1}))}^2&\le&\|\vw^{\gamma}\|_{L^\infty(]0,1[,L^2( B_{0,1}))}^2\\
&\leq &\|\gamma^2 \vw(t_0+\gamma^2 \cdot,\gamma \cdot)\|_{L^\infty(]0,1[,L^2( B_{0,1}))}^2=
\gamma\|\vw\|_{L^\infty(]t_0,t_0+\gamma^2[,L^2( B_{0,\gamma}))}^2.
\end{eqnarray*}
Since $\vw\in L^\infty_t L^2_x\cap L^2_t\dot H^1_x$ by hypothesis, we obtain 
$\|\vw^{\gamma}\|_{L^\infty(]0,1[,L^2( B_{0,1}))}^2\le \gamma\|\vw\|_{L^\infty_tL^2_x}^2\le C\gamma$. 
Moreover, since $\mathcal{T}-\delta^2<t_0<\mathcal{T}$ we have $\gamma=\sqrt{\frac{\mathcal{T}-t_0}{\mathfrak{S}}}< \sqrt{\frac{\delta^2}{\mathfrak{S}}}< r_0 \mathfrak{S}$ (recall $\delta=\min\{\sqrt{\frac{\mathcal{T}}{2}}, \frac{r_0\mathfrak{S}^\frac32}{2}\}$), and we can write
\begin{equation*}
\|\vw^{\gamma}\|_{L^\infty(]0,\mathfrak S[,L^2( B_{0,1}))}^2< C\mathfrak{S},
\end{equation*}
which is the wished estimate \eqref{InforotvW}. As pointed out in Remark \ref{Rem_SmallNessOmega}, the smallness condition to the ``external force'' is obtained by a suitable rescaling over the variable $\vw$.\\

\noindent For obtaining \eqref{Hypo3}, fix $\sqrt{\mathfrak{S}}<\mathfrak{r}<1$. Since $\vec{u}^\gamma(0,\cdot)=\gamma\vu(t_0,\gamma \cdot)$ and by  the change of variable $z=\gamma y$, we have
\begin{align*}
\int_{B_{0,\mathfrak{r}}}|\vec{u}^\gamma(0,y)|^2dy&=\int_{B_{0,\mathfrak{r}}}|\gamma\vu(t_0,\gamma y)|^2dy= \frac{1}{\gamma}\int_{B_{0,\gamma \mathfrak r}}|\vu(t_0,z)|^2dz.
\end{align*}
Moreover, since $\mathfrak{r}<1$ and since $\gamma=\sqrt{\frac{\mathcal T -t_0 }{\mathfrak{S}}}$ it follows from the assumption \eqref{controlSingularPoint} that
\begin{align*}
\int_{B_{0,\mathfrak{r}}}|\vec{u}^\gamma(0,y)|^2dy\le\frac{\sqrt{\mathfrak{S}}}{\sqrt{\mathcal T-t_0}}\int_{B_{0,\sqrt{\frac{\mathcal{T}-t_0}{\mathfrak{S}}}}}|\vu(t,x)|^2dx < \sqrt{\mathfrak{S}}\varepsilon_\ast,
\end{align*}
Therefore, since $\sqrt{\mathfrak{S}}<\mathfrak{r}$ one has $\displaystyle{\frac{1}{\mathfrak{r}}\int_{B_{0,\mathfrak{r}}}|\vec{u}^\gamma(0,y)|^2dy< \varepsilon_\ast}$. Since $ \sqrt{\mathfrak{S}}\le\mathfrak{r}<1$ is arbitrary, we conclude that 
\begin{equation*}
{\displaystyle{\sup_{\sqrt{\mathfrak{S}}<\mathfrak{r}\le 1}\frac{1}{\mathfrak{r}}\int_{B_{0,\mathfrak{r}}}|\vec{u}^\gamma_0|^2dy}}< \varepsilon_\ast.
\end{equation*}
Thus, we have proved that $(\vv,q,\vw^{\gamma})$ verifies the conditions \eqref{Hypo2}, \eqref{InforotvW} and \eqref{Hypo3}. Then, we can apply Lemma \ref{Lem_localSmoothing2} and Corollary \ref{Coro_localSmothing} and therefore there exists $\mathcal{T}_1=\mathcal{T}_1(T,\mathfrak{M})>0$, and a constant $\mathfrak{c}>1$  such that for any  $\sqrt{\mathfrak{S}}<r$ with $\lambda r^2\le \frac{\mathcal{T}_1}{\mathfrak{c}}$ and  $\lambda=\frac{1}{1+\varepsilon_\ast^2}$, we have
\begin{equation*}
\frac{1}{r^2} \int_0^{\lambda r^2}\int_{B_{0,r}}|\vv(s,y)|^3+|q-\mathfrak{h}_q|^\frac{3}{2}dyds \le C(\varepsilon_\ast+\varepsilon_\ast^\frac32).
\end{equation*}
Now, fix $\mathfrak{S}\ll1$ such that for $r=\sqrt{\frac{\mathfrak{S}}{\lambda}}$, we have $\sqrt{\mathfrak{S}}<r$ and $\lambda r^2< \frac{\mathcal{T}_1}{\mathfrak{c}}$. Therefore, by choosing $r=\sqrt{\frac{\mathfrak{S}}{\lambda}}$ in the expression above, we obtain
\begin{equation*}
\frac{\mathfrak{S}}{\lambda} \int_{0}^{\mathfrak{S}}\int_{B_{0,\sqrt{\frac{\mathfrak{S}}{\lambda}}}}|\vv|^3+|q-\mathfrak{h_q}|^\frac{3}{2}dyds \le C(\varepsilon_\ast+\varepsilon_\ast^\frac32).
\end{equation*}
Now, since  $\vec{u}^\gamma=\vv$ and $p^\gamma=q$  on $]0,\mathfrak{S}[\times \R$, we have
\begin{equation*}
\frac{\mathfrak{S}}{\lambda} \int_{0}^{\mathfrak{S}}\int_{B_{0,\sqrt{\frac{\mathfrak{S}}{\lambda}}}}|\vu^\gamma|^3+|p^\gamma-\mathfrak{h}^\gamma|^\frac{3}{2}dyds \le C(\varepsilon_\ast+\varepsilon_\ast^\frac32).
\end{equation*} \\
Thus, since $\vec{u}^\gamma(s,y)=\gamma \vu(t_0+\gamma^2 s,\gamma y),\;\;
p^{\gamma}(s,y)=\gamma^2 p(t_0+\gamma^2 s,\gamma y),$ we obtain (recall that $\gamma=\sqrt{\frac{\mathcal T -t_0 }{\mathfrak{S}}}$)
\begin{equation*}
\frac{\lambda}{\mathcal{T}-t_0} \int_{t_0}^{\mathcal{T}}\int_{B_{0,\sqrt{\frac{\mathcal{T}-t_0}{\lambda}}}}|\vu|^3+|p-\mathfrak{h}|^\frac{3}{2}dyds \le C(\varepsilon_\ast+\varepsilon_\ast^\frac32).
\end{equation*}
Since $\lambda<1$, we have $B_{0,\sqrt{\mathcal{T}-t_0}} \subset B_{0,\sqrt{\frac{\mathcal{T}-t_0}{\lambda}}}$ and we can write
\begin{equation*}
\frac{\lambda}{\mathcal{T}-t_0} \int_{\mathcal{T}-(\sqrt{\mathcal{T}-t_0})^2}^{\mathcal{T}}\int_{B_{0,\sqrt{\mathcal{T}-t_0}}}|\vu|^3+|p-\mathfrak{h}|^\frac{3}{2}dyds \le C(\varepsilon_\ast+\varepsilon_\ast^\frac32),
\end{equation*}
which can be rewritten as follows
\begin{equation*}
\frac{1}{(\sqrt{\mathcal{T}-t_0})^2} \int_{\mathcal{T}-(\sqrt{\mathcal{T}-t_0})^2}^{\mathcal{T}}\int_{B_{0,\sqrt{\mathcal{T}-t_0}}}|\vu|^3+|p-\mathfrak{h}|^\frac{3}{2}dyds \le C\frac{(\varepsilon_\ast+\varepsilon_\ast^\frac32)}{\lambda}.
\end{equation*}
Thus, since $\lambda=\frac{1}{1+\varepsilon_\ast^2}$, and $0<\varepsilon_\ast\ll1$ can be considered small enough, we  can find $0<\tilde \varepsilon\ll1$ such that $C(\varepsilon_\ast+\varepsilon_\ast^\frac32)\frac1\lambda<\tilde \varepsilon,$ and for $\rho=\sqrt{\mathcal{T}-t_0}$,  we obtain
\begin{equation*}
\frac{1}{ \rho^2} \int_{\mathcal{T}- \rho^2}^{\mathcal{T}}\int_{B_{0, \rho}}|\vu(s,y)|^3+|p-\mathfrak{h}|^\frac{3}{2}dyds <\tilde \varepsilon.
\end{equation*}
Now, since $(\vu,p-\mathfrak{h}, \vw)$ is also a partial suitable solution of the micropolar fluids equations, we can apply Theorem \ref{Theo_partialRegularity} in the appendix \ref{Sec_PartialRegularity} and it follows that $({\mathcal{T}},0)$ is a partial regular point in the sense of Definition \ref{Def_PartialRegularPoints}, which is a contradiction to the fact that $({\mathcal{T}},0)$ is partially singular  by hypothesis.\\

We thus have proved that there exists $\varepsilon_\ast>0$, $\mathfrak{S}=\mathfrak{S}(\mathfrak{M})$ and $\delta>0$ such that for all $t\in ]\mathcal{T}-\delta,\mathcal{T}[$, we have
\begin{equation*}
\frac{1}{\sqrt{{\mathcal{T}}-t}}\int_{B_{0,\sqrt{\frac{{\mathcal{T}}-t}{\mathfrak{S}}}}}|\vu(t,x)|^2dx\ge \varepsilon_\ast.
\end{equation*}
which as we mentioned before, implies the $L^3$-norm concentration effect of the velocity around the singular point $(\mathcal{T},0)$. This finishes the proof of Theorem \ref{Theo_BlowUpConcentration}. \hfill $\blacksquare$
\begin{Remark}\label{RemarkInterplay}
In all the previous computations, it might seem that the variable $\vw$ is only considered as an external force via the term $\rot \vw $ in the first equation (\ref{MicropolarFluidsEquationsEqua1}), however a very detailed study of its properties is essential to perform all the arguments given above. This fact will appear clearly in the pages below. 
\end{Remark}
\appendix
%%%%%%%%%%%%%%%%%%%%%%%%%%%%%%%%%%%%%%%%%%%%%%%%%%%
%%%%%%%%%%%%%%%%%%%%%%%%%%%%%%%%%%%%%%%%%%%%%%%%%%%
\mysection{A Serrin criterion for the micropolar fluid equations} \label{Sec_LocalregularityMP}
As pointed out in the page \pageref{ReferenceLocaleSerrin} of the introduction, in this appendix we establish a \emph{partial Serrin regularity criterion} for the micropolar fluids equations. The main idea consists in  deducing a gain of regularity for \emph{both} variables $\vu$ and $\vw$ by assuming only the local boundedness of $\vu$.\\

Since we are interested in the local behavior of a weak solution $(\vu, p,\vw)$ of the system \eqref{MicropolarFluidsEquationsEqua1} and \eqref{MicropolarFluidsEquationsEqua2} around a point $(t_0,x_0)\in ]0,+\infty[\times \R$, we will examine its regularity within the parabolic ball $Q_R(t_0,x_0)$ defined in (\ref{Def_Bolas2}) for some fixed $0<R^2<t_0$. We thus have:
%%%%%%%%%%%%%%%%%%%%%%%%%%%%%%%%%%%%%%%%%%%%%%%%%%%
\begin{Theorem}\label{Theo_SerrinMP}
Let $(\vu,p,\vw)$ be a weak solution of the micropolar fluids equations \eqref{MicropolarFluidsEquationsEqua1} and \eqref{MicropolarFluidsEquationsEqua2} over the parabolic ball $Q_R(t_0,x_0)$ given in \eqref{Def_Bolas2} such that
\begin{eqnarray*}
\vu,\vw \in L^\infty_tL_x^2 (Q_{R}(t_0,x_0))\cap L^2_t\dot H^1_x (Q_{R}(t_0,x_0))\quad\text{and }\quad p\in \mathcal{D}'_{t,x}(Q_R(t_0,x_0)).
\end{eqnarray*}
If we assume moreover that $\vu \in L^\infty_{t,x}(Q_R(t_0,x_0))$, then for all $0<r<R$, and for all $k\in \mathbb{N}$ we have 
$$\vu,\vw\in L^\infty(]t_0-r^2,t_0[,\dot H^k (B_{x_0,r}))\cap L^2(]t_0-r^2,t_0[,\dot H^{k+1}(B_{x_0,r})).$$
\end{Theorem}
%%%%%%%%%%%%%%%%%%%%%%%%%%%%%%%%%%%%%%%%%%%%%%%%%%%
\noindent Let us mention here that in order to obtain the wished gain of regularity, it will necessary to establish a dialogue between the variables $\vu$ and $\vw$ as it was pointed out in the Remark \ref{RemarkInterplay} above. Indeed, we will see first how to obtain a small gain of regularity for $\vu$, which will depend of the information we have over $\vw$. Then we will transfer this new information from $\vu$ to $\vw$, which in turn will imply a new gain of regularity of the velocity. Hence, by iterating this process we can obtain the whished conclusion.\\
 
%%%%%%%%%%%%%%%%%%%%%%%%%%%%%%%%%%%%%%%%%%%%%%%%%%%
\noindent\textbf{Proof of Theorem \ref{Theo_SerrinMP}.}
Let us study the regularity of $\vu$. For this, recall that  this variable satisfies  the  equation \eqref{MicropolarFluidsEquationsEqua1} \textit{i.e.,} we have
$$\partial_t\vu = \Delta \vu-(\vu\cdot \vn)\vu-\grad p +\frac{1}{2}\rot \vw.$$
Notice that the system above may be seen as the forced Navier-Stokes system where the external force is given by the term $\rot \vw $ which belongs to $L^2_{t}L^2_x(Q_R(t_0,x_0))$ (recall that by hypothesis we have ${\vw\in  L^2_t\dot H^1_x (Q_{R}(t_0,x_0))}$). Thus, since $\vu$ is bounded on $Q_R(t_0,x_0)$ by  hypothesis, we can apply the Serrin criterion of the Navier-Stokes equations (see for instance \cite[Theorem 13.1, pg 397]{PGLR1}), and therefore for some $0<r_1<R$, we have
\begin{equation}\label{RegulaU}
\vu \in L^\infty( ]t_0-r_1^2,t_0[,\dot H^1 (B_{x_0,r_1})\cap L^2(]t_0-r_1^2,t_0[,\dot H^2 (B_{x_0,r_1})).
\end{equation}
It is worth noting that we have obtained  a gain of regularity in the space variable for the velocity $\vu$, however we cannot expect any further information since the regularity of $\vu$ is linked to the external force represented here by the term $\frac{1}{2}\rot \vw $ and therefore we need to improve  the regularity of $\vw$.\\

Thus, let us prove  now that we can obtain the same gain of regularity given in \eqref{RegulaU} for $\vw$. For this,  we need some technical lemmas. First, we will recall a previous result given in \cite{ChLl22} that gives us an explicit gain of integrability for both variables $\vu$ and $\vw$ as long as $\vu$ belongs to some parabolic Morrey spaces. Next, we will establish that the divergence of $\vw$ belongs to $L^6_{t,x}(Q_r(t_0,x_0))$ for some $0<r<R$ and then we will show that the variable $\vw$ is bounded within $Q_{r}(t_0,x_0)$. Finally, with these informations at hand we will be able to deduce the gain of regularity for $\vw$ with respect to the space variable by considering the usual smoothing effects of the heat kernel.

\begin{Proposition}\label{Theo_MorreyTransfer}
Let $(\vu,p,\vw)$ be a weak solution of the micropolar fluids equations \eqref{MicropolarFluidsEquationsEqua1} and \eqref{MicropolarFluidsEquationsEqua2} over the parabolic ball ${\bf Q}_R(t_0,x_0)$ given in the expression \eqref{Def_Bolas1} such that we have the usual information ${\vu,\vw \in L^\infty( ]t_0-R^2,t_0+R^2[,L^2 (B_{x_0,R}))\cap L^2(]t_0-R^2,t_0+R^2[,\dot H^1 (B_{x_0,R}))}$ and $p\in \mathcal{D}'_{t,x}({\bf Q}_R(t_0,x_0))$.
If moreover we have the following local hypothesis
\begin{equation*}
\mathds{1}_{{\bf Q}_R(t_0,x_0)}\vu \in \M_{t,x}^{p_0,q_0}(\mathbb{R}\times\R)\quad with\;\;
2<p_0\le q_0,\; 5<q_0\le 6,
\end{equation*}
then
\begin{itemize}
\item[1)] for a parabolic ball ${\bf Q}_{\tau_1}(t_0,x_0)\subset {\bf Q}_R(t_0,x_0)$ we have
\begin{equation*}
\vu \in L_{t,x}^{q_0}({\bf Q}_{\tau_1}(t_0,x_0)),\qquad 5<q_0\le 6.
\end{equation*}
\item[2)] For a parabolic ball ${\bf Q}_{\tau_2}(t_0,x_0)\subset {\bf Q}_{\tau_1}(t_0,x_0) $, we have 
\begin{equation*}
\vw \in L_{t,x}^{q_0}({\bf Q}_{\tau_2}(t_0,x_0)),\qquad 5<q_0\le 6.
\end{equation*}
\end{itemize}
\end{Proposition}
%%%%%%%%%%%%%%%%%%%%%%%%%%%%%%%%%%%%%%%%%%%%%%%%%%%
\begin{Remark}\label{Rem_Morrey}
It is important to note that the previous proposition essentially extends O'Leary's result in \cite{OLeary} for the Navier-Stokes equations to the context of the micropolar fluid equations. In the statement above, the balls ${\bf Q}_r$ of the type \eqref{Def_Bolas1} are considered, however, the proof given in \cite{ChLl22} can be easily adapted to the balls $Q_r$ given in \eqref{Def_Bolas2} where only the time variable lies in a slightly different interval. With no loss of generality, we will frequently replace the balls ${\bf Q}_r$ by $Q_r$ in the sequel. See also \cite{Kukavica}, \cite{OLeary} and \cite[Theorem 13.3]{PGLR1} for a similar treatment. 
\end{Remark}
\noindent With this result at hand we now study the integrability of $\div(\vw)$ within the parabolic ball $Q_R(t_0,x_0)$.
%%%%%%%%%%%%%%%%%%%%%%%%%%%%%%%%%%%%%%%%%%%%%%%%%%%
\begin{Lemma}\label{Lem_integrabilityDivergence}
Under the hypotheses of Theorem \ref{Theo_SerrinMP}, for all ${Q_{r_1}(t_0,x_0)\subset Q_R(t_0,x_0)}$ we have ${\mathds{1}_{Q_{r_1}(t_0,x_0)}\div(\vw)\in L^{6}_{t,x}}(\mathbb{R}\times\mathbb{R}^3)$.
\end{Lemma} 
%%%%%%%%%%%%%%%%%%%%%%%%%%%%%%%%%%%%%%%%%%%%%%%%%%%
\noindent \textbf{Proof.} First, notice that since $\vu$ is bounded over the set $Q_R(t_0,x_0)$ by hypothesis, we obviously obtain that ${\mathds{1}_{Q_{R}}\vu\in\M_{t,x}^{3,6}(\mathbb{R}\times \R)}$. Therefore, by Remark \ref{Rem_Morrey} we can consider Proposition \ref{Theo_MorreyTransfer} over balls of the type $Q_r$, and thus there exists $0<r_0<R$ such that 
\begin{equation}\label{InforNomr6qOmega}
{\mathds{1}_{Q_{r_0}}\vw\in L^{6}_{t,x}(\mathbb{R}\times \R)}.
\end{equation}
Now, with this additional information over the variable $\vw$, we may study the local integrability of $\div(\vw)$. Let $\varphi: \mathbb{R}\times \R \longrightarrow \mathbb{R}$ be a test function such that for $0<r_1<r_0<R$,
\begin{align*}
\varphi &\equiv 1\;\; \text{on} \;\;]t_0-r_1^2,t_0+r_1^2[\times B_{x_0,r_1}\quad \text{and} \quad \supp(\varphi)\subset ]t_0-r_0^2,t_0+r_0^2[\times B_{x_0,r_0}.
\end{align*}
Given that we are interested in the local information of $\div(\vw)$, we set ${\mathbb{W}= \varphi \div(\vw)}$. By applying formally the divergence operator to the  equation \eqref{MicropolarFluidsEquationsEqua2} we obtain the following: 
\begin{equation}\label{equaDiv(w)}
\partial_t \div(\vw) = 2\Delta \div(\vw) -\div(\vw) -\div(\div(\vw\otimes \vu)).
\end{equation}
Moreover we easily deduce that (recall that we have ${\mathbb{W}= \varphi \div(\vw)}$):
$$\partial_t \mathbb{W}-2\Delta \mathbb{W}=(\partial_t \varphi -2\Delta \varphi) \div(\vw)+4\sum_{i=1}^3\partial_i((\partial_i \varphi) \div(\vw))+\varphi (\partial_t \div(\vw)-2\Delta \div(\vw)).$$
Hence, we get for any $t\in [0,t_0[$,
\begin{equation*}
\begin{cases}
\partial_t \mathbb{W}=2\Delta \mathbb{W}+(\partial_t\varphi-2\Delta \varphi-\varphi )\div(\vw)+4\displaystyle{\sum_{i=1}^3}\partial_i((\partial_i \varphi) \div(\vw))-\varphi \div(\div(\vw\otimes \vu)),\\[3mm]
\mathbb{W}(0, \cdot)=0.
\end{cases}
\end{equation*}
Thus, by the Duhamel's formula we obtain
\begin{eqnarray}\label{funcW}
\mathbb{W}(t,x)&=&\underbrace{\int_0^t e^{2(t-s)\Delta}\bigg((\partial_t\varphi-2\Delta \varphi-\varphi )\div(\vw)\bigg)ds}_{(I_\mathbb{W})}+4\sum_{i=1}^3\underbrace{\int_0^t e^{2(t-s)\Delta}\bigg(\partial_i((\partial_i \varphi) \div(\vw))ds\bigg)}_{(II_\mathbb{W})}\notag\\
&&-\underbrace{\int_0^t e^{2(t-s)\Delta}\bigg(\varphi \div(\div(\vw\otimes \vu))\bigg)ds}_{(III_\mathbb{W})}.
\end{eqnarray}
We shall prove that each term of the right-hand side of the expression above belongs to $L_{t,x}^{6}(]0,t_0[\times \R)$.
\begin{itemize}
\item[$\bullet$]For the first term $(I_\mathbb{W})$ in \eqref{funcW}, by setting $\psi= \partial_t\varphi-2\Delta \varphi-\varphi $ we can write 
\begin{align}
(I_\mathbb{W})&= \int_0^t e^{2(t-s)\Delta}\bigg((\partial_t\varphi-2\Delta \varphi-\varphi )\div(\vw)\bigg)ds =\sum_{i=1}^3\int_0^t e^{2(t-s)\Delta}\psi(\partial_i \omega _i) ds \notag\\
&=\sum_{i=1}^3 \int_0^t (\partial _ i e^{2(t-s)\Delta})(\psi \omega_i)ds-\int_0^t e^{2(t-s)\Delta} \big((\partial _ i\psi) \omega_i\big) ds. \label{FirstTermFuncW}
\end{align}
Now, since ${\psi \in \mathcal{C}_{0}^\infty(\mathbb{R}\times \R)}$ and $\supp(\psi)\subset  ]t_0-r_0^2,t_0+r_0^2[\times B_{x_0,r_0}$, for the first term of the right-hand side of \eqref{FirstTermFuncW} we obtain
\begin{eqnarray*}
\left\|\mathds{1}_{\{0<t<t_0\}}\int_0^t(\partial_i e^{2(t-s)\Delta}) (\psi \omega_i)(s,\cdot)ds \right\|_{L^{6}}&\leq & \mathds{1}_{\{0<t<t_0\}}\int_0^t\|\partial _i\mathfrak{g}_{2(t-s)}\|_{L^1_x}\|\psi \omega_i(s,\cdot)\|_{L^{6}}ds\\
&\leq &C\mathds{1}_{\{0<t<t_0\}}\int_0^{t} (t-s)^{-\frac{1}{2}} \left\|\mathds{1}_{Q_{r_0}}\psi w_i(s,\cdot) \right\|_{L^{6}}ds,
\end{eqnarray*}
where we have used the Young inequality for the convolution (recall that the action of the operator $e^{2(t-s)\Delta}$ is given by a convolution with the heat kernel $\mathfrak{g}_{2(t-s)}$) and the usual $L^p$-estimates of the heat kernel. 
Thus, by the Hölder inequality in the time variable with $1=\frac{1}{6}+\frac{5}{6}$, we obtain 
\begin{eqnarray*}
\left\|\mathds{1}_{\{0<t<t_0\}}\int_0^t(\partial _i e^{2(t-s)\Delta}) (\psi \omega_i)(s,\cdot)ds \right\|_{L^{6}}&\le &C\|\psi\|_{L^\infty_{t,x}} \left\|\mathds{1}_{Q_{r_0}}\vw\right\|_{L^{6}_{t,x}}\left(\int_0^{t} (t-s)^{-\frac{6}{10}}ds\right)^{\frac{5}{6}}\\
&\leq &C\|\mathds{1}_{Q_{r_0}}\vw\|_{L^{6}_{t,x}}.
\end{eqnarray*}
Hence, by taking the $L^6$-norm in the time variable and since ${\mathds{1}_{Q_{r_0}}\vw\in L^{6}_{t,x}(\mathbb{R}\times \R)}$ by \eqref{InforNomr6qOmega},  we have
\begin{equation}\label{termI.1}
\left\|\mathds{1}_{\{0<t<t_0\}}\int_0^t ( \partial _j e^{2(t-s)\Delta}) (\psi w_i)(s,\cdot)ds \right\|_{L^{6}_{t,x}}\le C\left\|\mathds{1}_{Q_{r_0}} \vw\right\|_{L^{6}_{t,x}}<+\infty.
\end{equation}
For the second term of \eqref{FirstTermFuncW}, again by the properties of the test function $\psi$ and the Young inequality for the convolution we obtain
\begin{eqnarray*}
\left\|\mathds{1}_{\{0<t<t_0\}}\int_0^te^{2(t-s)\Delta} \big((\partial _i\psi) \omega_i\big)(s,\cdot)ds \right\|_{L^{6}}&\leq & \mathds{1}_{\{0<t<t_0\}}\int_0^t \|\mathfrak{g}_{2(t-s)}\|_{L^1} \|(\partial _i\psi) \omega_i(s,\cdot)\|_{L^{6}}ds \\
&\le& C\int_0^{t_0}\|\mathds{1}_{Q_{r_0}} (\partial _i \psi)\vw(s,\cdot) \|_{L^{6}}ds.
\end{eqnarray*}
Hence, by taking the $L^6$-norm as well as the Hölder inequality in the time variable and from \eqref{InforNomr6qOmega} we obtain 
\begin{equation}\label{termI.2}
\left\|\mathds{1}_{\{0<t<t_0\}}\int_0^te^{2(t-s)\Delta}\big( (\partial _i\psi) \omega_i\big)(s,\cdot)ds \right\|_{L^{6}_{t,x}}\le C \|\partial _i \psi\|_{L^\infty_{t,x}} \left\|\mathds{1}_{Q_{r_0}}\vw \right\|_{L^{6}_{t,x}}< +\infty.
\end{equation}
Therefore from \eqref{FirstTermFuncW}, \eqref{termI.1} and the estimate above, we conclude that the term $(I_\mathbb{W})$ in \eqref{funcW} belongs to $L^{6}_{t,x}(]0,t_0[\times \R)$.

\item[$\bullet$] For the term $(II_\mathbb{W})$ in \eqref{funcW}, it is enough to study the following expression
\begin{equation}\label{termII}
\int_0^t \partial _i e^{2(t-s)\Delta}\big((\partial_i \varphi )\partial_j \omega_j\big)ds=
\int_0^t\partial _j \partial _i e^{2(t-s)\Delta} \big((\partial_i \varphi) \omega_j\big)ds-
\int_0^t\partial _i e^{2(t-s)\Delta} \big((\partial _j \partial_i \varphi) \omega_j\big)ds,
\end{equation}
for all $1\le i,j\le 3$. Thus, by the maximal regularity of the heat kernel (see \cite[Theorem 7.3]{PGLR0}) and by the support properties of the function $\varphi$, we have
\begin{eqnarray*}
\left\|\mathds{1}_{\{0<t<t_0\}}\int_0^t \partial_i \partial _j e^{2(t-s)\Delta} \big((\partial_i \varphi)
\omega_j(s,\cdot)\big)ds \right\|_{L_{t,x}^{6}}&\le& C \left\| \mathds{1}_{Q_{r_0}}(\partial _i \varphi)\vw \right\|_{L_{t,x}^{6}}\\
&\le& C \|\partial _i \varphi\|_{L^\infty_{t,x}}\left\| \mathds{1}_{Q_{r_0}}\vw \right\|_{L_{t,x}^{6}}<+\infty.
\end{eqnarray*}
Since the second term of the right-hand side of \eqref{termII} can be treated in a similar fashion as \eqref{termI.1}, by replacing $\psi$ for $\partial _j \partial_i \varphi$, 
we can conclude that the term $(II_\mathbb{W})$ in \eqref{funcW} belongs to $L^{6}_{t,x}(]0,t_0[\times \R)$.

\item[$\bullet$] For the third term $(III_\mathbb{W})$ of \eqref{funcW}, notice that for all $1\le i,j\le 3$, we have
\begin{eqnarray}
\int_0^t e^{2(t-s)\Delta}\big(\varphi \partial_i(\partial_j(\omega_ju_i))\big)ds
=\int_0^t \partial_j\partial_i e^{2(t-s)\Delta}\big(\varphi\omega_j u_i\big)ds-
\int_0^t \partial_i e^{2(t-s)\Delta}\big((\partial_j\varphi) \omega_ju_i\big)ds\nonumber\\
-\int_0^t \partial_j e^{2(t-s)\Delta}\big((\partial_i\varphi)\omega_ju_i\big)ds
+\int_0^te^{2(t-s)\Delta}\big((\partial_i\partial_j\varphi)\omega_ju_i\big)ds\label{div(uw)}.
\end{eqnarray}
For the first term of the expression above by using the maximal regularity of the heat kernel, the hypothesis $\mathds{1}_{Q_R}\vu \in L^\infty_{t,x}$ and \eqref{InforNomr6qOmega}, we can establish that
\begin{align*}
\left\|\mathds{1}_{\{0<t<t_0\}}\int_0^t \partial_j\partial_i e^{2(t-s)\Delta}(\varphi \omega_ju_i)ds\right\|_{L^{6}_{t,x}}&\le 
\|\mathds{1}_{Q_{r_0}}(\varphi\omega_ju_i)\|_{{L^{6}_{t,x}}}\\
&\le C\| \varphi\|_{L^\infty_{t,x}}\|\mathds{1}_{Q_{R}}\vu\|_{L^{\infty}_{t,x}}
\|\mathds{1}_{Q_{r_0}}\vw\|_{L^{6}_{t,x}}<+\infty.
\end{align*}
Since the second and third terms in the right hand-side of \eqref{div(uw)} share the same structure, 
it is enough to study only one of them. Thus, by the same arguments as in \eqref{termI.1}, we have
\begin{align*}
\left\|\mathds{1}_{\{0<t<t_0\}}\int_0^t\partial _j e^{2(t-s)\Delta} \big((\partial_j\varphi) \omega_ju_i\big)(s,\cdot)ds \right\|_{L^{6}_{t,x}}
&\le\left\|\mathds{1}_{Q_{r_0}}(\partial_j\varphi) \omega_ju_i\right\|_{L^{6}_{t,x}}\\
&\le C\|\partial_j\varphi\|_{L^\infty_{t,x}} \left\|\mathds{1}_{Q_{R}}\vu \right\|_{L^{\infty}_{t,x}} \left\|\mathds{1}_{Q_{r_0}} \vu\right\|_{L^{6}_{t,x}}<+\infty.
\end{align*} 
Finally the last term of \eqref{div(uw)}, we can use the estimate \eqref{termI.2} and
we obtain 
\begin{align*}
\left\|\int_0^te^{2(t-s)\Delta}\big((\partial_j\partial_i\varphi) \omega_ju_i\big)ds\right\|_{L^{6}_{t,x}}
&\le \left\|\mathds{1}_{Q_{r_0}}(\partial_j\partial_i\varphi) \omega_ju_i\right\|_{L^{6}_{t,x}}\\
&\le C\|\partial_j\partial_i\varphi\|_{L^\infty_{t,x}} \left\|\mathds{1}_{Q_{R}}\vu \right\|_{L^{\infty}_{t,x}} \left\| \vu\right\|_{L^{6}_{t,x}}<+\infty.
\end{align*}
Consequently, we find that $(III_\mathbb{W})$ in \eqref{funcW} belongs to $L^{6}_{t,x}(]0,t_0[\times \R)$. 
\end{itemize}

\noindent Hence, we have proved that the quantities $(I_\mathbb{W}), (II_\mathbb{W}), (III_\mathbb{W})$ given in (\ref{funcW}) belong to $L^6_{t,x}(]0,t_0[\times \R)$ and therefore we obtain that the function $\mathbb{W} \in L^6_{t,x}(]0,t_0[\times \R)$. By using the properties of the test function (recall $\varphi=1$ on $ ]t_0-r_1^2,t_0+r_1^2[\times B_{x_0,r_1}$) we finally conclude that
\begin{equation*}
\|\mathds{1}_{Q_{r_1}(t_0,x_0)}\div(\vw)\|_{L^6_{t,x}}<+\infty,
\end{equation*}
and this finishes the proof of the Lemma \ref{Lem_integrabilityDivergence}.\hfill $\blacksquare $\\

Having obtained this gain of information over $\div(\vw)$, we can now deduce the boundedness of $\vw$.
%%%%%%%%%%%%%%%%%%%%%%%%%%%%%%%%%%%%%%%%%%%%%%%%%%%
\begin{Proposition}\label{Pro_LinfinitymicropolarVitesse}
Let $(\vu,p,\vw)$ be a weak solution of the micropolar fluids equations \eqref{MicropolarFluidsEquationsEqua1} and \eqref{MicropolarFluidsEquationsEqua2} such that ${\vu,\vw \in L^\infty( ]t_0-R^2,t_0[,L^2 (B_{x_0,R}))\cap L^2(]t_0-R^2,t_0[,\dot H^1 (B_{x_0,R}))}$ and $p\in \mathcal{D}'_{t,x}(Q_R(t_0,x_0))$. If ${\vu \in L^{\infty}_{t,x}(Q_R)}$ then $\vw\in L^\infty_{t,x}(Q_{r_2})$ for all $Q_{r_2}\subset Q_R$.
\end{Proposition}
%%%%%%%%%%%%%%%%%%%%%%%%%%%%%%%%%%%%%%%%%%%%%%%%%%
\noindent \textbf{Proof.} Notice that from Remark \ref{Rem_Morrey}, Proposition \ref{Theo_MorreyTransfer} and Lemma \ref{Lem_integrabilityDivergence}, we can establish the existence of $r_0$ and $r_1$ such that  ${0<r_1<r_0<R}$ and such that
\begin{equation}\label{InformationL6}
\mathds{1}_{Q_{r_0}}\vw\in L^{6}_{t,x}\quad  \text{and} \quad{\mathds{1}_{Q_{r_1}}\div(\vw)\in L^{6}_{t,x}}.
\end{equation}  
We consider now a positive test function $ \phi: \mathbb{R}\times \R \longrightarrow \mathbb{R}$ such that for some radius $r_2$ we have $0<r_2< r_1<r_0 <R$ and 
\begin{align*}
\phi &\equiv 1\;\; \text{on} \;\;]t_0-r_2^2,t_0+r_2^2[\times B_{x_0,r_2}\quad \text{and} \quad \supp(\phi)\subset ]t_0-r_1^2,t_0+r_1^2[\times B_{x_0,r_1}.
\end{align*}
Define $\vec{\mathcal{W}}  =\phi \vw$. By using the dynamics of the variable $\vw$, \emph{i.e.} the equation  \eqref{MicropolarFluidsEquationsEqua2}, we have for any $t\in [0,t_0[,$
$$
\begin{cases}
\partial_t \vec{\mathcal{W}}=\Delta \vec{\mathcal{W}}+(\partial_t \phi-\Delta \phi)\vw+2\displaystyle{\sum_{i=1}^3}\partial_i ((\partial_i \phi)\vw)+\phi\left[\grad \div(\vw)-\vw -\div(\vw \otimes\vu )+\frac{1}{2}\rot\vu\right] ,\\
\vec{\mathcal{W}}(0, \cdot)=0.
\end{cases}
$$
By Duhamel's formula we obtain
\begin{eqnarray}\label{funcZ}
{\vec{\mathcal{W}}}(t,x)&=&\underbrace{\int_0^te^{(t-s)\Delta}(\partial_t\phi-\Delta \phi -\phi)\vw ds}_{(I_{\vec{\mathcal{W}}})}+2\sum_{i=1}^3\underbrace{\int_0^te^{(t-s)\Delta}\partial_i ((\partial_i \phi) \vw)ds}_{(II_{\vec{\mathcal{W}}})}+\underbrace{\int_0^te^{(t-s)\Delta}\phi\grad \div(\vw)ds}_{(III_{\vec{\mathcal{W}}})} \notag\\
&&-\underbrace{\int_0^te^{(t-s)\Delta}\phi \div(\vw \otimes\vu )ds}_{(IV_{\vec{\mathcal{W}}})}+\underbrace{\int_0^te^{(t-s)\Delta}\frac{\phi}{2}\rot\vu ds}_{(V_{\vec{\mathcal{W}}})}ds.
\end{eqnarray}
Thus, we shall prove that every term in the right-hand side of \eqref{funcZ} is bounded on $[0,t_0[\times \R$. We  study each term above separately.
\begin{itemize}

\item For the term $(I_{\vec{\mathcal{W}}})$ in \eqref{funcZ} by setting $\Phi= \partial_t\phi-\Delta \phi -\phi $ we can write 
$$(I_{\vec{\mathcal{W}}})= \int_0^t e^{(t-s)\Delta}(\partial_t\phi-\Delta \phi -\phi)\vw ds=\int_0^t e^{(t-s)\Delta}\Phi\vw ds.$$
Note that ${\Phi \in \mathcal{C}_{0}^\infty(\mathbb{R}\times \R)}$ and $\supp(\Phi)\subset  ]t_0-r_1^2,t_0+r_1^2[\times B_{x_0,r_1}$.  Recall that $e^{(t-s)\Delta}$ is given by a convolution operator with the heat kernel $\mathfrak{g}_{(t-s)}$, hence by the Young inequality for the convolution and the $L^p$-estimates of the heat kernel, we have
\begin{align*}
\left\|\mathds{1}_{\{0<t<t_0\}}\int_0^t e^{(t-s)\Delta} \Phi \vw(s,\cdot) ds \right\|_{L^\infty}
&\le C \mathds{1}_{\{0<t<t_0\}} \int_0^t \|\mathfrak{g}_{(t-s)}\|_{L^\frac{6}{5}} \|\Phi \vw(s,\cdot)\|_{L^{6}}ds \\
&\le  C\mathds{1}_{\{0<t<t_0\}}\int_0^{t} (t-s)^{-\frac{1}{4}} \| \mathds{1}_{Q_{r_0}}\Phi\vw(s,\cdot)\|_{L^{6}}ds.
\end{align*}
Moreover, by applying the H\"older inequality in the time variable  with $1=\frac{1}{6}+\frac{5}{6}$ we have 
\begin{align*}
\left\|\mathds{1}_{\{0<t<t_0\}}\int_0^t e^{(t-s)\Delta} \Phi \vw(s,\cdot) ds \right\|_{L^\infty}&\le C \mathds{1}_{\{0<t<t_0\}}
\| \mathds{1}_{Q_{r_0}}\Phi\vw\|_{L^{6}_{t,x}}\left(\int_0^{t} (t-s)^{-\frac{3}{10}}ds\right)^{\frac{5}{6}}\\
&\le C \mathds{1}_{\{0<t<t_0\}}  \| \Phi\|_{L^\infty_{t,x}}
\| \mathds{1}_{Q_{r_0}}\vw\|_{L^{6}_{t,x}} t ^{\frac{7}{12}}. 
\end{align*}
Therefore, by taking the supremum in the time variable and from \eqref{InformationL6}, one has
\begin{equation}\label{heatw}
\left\|\mathds{1}_{\{0<t<t_0\}}\int_0^t e^{(t-s)\Delta} \Phi \vw ds \right\|_{L_{t,x}^\infty}\le 
C\| \Phi\|_{L^\infty_{t,x}}\|\mathds{1}_{Q_{r_0}}\vw\|_{L^{6}_{t,x}}<+\infty.
\end{equation}
\item For the term $(II_{\vec{\mathcal{W}}})$ in \eqref{funcZ}, the Young inequality for the convolution and the $L^p-$ estimates of the heat kernel imply the following control for all $1\le i\le 3$
\begin{align*}
\left\|\mathds{1}_{\{0<t<t_0\}}\int_0^t \partial_i e^{(t-s)\Delta} (\partial_i\phi) \vw (s,\cdot)ds \right\|_{L^\infty}
&\le C \mathds{1}_{\{0<t<t_0\}}\int_0^{t}\|\partial_i\mathfrak{g}_{(t-s)}\|_{L^\frac{6}{5}}\| \mathds{1}_{Q_{r_0}} (\partial_i\phi)\vw(s,\cdot)\|_{L^6}ds\\
&\le C \mathds{1}_{\{0<t<t_0\}}\int_0^{t}(t-s)^{-\frac{3}{4}} \| \mathds{1}_{Q_{r_0}} (\partial_i\phi)\vw(s,\cdot)\|_{L^6}ds
\end{align*}
Thus, by the H\"older inequality in the time variable, we have
\begin{align*}
\left\|\mathds{1}_{\{0<t<t_0\}}\int_0^t \partial_i e^{(t-s)\Delta} (\partial_i\phi) \vw (s,\cdot)ds \right\|_{L^\infty}
&\le  C \mathds{1}_{\{0<t<t_0\}} \| \mathds{1}_{Q_{r_0}} \partial_i\phi \vw\|_{L^{6}_{t,x}} \left(\int_0^{t} (t-s)^{-\frac{9}{10}}ds\right)^{\frac{5}{6}}\\
&\le  C \mathds{1}_{\{0<t<t_0\}} \| \partial_i\phi\|_{L^\infty_{t,x}}\| \mathds{1}_{Q_{r_0}} \vw\|_{L^{6}_{t,x}} t^{\frac{1}{12}}.
\end{align*}
Therefore, by considering the supremum in time and from \eqref{InformationL6} we obtain
\begin{equation}\label{derivativeheatw}
\left\|\mathds{1}_{\{0<t<t_0\}}\int_0^t \partial_i e^{(t-s)\Delta}( \partial_i\phi )\vw ds \right\|_{L_{t,x}^\infty}\le C \| \partial_i\phi\|_{L^\infty_{t,x}} \|\mathds{1}_{Q_{r_0}}\vw\|_{L^{6}_{t,x}}<+\infty.
\end{equation}

\item Now, we consider the term $(III_{\vec{\mathcal{W}}})$ in \eqref{funcZ}. For this, we write
\begin{equation*}
\int_0^t e^{(t-s)\Delta} \phi \grad \div(\vw) ds=-\int_0^t e^{(t-s)\Delta} (\grad \phi) \div(\vw) ds+\int_0^t \grad e^{(t-s)\Delta} \phi \div(\vw) ds.
\end{equation*}
Therefore, since ${\mathds{1}_{Q_{ r_1}}\div(\vw)\in L^{6}_{t,x}}$, we can apply the same arguments as in \eqref{heatw} and \eqref{derivativeheatw} to obtain
\begin{align*}
\left\|\mathds{1}_{\{0<t<t_0\}}\int_0^t e^{(t-s)\Delta} (\grad\phi) \div(\vw) ds \right\|_{L_{t,x}^\infty}
&\le C\| \grad \phi\|_{L^\infty_{t,x}} \|\mathds{1}_{Q_{r_1}}\div(\vw)\|_{L^{6}_{t,x}}<+\infty,\\
\left\|\mathds{1}_{\{0<t<t_0\}}\int_0^t \grad e^{(t-s)\Delta} \phi \div(\vw) ds \right\|_{L_{t,x}^\infty}
&\le C\| \phi\|_{L^\infty_{t,x}}\| \mathds{1}_{Q_{r_1}}\div(\vw)\|_{L^{6}_{t,x}}<+\infty.
\end{align*}
\item Finally, we turn our attention to the terms $(IV_{\vec{\mathcal{W}}})$ and $(V_{\vec{\mathcal{W}}})$ in \eqref{funcZ}, both of which involve the presence of the velocity $\vu$. Let us begin with the term $(IV_{\vec{\mathcal{W}}})$, for which is enough to study the expression
\begin{equation}\label{termuW}
\int_0^t e^{(t-s)\Delta}\phi \partial_j(\omega_i u_j )ds=
\int_0^t \partial_j e^{(t-s)\Delta}\phi \omega_i u_j ds-\int_0^t e^{(t-s)\Delta}(\partial_j\phi )\omega_i u_j ds,
\end{equation}
for all $1\le i, j\le 3$. For the first term of the right-hand side in the expression above, by using the same arguments as in \eqref{derivativeheatw} we obtain 
\begin{align*}
\left\|\mathds{1}_{\{0<t<t_0\}}\int_0^t \partial_j e^{(t-s)\Delta} \phi \omega_i u_jds \right\|_{L_{t,x}^\infty}
&\le C\| \phi\|_{L^\infty_{t,x}}\| \mathds{1}_{Q_{r_0}} \omega_i u_j \|_{L^{6}_{t,x}}\\
&\le C\| \phi\|_{L^\infty_{t,x}}\|\mathds{1}_{Q_{R}}\vu\|_{L^{\infty}_{t,x}}
\|\mathds{1}_{Q_{r_0}}\vw\|_{L^{6}_{t,x}}<+\infty.
\end{align*}
In addition for the second term of \eqref{termuW}, by the same arguments of \eqref{heatw} we have 
\begin{align*}
\left\|\mathds{1}_{\{0<t<t_0\}}\int_0^te^{(t-s)\Delta} (\partial_j\phi) \omega_i u_j ds \right\|_{L_{t,x}^\infty}&\le C\| \partial_j\phi\|_{L^\infty_{t,x}}\| \mathds{1}_{Q_{r_0}}\omega_i u_j \|_{L^{6}_{t,x}}\\
& \le C\|\mathds{1}_{Q_R}\vu\|_{L^{\infty}_{t,x}}
\|\mathds{1}_{Q_{r_0}} \vw\|_{L^{6}_{t,x}}<+\infty.
\end{align*}
The last term $(V_{\vec{\mathcal{W}}})$ of \eqref{funcZ} follows easily. Indeed, since  we can write
\begin{equation*}
\int_0^te^{(t-s)\Delta} \phi\rot\vu ds=-\int_0^te^{(t-s)\Delta} (\grad\phi)\wedge\vu ds+\int_0^t \rot e^{(t-s)\Delta} \phi\vu ds,
\end{equation*}
and since $\mathds{1}_{Q_R}\vu\in L^6_{t,x}(\mathbb{R}\times \R)$, we can treat each  term above by using the same arguments as in \eqref{heatw} and \eqref{derivativeheatw} respectively. 
\end{itemize}
Therefore by the previous points, the terms $(I_{\vec{\mathcal{W}}})$-$(V_{\vec{\mathcal{W}}})$ given in \eqref{funcW} are bounded, and thus ${\vec{\mathcal{W}}\in L^\infty_{t,x}([0,t_0[\times \R)}$. Hence, by using the properties of the test function ($\phi =1$ on $]t_0-r_2^2,t_0+r_2^2[\times B_{x_0,r_2}$),  we obtain that $\vw$ is bounded on $Q_{r_2}(t_0,x_0)$ and this finish the proof of Proposition \ref{Pro_LinfinitymicropolarVitesse}. \hfill $\blacksquare $\\

We now state the last technical result. Here, we will prove  that whenever $\vu$ is more regular than $\vw$  we can transfer this information to $\vw$ in smaller balls.
%%%%%%%%%%%%%%%%%%%%%%%%%%%%%%%%%%%%%%%%%%%%%%%%%%%
\begin{Proposition}\label{Propo_gainRegularityW}
Under the general hypothesis of Theorem \ref{Theo_SerrinMP} if we assume that
\begin{align*}
\vu\in L^\infty( ]t_0-r_1^2,t_0[,\dot H^{1} (B_{x_0,r_1}))\cap L^2(]t_0-r_1^2,t_0[,\dot H^{2} (B_{x_0,r_1})),
\end{align*}
then for some radius $r_3$ such that $0<r_3<r_2<r_1<R$, we have
\begin{equation*}
\vw \in L^\infty( ]t_0-r_3^2,t_0[,\dot H^{1} (B_{x_0,r_3}))\cap L^2(]t_0-r_3^2,t_0[,\dot H^{2} (B_{x_0,r_3})).
\end{equation*}
\end{Proposition}
%%%%%%%%%%%%%%%%%%%%%%%%%%%%%%%%%%%%%%%%%%%%%%%%%%%
\noindent \textbf{Proof.}  First, notice that from Proposition \ref{Pro_LinfinitymicropolarVitesse}, we have for  ${0<r_2<r_1<R}$,
\begin{equation}\label{InformationLInfinty}
\mathds{1}_{Q_{r_2}}\vw\in L^{\infty}_{t,x}.
\end{equation}  
Let $ \phi: \mathbb{R}\times \R \longrightarrow \mathbb{R}$ be a test function such that for $0<r_3<\mathfrak{r}<r_2<r_1<R$, 
\begin{align*}
\phi&\equiv 1\;\; \text{on} \;\;]t_0-r_3^2,t_0+r_3^2[\times B_{x_0,r_3}\quad \text{and} \quad \supp(\phi)\subset ]t_0-\mathfrak{r}^2,t_0+\mathfrak{r}^2[\times B_{x_0,\mathfrak{r}}.
\end{align*}
By using the equality $\Delta (\phi \omega_i)= \Delta \phi \omega_i +2\div(\grad \phi\omega_i)- \phi \Delta \omega_i$, we have for any $1\le i\le 3$, 
\begin{equation}\label{TermToStudy}
\phi \omega_i=\frac{1}{(-\Delta)}\left[- \Delta \phi \omega_i -2\div(\grad \phi \omega_i)+ \phi \Delta \omega_i\right].
\end{equation}
Thus, in order to improve the regularity of $\vw$, we may prove that the expression above belongs to ${L^\infty_t\dot H^{1}_x\cap L^2_t\dot H_x^{2}}$. For this,  we will deduce a gain of information for the laplacian of $\vw$, and later we will study the regularity of $\vw$. 
%%%%%%%%%%%%%%%%%%%%%%%%%%%%%%%%%%%%%%%%%%%%%%%%%%%
\begin{itemize}
\item[$\ast$] \textit{\bf A local gain of information for the laplacian of $\vw$.} 
By considering the identity
\begin{equation}\label{InfoLaplacian}
\Delta \vw =\grad \div(\vw)-\rot (\rot \vw),
\end{equation}
it is clear that we can obtain information for the laplacian of $\vw$ from the its divergence and its curl. Thus, let $\psi: \mathbb{R}\times \R \longrightarrow \mathbb{R}$ be a test function such that for $0<r_3<\mathfrak{r}<r_2<R$,  $$\psi \equiv 1 \quad \text{on}\;\;]t_0-\mathfrak{r}^2,t_0+\mathfrak{r}^2[\times B_{x_0,\mathfrak{r}}\quad \text{and} \quad \supp(\psi)\subset ]t_0-r_2^2,t_0+r_2^2[\times B_{x_0,r_2}.$$
Define now $\vec{\mathfrak{W}}=\psi\rot \vw$ and $\mathbb{W}=\psi \div(\vw)$. Note that the dynamics of these variables are straightforward to compute, indeed, by taking formally the curl operator to the equation \eqref{MicropolarFluidsEquationsEqua2}, we obtain $\partial_t\rot \vw=\Delta \rot \vw -\rot((\vu\cdot\grad)\vw)-\rot\vw+\frac{1}{2}\rot \rot \vu$ and recall the dynamics of  $\div(\vw)$ was already obtained in \eqref{equaDiv(w)}. Hence,  from these equations we can deduce
{\small
\begin{align}\label{equaCurlAndDiv}
\begin{aligned}
&\partial_t \vec{\mathfrak{W}}=\Delta \vec{\mathfrak{W}}+\underbrace{(\partial_t\psi-\Delta \psi -\psi)\rot\vw}_{(I_a)}+2\displaystyle{\sum_{i=1}^3}\underbrace{\partial_i((\partial_i \psi) \rot\vw)}_{(II_a)}-\underbrace{\psi \rot ((\vu\cdot\grad )\vw)}_{(IIII_a)}+\underbrace{\frac{1}{2}\psi\rot( \rot \vu)}_{(IV_a)},\\[1mm]
&\partial_t \mathbb{W}=2\Delta \mathbb{W}+\underbrace{(\partial_t\psi-2\Delta \psi -\psi)\div(\vw)}_{(I_\alpha)}+4\displaystyle{\sum_{i=1}^3}\underbrace{\partial_i((\partial_i \psi) \div(\vw))}_{(II_\alpha)}-\underbrace{\psi \div((\vu\cdot\grad )\vw)}_{(III_\alpha)},
\end{aligned}
\end{align}
}\noindent such that $\vec{\mathfrak{W}}(0,\cdot)=\mathbb{W}(0,\cdot)=0$ due to the properties of the test function.\\

\noindent We claim now that the each one of the term of right-hand side of \eqref{equaCurlAndDiv} belongs to $L^2_t\dot H^{-1}_x$. Indeed we have the following points:
%%%%%%%%%%%%%%%%%%%%%%%%%%%%%%%%%%%%%%%%%%%%%%%%%%%
\begin{itemize}
\item [$\bullet$] First, we consider the terms $(I_a)$ and $(I_\alpha)$. Note that they share the same structure, therefore we study only the first one. Hence, we have
\begin{align*}
\|(\partial_t\psi-\Delta \psi -\psi)\rot \vw(t,\cdot)\|_{\dot H^{-1}}
&\le \|(\partial_t\psi-\Delta \psi-\psi )\rot \vw(t,\cdot)\|_{L^\frac{6}{5}}\\
&\le C\|\partial_t\psi-\Delta \psi-\psi(t,\cdot)\|_{L^3}\|\rot \vw(t,\cdot)\|_{L^2(B_{r_2})},
\end{align*} 
where we have used the embedding $L^{\frac{6}{5}}(\R)\subset \dot H^{-1}(\R)$ and the Hölder inequality ($\frac{5}{6}=\frac{1}{3}+\frac{1}{2}$).
Moreover, by taking the $L^2$-norm in the time variable we obtain
\begin{align*}
\|(\partial_t\psi-\Delta \psi -\psi)\rot \vw\|_{L^2_t\dot H_x^{-1}}
\le C\|\partial_t\psi-\Delta \psi-\psi\|_{L^\infty_tL^3_x}\|\rot \vw\|_{L^2_{t}L^2_x(Q_{R})}<+\infty.
\end{align*}

\item [$\bullet$] For the terms $(II_a)$ and $(II_\alpha)$ in \eqref{equaCurlAndDiv}, again since they share the same structure, it is enough to study only the first one. Notice that for any $1\le i\le 3$,
\begin{equation*}
\|\partial_i(\partial_i \psi \rot \vw)\|_{L^2_t\dot H_x^{-1}}\le 
\|\partial_i \psi \rot \vw\|_{L^2_{t}L^2_x}\le\|\partial_i \psi\|_{L^\infty_{t,x}} \| \rot \vw\|_{L^2_{t}L^2_x(Q_{R})}<+\infty,
\end{equation*}
since $\vw\in L^2_t\dot H^1_x(Q_R)$.\\
\item [$\bullet$]For the terms $(III_a)$ and $(III_\alpha)$ in \eqref{equaCurlAndDiv}, since we can write
$(\vu\cdot\grad )\vw=\div(\vw\otimes \vu)$, it is enough to study the following expression for any $1\le i,m,j\le 3$
\begin{eqnarray}\label{DescompositionBilineal}
\psi \partial_i(\partial _j (\omega_m u_j ))
&=&\partial_i\partial _j(\psi\omega_m u_j )-\partial_j((\partial_i \psi) \omega_m u_j )\notag\\
&&-\partial_i((\partial_j \psi) \omega_m u_j )+(\partial_j \partial_i\psi)(\omega_m u_j ).
\end{eqnarray}
By taking the $L^2_t\dot H^{-1}_x$ norm in the expression above we obtain
\begin{eqnarray}\label{TermsBilineal}
\|\psi \partial_i(\partial _j (\omega_m u_j ))\|_{L^2_t\dot H_x^{-1}}
&\le&\underbrace{\|\partial_i\partial _j(\psi\omega_m u_j )\|_{L^2_t\dot H_x^{-1}}}_{(1)}+\underbrace{\|\partial_j((\partial_i \psi)( \omega_m u_j ))\|_{L^2_t\dot H_x^{-1}}}_{(2)}\\\notag
&&+\underbrace{\|\partial_i((\partial_j \psi)( \omega_m u_j ))\|_{L^2_t\dot H_x^{-1}}}_{(3)}+\underbrace{\|(\partial_j \partial_i\psi)(\omega_m u_j )\|_{L^2_t\dot H_x^{-1}}}_{(4)}.
\end{eqnarray}
For the term $(1)$ in the expression above, since  $\vu$ is bounded on $Q_R(t_0,x_0)$  by hypothesis, $\vw$ is bounded on $Q_{r_2}(t_0,x_0)$ by \eqref{InformationLInfinty} and since $\supp(\psi)\subset ]t_0-r_2^2,t_0+r_2^2[\times B_{x_0,r_2} $, we easily observe
\begin{eqnarray}
\|\partial_i\partial _j(\psi\omega_m u_j )\|_{L^2_t\dot H_x^{-1}}&\le& \|\psi\omega_m u_j \|_{L^2_t\dot H_x^{1}}=\sum_{|\alpha|=1}\|D^\alpha(\psi\omega_m u_j )\|_{L^2_t L_x^{2}}\notag\\
&=&\sum_{|\alpha|=1}\|(D^\alpha \psi)\omega_m u_j+ \psi(D^\alpha\omega_m) u_j+\psi\omega_m(D^\alpha u_j) \|_{L^2_t L_x^{2}}\notag\\
&\le&C\sum_{|\alpha|=1}\|(D^\alpha \psi)\omega_m\|_{L_{t,x}^{\infty}}\|u_j \|_{L_{t,x}^2(Q_{R})}+ C   \|\psi u_j\|_{L_{t,x}^{\infty}}\|\omega_m\|_{L^2_t\dot H_x^{1}(Q_{R})}\notag\\
&&+\|\psi \omega_m\|_{L_{t,x}^{\infty}}\|u_j\|_{L^2_t\dot H_x^{1}(Q_{R})}<+\infty.\label{termIBilineal}
\end{eqnarray}
For the terms $(2)$ and $(3)$ in \eqref{TermsBilineal} since they have the same structure, we only study  the first one. Thus, again by using the boundedness of $\vw$ on $Q_{r_2}(t_0,x_0)$, we have
\begin{eqnarray}
\|\partial_j((\partial_j \psi)( \omega_m u_j ))\|_{L^2_t\dot H_x^{-1}}&\le& \| (\partial_j \psi) \omega_m u_j \|_{L^2_t L_x^{2}}\notag\\
&\le& C\| (\partial_j \psi)\omega_m\|_{L_{t,x}^{\infty}}\| u_j \|_{L^\infty_tL_x^{2}(Q_{R})}<+\infty.\label{termIIBilineal}
\end{eqnarray}
For the last term $(4)$ of \eqref{TermsBilineal}, by the embedding $L^\frac{6}{5}(\R)\subset \dot H^{-1}(\R)$, and the H\"older inequality ($\frac{5}{6}=\frac{1}{2}+\frac{1}{3}$), one has
\begin{eqnarray}
\|(\partial_j \partial_i\psi)\omega_m u_j \|_{L^2_t\dot H_x^{-1}}&\le &\|(\partial _j\partial_i \psi) \omega_m u_j \|_{L^2_t L_x^{\frac{6}{5}}}\notag\\
&\le &C\| u_j \|_{L^\infty_tL_x^{2}(Q_{R})}\|(\partial_j \partial_i\psi) \omega_m\|_{L_{t}^{2}L_x^3}\notag\\
&\le& C\| u_j \|_{L^\infty_tL_x^{2}(Q_{R})}\|\vw\|_{L^\infty_{t,x}(Q_{r_2})}\|\partial_j \partial_i\psi\|_{L_{t}^{2}L_x^3}<+\infty.\label{termIVBilineal}
\end{eqnarray}
Hence, from \eqref{termIBilineal}-\eqref{termIIBilineal}-\eqref{termIVBilineal}, using the expression (\ref{DescompositionBilineal}) we can see that the terms $(III_a)$ and $(III_\alpha)$ in \eqref{equaCurlAndDiv} belong to $L^2_t\dot H^{-1}_x$ (recall that by Proposition \ref{Pro_LinfinitymicropolarVitesse} we have $\|\vw\|_{L^\infty_{t,x}(Q_{r_2})}<+\infty$).\\
\item [$\bullet $]For the last term $(IV_a)$ of \eqref{equaCurlAndDiv} we have
\begin{align*}
\psi \rot (\rot \vu) &=\rot(\psi \rot \vu) -(\grad \psi )\wedge (\rot \vu).
\end{align*}
Thus, by taking the $\dot H^{-1}(\R)$ norm in space variable we obtain
\begin{align*}
\|\psi \rot (\rot \vu)(t,\cdot) \|_{\dot H^{-1}}&\le \|\rot(\psi \rot \vu)(t,\cdot)\|_{\dot H^{-1}}+\|(\grad \psi )\wedge (\rot \vu)(t,\cdot)\|_{\dot H^{-1}}\\
&\le \|\psi \rot \vu(t,\cdot)\|_{L^2}+\|(\grad \psi )\wedge (\rot \vu)(t,\cdot)\|_{L^\frac{6}{5}}
\end{align*}
where we have used the embedding $L^\frac{6}{5}(\R)\subset \dot H^{-1}(\R)$. Now, by integrating in time and by the Hölder inequality in space with $\frac{5}{6}=\frac{1}{3}+\frac{1}{2}$ we conclude
\begin{eqnarray*}
\|\psi \rot (\rot \vu)\|_{L^2_t\dot H_x^{-1}}&\le &C \|\psi\|_{L^\infty_{t,x}} \| \rot \vu\|_{L_{t,x}^2(Q_{R})}\notag\\
&& + C\|\grad \psi\|_{L^2_{t}L^3_x}\|\rot\vu\|_{L^2_{t}L^2_x(Q_{R})}<+\infty.
\end{eqnarray*}
\end{itemize}
%%%%%%%%%%%%%%%%%%%%%%%%%%%%%%%%%%%%%%%%%%%%%%%%%%%
Therefore, from the previous points, we have proven that each term of the right-hand side of \eqref{equaCurlAndDiv} belong to $L^2([0,t_0[,\dot H^{-1}(\R))$. Thus, by the theory developed in \cite[Section 13, page 398]{PGLR1} (which is essentially the Serrin regularity criterion for the Navier-Stokes equations) we have 
\begin{equation*}
\vec{\mathfrak{W}},  \mathbb{W}\in L^\infty([0,t_0[, L^2(\R))\cap L^2([0,t_0[,\dot H^{1}(\R)).
\end{equation*} 
Furthermore,  from the identity \eqref{InfoLaplacian}, we can deduce that
\begin{equation}\label{LaplacianGainInfo}
\phi \Delta \vw \in L^\infty([0,t_0[, \dot H^{-1}(\R))\cap L^2([0,t_0[,L^2(\R)).
\end{equation}
%%%%%%%%%%%%%%%%%%%%%%%%%%%%%%%%%%%%%%%%%%%%%%%%%%%
\item[$\ast$]\textit{\bf A gain of regularity in the space variable for $\vw$}. 
Recall that from \eqref{TermToStudy} we have for all $1\le i \le 3$
\begin{equation}\label{OmegaDescomposition}
\phi \omega_i=\frac{1}{(-\Delta)}\left[- \Delta \phi \omega_i -2\div(\grad \phi \omega_i)+ \phi\Delta \omega_i\right].
\end{equation}
Now, we will prove that each term in the expression above belongs to $L^\infty_t\dot H_x^{1}\cap L^2_t\dot H_x^{2} $.  Firs, by considering the $L^2_t\dot H^{2}_x$-norm in \eqref{OmegaDescomposition}, we have
\begin{equation*}
\|\phi \omega _i\|_{L^2_t\dot H_x^{2}}\le \|-\Delta \phi \omega_i -2\div(\grad \phi \omega_i)+ \phi \Delta \omega_i\|_{L^2_tL^2_x}.
\end{equation*}
Notice that 
\begin{equation*}
\|2\div(\grad \phi \omega_i)\|_{L^2_tL^2_x}\le2\|\Delta \phi \omega_i\|_{L^2_tL^2_x}+ 2\|\grad \phi \cdot \grad\omega_i\|_{L^2_tL^2_x},
\end{equation*}
and thus by the triangular inequality and since $\supp(\phi)\subset ]t_0-R^2,t_0+r^2[\times B_{x_0,R}$, it follows that
\begin{eqnarray}
\|\phi \omega _i\|_{L^2_t\dot H_x^{2}}&\le& C\| \phi \Delta \omega_i\|_{L^2_tL^2_x}+ C\|\Delta \phi \omega_i\|_{L^2_tL^2_x}+C\|\grad \phi \cdot \grad\omega_i\|_{L^2_tL^2_x}\notag\\
&\le & C\| \phi \Delta \omega_i\|_{L^2_tL^2_x}+ C\| \vw\|_{L^2_tL^2_x(Q_{R})}+C\| \grad\vw\|_{L^2_tL^2_x(Q_{R})}<+\infty,\label{L2H2norm}
\end{eqnarray}
where we have used  \eqref{LaplacianGainInfo}  and the fact that $\vw \in L^\infty_tL_x^2 (Q_{R}(t_0,x_0))\cap L^2_t\dot H^1_x (Q_{R}(t_0,x_0))$.\\

On the other hand, by  considering the $ \dot H^{1}(\R)$-norm in  \eqref{OmegaDescomposition}, we obtain
\begin{align}
\|\phi \omega _i(t,\cdot)\|_{\dot H^{1}}&\le \|\Delta \phi \omega_i (t,\cdot)+2\div(\grad \phi \omega_i)(t,\cdot)+ \phi \Delta \omega_i(t,\cdot)\|_{\dot H^{-1}}\notag\\
&\le 
\|\Delta \phi \omega_i (t,\cdot)\|_{\dot H^{-1}}+ \|\grad \phi \omega_i(t,\cdot)\|_{L^2}+\| \phi \Delta \omega_i(t,\cdot)\|_{\dot H^{-1}}.\label{estimateaux1}
\end{align}
From the embedding $\dot H^{-1}(\R)\subset L^\frac65(\R)$ and the Hölder inequality, we have
$$\|\Delta \phi \omega_i (t,\cdot)\|_{\dot H^{-1}}\le\|\Delta \phi \omega_i (t,\cdot)\|_{L^\frac65}\le \|\Delta \phi(t,\cdot) \|_{L^3}\|\omega_i(t,\cdot)\|_{L^2(B_{x_0,R})},$$
 and therefore  by taking the supremum in time in \eqref{estimateaux1} and by \eqref{InfoLaplacian}, one has
\begin{equation*}
\|\phi \omega _i\|_{L^\infty_t\dot H_x^{1}}\le C\|\omega_i\|_{L^\infty_tL^2_x(Q_{R})}+C\|\omega_i\|_{L^\infty_tL_x^2(Q_{R})}+C\|\phi \Delta \omega_i\|_{L^\infty_t\dot H_x^{-1}}<+\infty.
\end{equation*}
Therefore, from \eqref{L2H2norm} and the expression above we obtain $\phi \vw \in L^\infty_t\dot H_x^{1}\cap L^2_t\dot H_x^{2}$, \textit{i.e.},
\begin{equation*}
\vw \in L^\infty( ]t_0-r_3^2,t_0[,\dot H^{1} (B_{x_0,r_3}))\cap L^2(]t_0-r_3^2,t_0[,\dot H^{2} (B_{x_0,r_3})),
\end{equation*} 
and thus the proof of Proposition \ref{Propo_gainRegularityW} is finished.\hfill $\blacksquare$
\end{itemize}
%%%%%%%%%%%%%%%%%%%%%%%%%%%%%%%%%%%%%%%%%%%%%%%%%%%
\subsection*{End of the proof of Theorem \ref{Theo_SerrinMP}}
Recall that we have proved that for some $0<r_1<R$, we have
\begin{equation*}
\vu \in L^\infty( ]t_0-r_1^2,t_0[,\dot H^1 (B_{x_0,r_1})\cap L^2(]t_0-r_1^2,t_0[,\dot H^2 (B_{x_0,r_1})).
\end{equation*}
Thus, we can apply Proposition \ref{Propo_gainRegularityW} and therefore it follows that  for some $0<r_3<r_2<r_1<R$  we have
\begin{equation*}
\vw \in L^\infty( ]t_0-r_3^2,t_0[,\dot H^{1} (B_{x_0,r_3}))\cap L^2(]t_0-r_3^2,t_0[,\dot H^{2} (B_{x_0,r_3})).
\end{equation*}
In particular, since  $\rot\vw\in L^2_t\dot H^1_x(Q_{r_3})$, we can apply again the Serrin criterion for the Navier-Stokes equations to $\vu$ and therefore it follows that for some radius $r_4$ such that $0<r_4<r_3<R$, we have
\begin{equation*}
\vu \in L^\infty( ]t_0-r_4^2,t_0[,\dot H^2(B_{x_0,r_4})\cap L^2(]t_0-r_4^2,t_0[,\dot H^3 (B_{x_0,r_4})).
\end{equation*}
Thus, by following the same arguments given as in Proposition \ref{Propo_gainRegularityW}, we can improve as well the regularity of $\vw$ and since we can iterate this process, we obtain the wished regularity for $(\vu,\vw)$ and hence the proof of Theorem \ref{Theo_SerrinMP} is finished. 
\hfill $\blacksquare $

%%%%%%%%%%%%%%%%%%%%%%%%%%%%%%%%%%%%%%%%%%%%%%%%%%%
%%%%%%%%%%%%%%%%%%%%%%%%%%%%%%%%%%%%%%%%%%%%%%%%%%%
\mysection{Partial regularity theory for the micropolar equations} \label{Sec_PartialRegularity}
This section is devoted to the partial regularity theory of the micropolar equations \eqref{MicropolarFluidsEquationsEqua1} and \eqref{MicropolarFluidsEquationsEqua2}. Let us point out that in \cite{ChLl23}, this theory was already developed in the framework of the micropolar system when considering partial  suitable solutions in the sense of Definition \ref{Def_PartialSuitable}. Indeed, it was proven that if for some $0<\varepsilon\ll1$ we  have
\begin{equation*}
\limsup_{r\to 0}\int _{{\bf Q}_r(t_0,x_0)}| \grad\otimes\vu|^2dyds<\varepsilon,
\end{equation*}
then the variables $(\vu,\vw)$ are Hölder continuous in time and space around the point $(t_0,x_0)$. \\

Thus, following essentially the same ideas, we will prove in this section that the ``second'' criterion of the Caffarelli, Kohn and Nirenberg theory remains valid for the micropolar fluid equations \textit{i.e.,} we will deduce a gain of regularity for the variables $(\vu,\vw)$ when only assuming some conditions over the velocity $\vu$ and the pressure $p$. 
\begin{Theorem}\label{Theo_partialRegularity}
Let $(\vu,p, \vw)$ be a partial suitable solution in the sense of Definition \ref{Def_PartialSuitable} for the micropolar equations \eqref{MicropolarFluidsEquationsEqua1} and \eqref{MicropolarFluidsEquationsEqua2} in $Q_1(t_0,x_0)$. Assume there exists a constant $\varepsilon>0$ small enough such that for some $0<R^2<\min \{1,t_0\}$, we have
\begin{equation}\label{TeohyL3}
\frac{1}{R^2}\int_{t_0-R^2}^{t_0}\int _{B_{x_0,R}}| \vu|^3+|p|^\frac{3}{2}dx ds<\varepsilon.
\end{equation}
Then, there exists some $0<r<R$ such that $\vu,\vw\in L_{t,x}^{\infty}(Q_r(t_0,x_0))$.
\end{Theorem} 
%%%%%%%%%%%%%%%%%%%%%%%%%%%%%%%%%%%%%%%%%%%%%%%%%%%
\noindent Following the same ideas than \cite{ChLl23}, in order to prove the previous theorem, we will first deduce from the hypothesis \eqref{TeohyL3} a gain of Morrey information for the velocity $\vu$. Indeed, we have
%%%%%%%%%%%%%%%%%%%%%%%%%%%%%%%%%%%%%%%%%%%%%%%%%%%
\begin{Proposition}\label{Pro_partialRegualrity}
Let $(\vu,p, \vw)$ be a partial suitable solution in the sense of Definition \ref{Def_PartialSuitable}  of the micropolar system \eqref{MicropolarFluidsEquationsEqua1} and \eqref{MicropolarFluidsEquationsEqua2} over the parabolic ball $Q_1(t_0,x_0)$. Assume there exists a sufficiently small constant $\varepsilon>0$ such that for some $0<R^2\le\min\{1,t_0\}$, we have
\begin{equation*}
\frac{1}{R^2}\int_{Q_R(t_0,x_0)}| \vu|^3+|p|^\frac{3}{2}dx ds<\varepsilon.
\end{equation*}
Then, there exists a radius $0<\mathfrak{r}<R$ such that for any $\frac{5}{1-\alpha}<\tau_0\le \frac{15}{2}$ with $0<\alpha<\frac{1}{6}$, we have ${\mathds{1}_{Q_{\mathfrak{r}}(t_0,x_0)}\vu\in\M_{t,x}^{3,\tau_0}(\mathbb{R}\times \R)}$ and $\mathds{1}_{Q_{\mathfrak{r}}(t_0,x_0)}p\in \M_{t,x}^{\frac{3}{2},\frac{\tau_0}{2}}(\mathbb{R}\times \R)$.
\end{Proposition} 
%%%%%%%%%%%%%%%%%%%%%%%%%%%%%%%%%%%%%%%%%%%%%%%%%%%
\noindent {\bf Proof.} Our aim consists in proving that for some $0<\mathfrak{r}<R$,  we have for all  $0<r\le\mathfrak{r}$ and ${(t,x)\in Q_{\mathfrak{r}}(t_0,x_0)}$, 
\begin{equation}\label{GainMorrey}
\int_{{\bf Q}_r(t,x)}{\mathds{1}_{Q_{\mathfrak{r}}(t_0,x_0)}}|\vu|^3 dy ds \le C r^{5(1-\frac{3}{\tau_0})}\quad\text{and}\quad\int_{{\bf Q}_r(t,x)}{\mathds{1}_{Q_{\mathfrak{r}}(t_0,x_0)}}|p|^\frac32 dy ds \le C r^{5(1-\frac{3}{\tau_0})}
\end{equation}
which is the definition of Morrey spaces (see for instance \eqref{DefMorreyparabolico}). For this, we will consider the following  quantities: for a point $(t,x)\in \mathbb{R}\times \R$ and for $r>0$ we write
\begin{equation}\label{Def_Invariants}
\begin{split}
\mathcal{A}_r(t,x)&=\sup_{t-r^2<s<t}\frac{1}{r}\int_{B_{x_0,r}}|\vu(s,y)|^2dy, 
\qquad\qquad \alpha_r(t,x)=\frac{1}{r} \int _{Q_r(t,x)}|\grad\otimes \vu(s,y)|^2dyds,\\
\lambda_r(t,x)&=\frac{1}{r^2} \int _{Q_r(t,x)}|\vu(s,y)|^3dyds, 
\qquad\qquad \qquad\mathcal{P}_r(t,x)=\frac{1}{r^2} \int _{Q_r(t,x)}|p (s,y)|^{\frac{3}{2}}dyds.
\end{split}
\end{equation}
If $Q_r(t,x)\cap Q_R^c\neq \emptyset$ then the above quantities are replaced by $Q_r(t,x)\cap Q_R$. Moreover, for simplicity we introduce the following notations
\begin{equation}\label{Def_QuantiteIteration}
\Lambda_r=\frac{1}{r^{3(1-\frac{5}{\tau_0})}}\lambda_r,\quad \mathbb{P}_r=\frac{1}{r^{3(1-\frac{5}{\tau_0})}}\mathcal{P}_r\quad\mbox{and}\quad \mathbb{O}_r=\Lambda_r+\kappa^6\mathbb{P}_r.
\end{equation}
where $0<\kappa\ll1$ is a small fixed parameter to be defined later. Therefore it is easy to see that \eqref{GainMorrey} is equivalent to prove that for all $0<r\le\mathfrak{r}$ and for all $(t,x)\in Q_{r}(t_0,x_0)$, we have
\begin{equation}\label{IterativeO1}
\mathbb{O}_r(t,x)\le C.
\end{equation}
Now, in order to obtain \eqref{IterativeO1}, we will proceed by an iterative argument for which  we  need to introduce some technical  lemmas. First, let us point out the following relationship between the quantities given in \eqref{Def_Invariants}.
%%%%%%%%%%%%%%%%%%%%%%%%%%%%%%%%%%%%%%%%%%%%%%%%%%%
\begin{Lemma}\label{Estimation_normeLambda3}
For any $0<r\le R$, the quantities defined in \eqref{Def_Invariants} verify that
\begin{equation*}
\lambda_r^{\frac{1}{3}}\le C(\mathcal{A}_r+\alpha_r)^{\frac{1}{2}},
\end{equation*}
where $C$ is a constant that does not depend on $r$
\end{Lemma}
%%%%%%%%%%%%%%%%%%%%%%%%%%%%%%%%%%%%%%%%%%%%%%%%%%%
\noindent \textbf{Proof.} By using the definition of $ \lambda_r$ given in (\ref{Def_Invariants}) above and by the H\"older inequality $(\frac{1}{3}=\frac{1}{30}+\frac{3}{10})$ we have
$$\lambda_r^{\frac{1}{3}}=\frac{1}{r^{\frac{2}{3}}}\|\vu\|_{L_{t,x}^{3}(Q_r)}\le\frac{C}{r^\frac{2}{3}}r^{\frac5{30}}\|\vu\|_{L_{t,x}^{\frac{10}{3}}(Q_r)}=C\frac1{r^\frac{1}{2}}\|\vu\|_{L_{t,x}^{\frac{10}{3}}(Q_r)}.$$
Since by interpolation we have $\|\vu(t,\cdot)\|_{L^{\frac{10}{3}}(B_{x_0,r})}\leq \|\vu(t,\cdot)\|_{L^{2}(B_{x_0,r})}^{\frac{2}{5}}\|\vu(t,\cdot)\|_{L^{6}(B_{x_0,r})}^{\frac{3}{5}}$, we can easily deduce that $\|\vu\|_{L_{t,x}^{\frac{10}{3}} (Q_r)} \leq \|\vu\|_{L_t^{\infty}L_x^{2} (Q_r)}^{\frac25}\|\vu\|_{L_t^{2}L_x^{6} (Q_r)}^{\frac{3}{5}}$. Now, we use the classical Gagliardo-Nirenberg inequality (see \cite{Brezis}) to obtain $\|\vu\|_{L_t^{2}L_x^{6} (Q_r)} \leq C\big(\|\grad\otimes\vu\|_{L_t^{2}L_x^{2} (Q_r)} +\|\vu\|_{L_t^{\infty}L_x^{2} (Q_r)}\big)$ and by using Young's inequalities we have
\begin{eqnarray*}
\|\vu\|_{L_{t,x}^{\frac{10}{3}} (Q_r)} &\leq& C \|\vu\|_{L_t^{\infty}L_x^{2} (Q_r)}^{\frac 25}\big(\|\grad\otimes \vu\|_{L_t^{2}L_x^{2} (Q_r)}^{\frac 35}+\|\vu\|_{L_t^{\infty}L_x^{2} (Q_r)}^{\frac35} \big) \leq C\big(\|\vu\|_{L_t^{\infty}L_x^{2} (Q_r)}+\|\grad\otimes \vu\|_{L_t^{2}L_x^{2} (Q_r)}\big).
\end{eqnarray*}
By noting that $\|\vu\|_{L_t^{\infty}L_x^{2} (Q_r)}=r^{\frac12}\mathcal{A}_r^{\frac12}$ and $\|\grad\otimes \vu\|_{L^2_tL_x^{2} (Q_r)}=r^{\frac12}\alpha_r^{\frac12}$, we finally obtain the desired estimate. \hfill $\blacksquare $\\

We now present a first  estimate linked to the local energy inequality that allows us to control the terms in \eqref{Def_Invariants} within smaller balls.
%%%%%%%%%%%%%%%%%%%%%%%%%%%%%%%%%%%%%%%%%%%%%%%%%%%
\begin{Lemma}\label{lema_FirstEstimate}
Under the hypotheses of Proposition \ref{Pro_partialRegualrity}, for any radius $\displaystyle{0<r\le \tfrac\rho2\le R}$ we have the inequality
\begin{equation*}
\mathcal{A}_r+\alpha_r \le C\frac{r^2}{\rho ^2} \lambda_\rho ^\frac{2}{3}+C\frac{\rho ^2}{r^2}\left(\mathcal{P}_\rho+\lambda_\rho\right)+C\frac{\rho^{\frac{3}{2}}}{r}\lambda_\rho^{\frac{1}{3}}.
\end{equation*}
\end{Lemma}
%%%%%%%%%%%%%%%%%%%%%%%%%%%%%%%%%%%%%%%%%%%%%%%%%%%
\noindent \textbf{Proof.} The main idea for proving this lemma consists in plugging a well chosen test function in the local energy inequality which we recall in the following lines: for all ${\psi\in \mathcal{D}_{t,x}(\mathbb{R}\times \R)}$ 
\begin{eqnarray}
\int_{\R} |\vu|^2\psi(t,\cdot) dx+2 \int_{s<t} \int _{\R}|\grad\otimes \vu|^2\psi\, dx ds
\leq \int_{s<t}\int _{\R} (\partial_t\psi +\Delta \psi)|\vu|^2dyds+2\int_{s<t}\int_{\R}p (\vu\cdot \grad\psi )dx ds\notag\\
+\int_{s<t}\int _{\R} |\vu|^2(\vu \cdot \grad)\psi dx ds+\int_{s<t}\int_{\R}(\rot \vw)\cdot (\psi\vu)dyds.\qquad \qquad\label{Ineq_Energie1}
\end{eqnarray}
Regarding the test function to be chosen, we can mention Scheffer's work in \cite{She77} where it was introduce the following one: consider $\phi \in \mathcal{C}_0^{\infty}(\mathbb{R}\times \R)$  a test function such that
\begin{equation*}
\phi(s,y)=r^2\gamma\left(\frac{s-t_0}{\rho^2},\frac{y-x_0}{\rho}\right)
\theta\left(\frac{s-t_0}{r^2}\right)\mathfrak{g}_{(4r^2+t_0-s)}(x_0-y),
\end{equation*}
where $\gamma \in \mathcal{C}_0^{\infty}(\mathbb{R}\times \R)$ is positive function whose support is in ${ Q}_1(0,0)$ and equal to 1 in ${ Q}_{\frac{1}{2}}(0,0)$. In addition $\theta$ is a  non negative smooth function such that $\theta=1$ over $]-\infty,1[$ and $\theta=0$ over $]2,+\infty[$ and $\mathfrak{g}_t(\cdot)$ is the usual heat kernel. Then, we have the following points.
\begin{itemize}
\item[1)]the function $\phi$ is a bounded non-negative function, and its support is contained in the parabolic ball ${Q}_\rho$, and for all $(s,y)\in {\ Q}_r(t_0,x_0)$ we have the lower bound $\phi(s,y)\ge \frac{C}{r}$,
\item[2)] for all $(s,y)\in { Q}_\rho$ we have $\phi(s,y)\le \frac{C}{r}$,
\item[3)] for all $(s,y)\in { Q}_\rho$ we have $|\grad\phi(s,y)|\le \frac{C}{r^2}$,
\item [4)] moreover, for all $(s,y)\in { Q}_\rho(t,x) $ we have $|(\partial_s+\Delta)\phi(s,y)|\le C\frac{r^2}{\rho^5}$.
\end{itemize}
%%%%%%%%%%%%%%%%%%%%%%%%%%%%%%%%%%%%%%%%%%%%%%%%%%%
A detailed proof of the  properties above can be found for instance in  \cite{PGLR1}.\\ 

Now, by considering the aforementioned function $\phi$ in the local energy inequality \eqref{Ineq_Energie1}, we easily obtain
\begin{eqnarray}
\mathcal{A}_r+\alpha_r &\leq &\underbrace{\int_{s<t}\int_{\R}(\partial_t\phi +\Delta \phi)|\vu|^2dyds}_{(1)}+2\underbrace{\int_{s<t}\int_{\R}p (\vu\cdot \grad\phi )dyds}_{(2)}\notag\\
&&+\underbrace{\int_{s<t}\int _{\R} |\vu|^2(\vu \cdot \grad)\phi dx ds}_{(3)}+\underbrace{\int_{s<t}\int_{\R}(\rot \vw)\cdot (\phi\vu)dyds}_{(4)}. \label{Ineq_Energie2}
\end{eqnarray}
Let us study each term of the right-hand side above. 
\begin{itemize}
\item[$\bullet$] For the first term  $(1)$ in (\ref{Ineq_Energie2}), by the forth property of the function $\phi$ and by  the Hölder inequality $(1=\frac{1}{3}+\frac{2}{3})$ we have
$$\int_{s<t}\int_{\R}(\partial_t\phi +\Delta \phi)|\vu|^2dyds\leq C \frac{r^2}{\rho^5}\int_{Q_\rho}|\vu|^2dyds\le C\frac{r^2}{\rho ^5}\rho^\frac{5}{3} \|\vu\|_{L^3_{t,x}(Q_\rho)}^2.$$
Moreover, by (\ref{Def_Invariants}) we have  $\|\vu\|^2_{L^3_{t,x}(Q_\rho)} =\rho^\frac{4}{3}\lambda_\rho^\frac{2}{3} $, and then
$$\int_{s<t}\int_{\R}(\partial_t\phi +\Delta \phi)|\vu|^2dyds\leq  C\frac{r^2}{\rho ^2}\lambda_\rho ^\frac{2}{3}.$$

\item[$\bullet$] For the term (2) in (\ref{Ineq_Energie2}), by the third property of the test function $\phi$ and by the H\"older inequality, we obtain
$$\int_{s<t}\int_{\R}p (\vu\cdot \grad\phi )dyds\leq \frac{C}{r^2}\int_{Q_\rho} | p | |\vu|dyds\le \frac{C}{r^2}\|p\|_{L_{t,x}^{\frac{3}{2}}(Q_\rho)}\|\vu\|_{L_{t,x}^{3}(Q_\rho)}.$$
By  (\ref{Def_Invariants}) we have $\|p\|_{L_{t,x}^{\frac{3}{2}}(Q_\rho)}=\rho^{\frac{4}{3}} \mathcal{P}_\rho^{\frac{2}{3}}$ and $\|\vu\|_{L_{t,x}^{3}(Q_\rho)}=\rho^{\frac{2}{3}}\lambda_\rho^{\frac{1}{3}}$, we can thus write by the Young inequality that $$\int_{s<t}\int_{\R}p (\vu\cdot \grad\phi )dyds\leq\frac{C}{r^2}\left(\rho^{\frac{4}{3}} \mathcal{P}_\rho^{\frac{2}{3}}\right)\left(\rho^{\frac{2}{3}}\lambda_\rho^{\frac{1}{3}}\right)\le C \frac{\rho ^2}{r^2}\mathcal{P}_\rho^{\frac{2}{3}}\lambda_\rho^{\frac{1}{3}}\le C\frac{\rho ^2}{r^2}( \mathcal{P}_\rho+\lambda_\rho).$$
\item[$\bullet$] For the term (3) in (\ref{Ineq_Energie2}), by the second property of the function $\phi$, one has
\begin{eqnarray*}
\int_{s<t}\int_{\R} |\vu|^2(\vu \cdot \grad)\phi dyds&\le&\frac{C}{r^2}\int_{Q_\rho} |\vu|^3 dx ds=C\frac{\rho ^2}{r^2}\lambda_\rho.
\end{eqnarray*}
where by (\ref{Def_Invariants}) we can write $\|\vu\|_{L_{t,x}^3(Q_\rho)}^3=\rho^{2}\lambda_\rho.$
\item Finally, for the term (4) in (\ref{Ineq_Energie2}), by the properties of the function $\phi$  and by the H\"older inequality $(1=\frac{1}{6}+\frac12+\frac 13)$, we write
\begin{eqnarray*}
\int_{s<t}\int_{\R}(\rot \vw)\cdot(\phi \vu)dyds
&\leq &C\frac{\rho^{\frac{5}{6}}}{r}\|\rot \vw\|_{L_{t,x}^2(Q_\rho)}\|\vu\|_{L_{t,x}^3 (Q_\rho)}
\le C\frac{\rho^{\frac{3}{2}}}{r}\lambda_\rho^{\frac{1}{3}},
\end{eqnarray*}
where we have used the fact that $\|\vu\|_{L_{t,x}^3(Q_\rho)}=\rho^{\frac{2}{3}}\lambda_\rho^{\frac{1}{3}}$ and $\|\rot \vw\|_{L_{t,x}^2(Q_\rho)}\le \|\rot \vw\|_{L_{t,x}^2(Q_R)}<+\infty,$ since $\vw\in L^2_t\dot H^1_x(Q_R)$.
\end{itemize}
By gathering all the previous estimates we obtain
\begin{equation*}
\mathcal{A}_r+\alpha_r \le C\frac{r^2}{\rho ^2} \lambda_\rho ^\frac{2}{3}+C\frac{\rho ^2}{r^2}\left(\mathcal{P}_\rho+\lambda_\rho\right)+C\frac{\rho^{\frac{3}{2}}}{r}\lambda_\rho^{\frac{1}{3}},
\end{equation*} 
and this ends the proof of Lemma \ref{lema_FirstEstimate}. \hfill$\blacksquare $\\

As it was pointed out in the $\varepsilon$-regularity theory for the Navier-Stokes equations (see \cite{CKN}, \cite{Kukavica} or \cite{PGLR1}), we need to study more in detail the pressure $p$, which only appears in the first equation of the micropolar system. Following the same ideas presented in our previous works \cite{ChLl21, ChLl23} (also refer to \cite{Kukavica}, \cite[Lemma 13.3]{PGLR1}), we derive the following lemma.
%%%%%%%%%%%%%%%%%%%%%%%%%%%%%%%%%%%%%%%%%%%%%%%%%%%
\begin{Lemma}\label{Lema_EstimatePressure}
Under the hypotheses of Proposition \ref{Pro_partialRegualrity} for any $\displaystyle{0<r\le\tfrac\rho{2}\le R}$, we have the inequality
\begin{equation}\label{SecondEstimate}
\mathcal{P}_r^{\frac{2}{3}}\leq C\left(\left(\frac{\rho}{r}\right)^{\frac{4}{3}}\lambda_\rho^{\frac{2}{3}}+\left(\frac{r}{\rho}\right)^{\frac{2}{3}}\mathcal{P}_\rho^{\frac{2}{3}}\right).
\end{equation}
\end{Lemma}
%%%%%%%%%%%%%%%%%%%%%%%%%%%%%%%%%%%%%%%%%%%%%%%%%%%
\noindent {\bf Proof.} First, let us prove the following estimate
\begin{equation}\label{Estimation_Pression_Intermediaire}
\|p\|_{L_{t,x}^{\frac{3}{2}}(Q_\sigma)} \leq C \left( \|\vu\|_{L_{t,x}^{3} (Q_1)}+\sigma^{2}\|p\|_{L_{t,x}^{\frac{3}{2}} (Q_1)}\right),
\end{equation}
where $Q_\sigma$ and $Q_1$ are parabolic balls of radius $\sigma$ and $1$ respectively. Later, we will derive \eqref{SecondEstimate} by a change of variable. \\

In order to obtain \eqref{Estimation_Pression_Intermediaire}, we introduce $\eta : \R \longrightarrow \mathbb{R}$ a smooth function supported in the ball $B_{0,1}$ such that $\eta \equiv 1$ on the ball $B_{0,\frac{3}{5}}$ and $\eta \equiv 0$ outside the ball $B_{0,\frac{4}{5}}$. Fix $0<\sigma\le \frac{1}2$ and notice that $p=\eta p$ in  $B_{0,\sigma}$. Now, by using the identity
$$ - \Delta (\eta p) = -\eta \Delta p + (\Delta \eta)p - 2\sum^3_{i= 1}\partial_i ((\partial_i \eta) p),$$
we deduce the inequality
\begin{eqnarray}
\|p\|_{L_{t,x}^{\frac{3}{2}}(Q_\sigma)} =\|\eta p\|_{L_{t,x}^{\frac{3}{2}}(Q_\sigma)} &\leq& \underbrace{\left\|\frac{\big(-\eta \Delta p\big)}{(-\Delta)}\right\|_{L_{t,x}^{\frac{3}{2}}(Q_\sigma)}}_{(p_1)} + \underbrace{\left\|\frac{(\Delta \eta) p}{(-\Delta)}\right\|_{L_{t,x}^{\frac{3}{2}}(Q_\sigma)} }_{(p_2)} \notag\\
&&+2\sum^3_{i= 1}\underbrace{\left\|\frac{\partial_i ( (\partial_i \eta) p)}{(-\Delta)}\right\|_{L_{t,x}^{\frac{3}{2}}(Q_\sigma)}}_{(p_3)}.\label{FormuleEtaPression}
\end{eqnarray}
For the first term of (\ref{FormuleEtaPression}), since we have the equation $\Delta p=-\displaystyle{\sum^3_{i,j= 1}}\partial_i \partial_j (u_i u_j)$, we can write
\begin{eqnarray}
(p_1)&=&\left\|\frac{\big(-\eta \Delta p\big)}{(-\Delta)}\right\|_{L_{t,x}^{\frac{3}{2}}(Q_\sigma)}\leq C\left\|\frac{1}{(-\Delta)}\Big(\eta \,\sum^3_{i,j= 1}\partial_i \partial_j u_i u_j\Big)\right\|_{L_{t,x}^{\frac 32}(Q_\sigma)}\notag\\
&\leq &C\sum^3_{i,j= 1}\left\|\frac{1}{(-\Delta)} \left(\partial_i \partial_j(\eta  u_i u_j ) - \partial_i \big((\partial_j \eta)  u_i u_j \big) - \partial_j \big((\partial_i \eta)  u_i u_j \big) + (\partial_i \partial_j \eta)  u_i u_j\right)\right\|_{L_{t,x}^{\frac 32}(Q_\sigma)}\notag\\
&\leq &C\sum^3_{i,j= 1}\left\|\frac{1}{(-\Delta)}\partial_i \partial_j(\eta  u_i u_j )\right\|_{L_{t,x}^{\frac 32}(Q_\sigma)}+C\sum^3_{i,j= 1}\left\|\frac{1}{(-\Delta)} \partial_i \big((\partial_j \eta)  u_i u_j \big)\right\|_{L_{t,x}^{\frac 32}(Q_\sigma)}\notag\\&&+C\sum^3_{i,j= 1}\left\|\frac{1}{(-\Delta)}  \partial_j \big((\partial_i \eta)  u_i u_j \big) \right\|_{L_{t,x}^{\frac 32}(Q_\sigma)}+C\sum^3_{i,j= 1} \left\|\frac{1}{(-\Delta)} (\partial_i \partial_j \eta)  u_i u_j\right\|_{L_{t,x}^{\frac 32}(Q_\sigma)}.\label{FormulePression_I}
\end{eqnarray}
Let us study each term of the expression above. Denoting by $\mathcal{R}_i=\frac{\partial_i}{\sqrt{-\Delta}}$ the usual Riesz transforms on $\mathbb{R}^3$, by the boundedness of these operators in $L^\frac32(\R)$, and using the support properties of the auxiliary function $\eta$, we have for the first term above:
\begin{eqnarray*}
\left\|\frac{\partial_i \partial_j}{(-\Delta)} \eta  u_i u_j(t,\cdot) \right\|_{L^{\frac32} (B_{0,\sigma})}
&\leq &\|\mathcal{R}_i \mathcal{R}_j (\eta  u_i u_j )(t,\cdot) \|_{L^{\frac32}(\R)}\leq C\|\eta  u_i u_j(t,\cdot)\|_{L^{\frac32}(\R) }\\ 
&\le &C\|\vu(t,\cdot) \|^2_{L^{3}(B_{0,1})},
\end{eqnarray*}
By taking the $L^{\frac32}$-norm in  the time variable in the previous inequality we obtain
\begin{equation}\label{p11}
\left\|\frac{\partial_i \partial_j}{(-\Delta)} \eta  u_i u_j \right\|_{L_{t,x}^{\frac32} (Q_\sigma)}\leq C\|\vu(t,\cdot) \|^2_{L^{3}(Q_1)}.
\end{equation}
The remaining terms of (\ref{FormulePression_I}) can all be studied in a similar manner. Indeed, noting that $\partial_i \eta$ vanishes on $B_{\frac35} \cup B^c_{\frac45}$ by using the integral representation for the operator $\frac{\partial_i}{(- \Delta )}$ we have for the second term of (\ref{FormulePression_I}) the estimate
\begin{eqnarray}
\left\|\frac{\partial_i}{(- \Delta )}\big((\partial_j\eta) u_i u_j\big)(t,\cdot)\right\|_{L^{\frac32}(B_{0,\sigma})} &\leq &C\sigma^2\left\|\frac{\partial_i}{(- \Delta )}\big((\partial_j\eta) u_i u_j\big)(t,\cdot)\right\|_{L^{\infty}(B_{0,\sigma})}\notag\\
&\leq &C \, \sigma^2 \left\|\int_{\{\frac35<|y|<\frac45\}} \frac{x_i - y_i }{|x-y|^3} \big((\partial_j \eta)  u_i u_j \big)(t,y)\, dy\right\|_{L^{\infty}(B_{0,\sigma})}\notag
\end{eqnarray}
Now, since $x\in B_{0,\sigma}$ and $\sigma\le \frac12$, we have for any $\frac35<|y|<\frac45$ that $\frac1{10}<|x-y|$ and since $\supp(\eta_j)\subset B_{0,1}$,  it follows that 
\begin{eqnarray}
\left\|\frac{\partial_i}{(- \Delta )}\big((\partial_j\eta) u_i u_j\big)(t,\cdot)\right\|_{L^{\frac32}(B_{0,\sigma})} 
& \leq &C  \| u_i u_j(t,\cdot)\|_{L^{1} (B_{0,1})}\label{KernelEstimate1}\\
& \leq &C \|\vu (t,\cdot) \|_{L^{3} (B_{0,1})}^2.\notag
\end{eqnarray}
By  taking the $L^{\frac32}$-norm in the time variable in the expression above, we obtain
\begin{equation}\label{p12}
\left\|\frac{\partial_i}{(- \Delta )}\big((\partial_j\eta) u_i u_j\big)\right\|_{L^{\frac32}_{t,x}(Q_\sigma)} \leq C \|\vu\|_{L_{t,x}^{3} (Q_1)}^2.
\end{equation}
A symmetric argument gives 
\begin{equation}\label{p13}
\left\|\frac{\partial_{j}}{(- \Delta )}\big((\partial_i\eta) u_i u_j\big)\right\|_{L^{\frac32}_{t,x}(Q_\sigma)} \leq C \|\vu\|_{L_{t,x}^{3} (Q_1)}^2.
\end{equation}
Since the convolution kernel associated to the operator $\frac{1}{(-\Delta)}$ is $\frac{C}{|x|}$, by following the same ideas as in \eqref{KernelEstimate1}, we easily obtain for the last term of (\ref{FormulePression_I}) that
\begin{equation}\label{p14}
\left\|\frac{(\partial_i \partial_j \eta)  u_i u_j}{(-\Delta)}\right\|_{L_{t,x}^{\frac 32}(Q_\sigma)}\leq C \|\vu\|_{L_{t,x}^{3} (Q_1)}^2
\end{equation}
By merging the estimates \eqref{p11}, \eqref{p12},\eqref{p13} and\eqref{p14} in \eqref{FormulePression_I}, we obtain 
\begin{eqnarray}
(p_1)=\left\|\frac{\big(-\eta \Delta p\big)}{(-\Delta)}\right\|_{L_{t,x}^{\frac 32}(Q_\sigma)}\leq C \|\vu\|_{L_{t,x}^{3} (Q_1)}^2.\label{FormulePression_I1}
\end{eqnarray}
For treating the term $(p_2)$ in (\ref{FormuleEtaPression}), by the properties of the auxiliary function $\eta$ and  the convolution kernel associated to the operator $\frac{1}{(-\Delta)}$, we can write (see (\ref{KernelEstimate1})):
$$\left\|\frac{(\Delta \eta) p(t,\cdot)}{(-\Delta)}\right\|_{L^{\frac{3}{2}}(B_{0,\sigma})}\leq C\sigma^{2} \| p(t,\cdot)\|_{L^1(B_{0,1})}\leq C\sigma^{2} \| p(t,\cdot)\|_{L^{\frac{3}{2}}(B_{0,1})},$$
and thus, taking the $L^{\frac{3}{2}}$-norm in the time variable we obtain:
\begin{equation}\label{FormulePression_II}
(p_2)=\left\|\frac{(\Delta \eta) p}{(-\Delta)}\right\|_{L^{\frac{3}{2}}_{t}L^{\frac{3}{2}}_{x}(Q_\sigma)}\leq C\sigma^{2} \|p\|_{L^{\frac{3}{2}}_{t}L^{\frac{3}{2}}_{x}(Q_1)}.
\end{equation}
For the last term  $(p_3)$ of  (\ref{FormuleEtaPression}), following the same ideas developed in (\ref{KernelEstimate1}) we can write
$$\left\|\frac{\partial_i}{(-\Delta)}(\partial_i \eta) p(t,\cdot) \right\|_{L^{\frac{3}{2}}(B_{0,\sigma})}\leq C\sigma^{2} \| p(t,\cdot)\|_{L^1(B_{0,1})}\leq C\sigma^{2} \|p(t,\cdot)\|_{L^{\frac{3}{2}}(B_{0,1})},$$
and therefore 
\begin{equation}\label{FormulePression_III}
(p_3)=\left\|\frac{\partial_i ( (\partial_i \eta) p )}{(-\Delta)}\right\|_{L_{t,x}^{\frac{3}{2}}(Q_\sigma)}\leq C\sigma^{2} \|p\|_{L^{\frac{3}{2}}_{t,x}(Q_1)}.
\end{equation}
Now, gathering the estimates (\ref{FormulePression_I1}), (\ref{FormulePression_II}) and (\ref{FormulePression_III}) we obtain the inequality 
$$\|p\|_{L_{t,x}^{\frac{3}{2}}(Q_\sigma)} \leq C \left( \|\vu\|^2_{L^3_{t,x}(Q_1)} +
\sigma^{2}\|p\|_{L_{t,x}^{\frac{3}{2}} (Q_1)}\right).$$

Now, with this estimate at hand, it is straightforward to deduce inequality \eqref{SecondEstimate}. Indeed, if we fix $\sigma = \frac{r}{\rho} \leq \frac12$ and by introducing the functions $p_\rho(t,x)=p(\rho^2t, \rho x)$ and $\vu_\rho(t,x)=\vu(\rho^2t, \rho x)$ then  the previous estimate we have
$$\|p_\rho\|_{L_{t,x}^{\frac{3}{2}}(Q_{\frac{r}{\rho}})} \leq C \left(\|\vu_\rho\|_{L_{t,x}^{3} (Q_1)}^2 +\left(\frac{r}{\rho}\right)^{2}\|p_\rho\|_{L_{t,x}^{\frac{3}{2}} (Q_1)}\right).$$
Hence by a convenient change of variable we obtain
$$\|p\|_{L_{t,x}^{\frac{3}{2}}(Q_{r})}\rho^{-\frac{10}{3}} \leq C \left( \rho^{-\frac{10}{3}}\|\vu\|^2_{L^3_{t,x}(Q_\rho)} +\left(\frac{r}{\rho}\right)^{2}\rho^{-\frac{10}{3}}\|p\|_{L_{t,x}^{\frac{3}{2}} (Q_\rho)}\right).$$
Moreover, by (\ref{Def_Invariants}) we have the identities $r^{\frac{4}{3}}\mathcal{P}_r^{\frac{2}{3}}=\|p\|_{L_{t,x}^{\frac{3}{2}}(Q_r)}$ and $\rho^{\frac{4}{3}}\lambda_\rho^{\frac{2}{3}}=\| \vu\|^2_{L_{t,x}^{3} (Q_\rho)},$
and therefore we obtain $\displaystyle{\mathcal{P}_r^{\frac{2}{3}}\leq C\left(\frac{\rho}{r}\right)^{\frac{4}{3}}\lambda_\rho^{\frac{2}{3}}+C\left(\frac{r}{\rho}\right)^{\frac{2}{3}}\mathcal{P}_\rho^{\frac{2}{3}}}$ and this finishes the proof of Lemma \ref{Lema_EstimatePressure}.
\hfill $\blacksquare $\\

%%%%%%%%%%%%%%%%%%%%%%%%%%%%%%%%%%%%%%%%%%%%%%%%%%%
\noindent \textbf{End of the proof of the Proposition \ref{Pro_partialRegualrity}} Now, in order to deduce \eqref{IterativeO1}, we remark that it is equivalent to say that  there exists $0<\mathfrak{r}<R$ and  $0<\kappa<\frac12$  such that for all $n\in \mathbb{N}$ and for all $(t,x)\in Q_{\kappa^n\mathfrak{r}}(t_0,x_0)$,  we have
\begin{equation}\label{IterativeO}
\mathbb{O}_{\kappa^n\mathfrak r}(t,x)=\Lambda_{\kappa^n\mathfrak{r}}(t,x)+\kappa^6\mathbb{P}_{\kappa^n\mathfrak{r}}(t,x)\le C,
\end{equation}
where
\begin{equation}\label{TermLambdaPressure}
\Lambda_r(t,x)=\frac{1}{r^{3(1-\frac{5}{\tau_0})}}\lambda_r(t,x)\quad\text{and}\quad \mathbb{P}_r(t,x)=\frac{1}{r^{3(1-\frac{5}{\tau_0})}}\mathcal{P}_r(t,x).
\end{equation}
Thus, for deducing \eqref{IterativeO}, we will apply an iterative argument and to do so we need to estimate  $\Lambda_r$ and $\mathbb{P}_r$ in terms of $\Lambda_\rho$ and $\mathbb{P}_\rho$ for any radius $0<r<\frac \rho 2\le R$.  Thus, by  Lemmas \ref{Estimation_normeLambda3} and \ref{lema_FirstEstimate} we have
\begin{eqnarray}
\Lambda_r=\frac{1}{r^{3(1-\frac{5}{\tau_0})}}\lambda_r&\le& \frac{C}{r^{3(1-\frac{5}{\tau_0})}}(\mathcal{A}_{r}+\alpha_r)^{\frac{3}{2}}\notag\\
&\leq& \frac{C}{r^{3(1-\frac{5}{\tau_0})}}\frac{r^3}{\rho ^3} \lambda_\rho+\frac{C}{r^{3(1-\frac{5}{\tau_0})}}\frac{\rho ^3}{r^3}\left(\mathcal{P}_\rho+\lambda_\rho\right)^{\frac{3}{2}}+\frac{C}{r^{3(1-\frac{5}{\tau_0})}}\frac{\rho^{\frac{9}{4}}}{r^{\frac{3}{2}}}\lambda_\rho^{\frac{1}{2}}.\label{EstimationPourIteration}
\end{eqnarray}
Let us study more in detail each term of the right-hand side above. 
\begin{itemize}
\item For the first term of \eqref{EstimationPourIteration}, since $\lambda_\rho={\rho^{3(1-\frac{5}{\tau_0})}}\Lambda_\rho$  by \eqref{TermLambdaPressure}, we have
$$\frac{C}{r^{3(1-\frac{5}{\tau_0})}}\frac{r^3}{\rho ^3}\lambda_\rho =C \left(\frac{r}{\rho}\right)^{\frac{15}{\tau_0}}\Lambda_\rho.$$
\item For the second term of (\ref{EstimationPourIteration}), by the definition of $\mathbb{P}_\rho$ and $\Lambda_\rho,$ given in  by \eqref{TermLambdaPressure}, we obtain
\begin{eqnarray*}
\frac{1}{r^{3(1-\frac{5}{\tau_0})}}\frac{\rho ^3}{r^3}\left(\mathcal{P}_\rho+\lambda_\rho\right)^{\frac{3}{2}}
&=& C\left(\frac{\rho}{r}\right)^{6-\frac{15}{\tau_0}}\rho^{-3+\frac{15}{\tau_0}}\rho^{\frac{9}{2}(1-\frac{5}{\tau_0})}\left(\mathbb{P}_\rho+\Lambda_\rho\right)^{\frac{3}{2}}\\
&=& C\left(\frac{\rho}{r}\right)^{6-\frac{15}{\tau_0}}\rho^{\frac32-\frac{15}{2\tau_0}}(\mathbb{P}_\rho+\Lambda_\rho)^{\frac{3}{2}}\\
&\le&  C\left(\frac{\rho}{r}\right)^{6-\frac{15}{\tau_0}}\left(\mathbb{P}_\rho^{\frac{3}{2}}+\Lambda_\rho^{\frac{3}{2}}\right),
\end{eqnarray*}
where we have used the fact that $\rho^{\frac32-\frac{15}{2\tau_0}}<1$ since $\frac32-\frac{15}{2\tau_0}>0$ due to $\tau_0>5$.
\item Finally, for the last term of (\ref{EstimationPourIteration}),  by \eqref{TermLambdaPressure}, we have
\begin{align*}
\frac{1}{r^{3(1-\frac{5}{\tau_0})}}\frac{\rho^{\frac{9}{4}}}{r^{\frac{3}{2}}} \lambda_{\rho}^{\frac{1}{2}}&= \frac{1}{r^{3(1-\frac{5}{\tau_0})}}\rho^{\frac{3}{2}(1-\frac{5}{\tau_0})}\left(\frac{\rho^{\frac94}}{r^{\frac32}}\right) \Lambda_{\rho}^{\frac{1}{2}}\\
&=\rho^{\frac{3}{2}(-\frac12+\frac{5}{\tau_0})}\left(\frac{\rho}{r}\right)^{\frac92-\frac{15}{\tau_0}}\Lambda_{\rho}^{\frac{1}{2}}.
\end{align*}
\end{itemize}
Thus, by gathering all the previous estimates, we have 
\begin{align}\label{estimateA}
\Lambda_r\leq C\Bigg(\left(\frac{r}{\rho}\right)^{\frac{15}{\tau_0}}\Lambda_\rho+\left(\frac{\rho}{r}\right)^{6-\frac{15}{\tau_0}}\left(\mathbb{P}_\rho^{\frac{3}{2}}+\Lambda_\rho^{\frac{3}{2}}\right)+\rho^{\frac{3}{2}(-\frac12+\frac{5}{\tau_0})}\left(\frac{\rho}{r}\right)^{\frac92-\frac{15}{\tau_0}}\Lambda_{\rho}^{\frac{1}{2}}\Bigg).
\end{align}
Let us study now the pressure term in \eqref{TermLambdaPressure}. From the estimate \eqref{SecondEstimate}, we can write for any $0<r\le \frac{\rho}{2}$
$$\mathbb{P}_r=\frac{1}{r^{3(1-\frac{5}{\tau_0})}}\mathcal{P}_r\le  \frac{C}{r^{3(1-\frac{5}{\tau_0})}}\left(\frac{\rho}{r}\right)^{2}\lambda_\rho+\frac{C}{r^{3(1-\frac{5}{\tau_0})}}\frac{r}{\rho}\mathcal{P}_\rho.$$
For the first term of the right-hand side above, since $\lambda_\rho={\rho^{3(1-\frac{5}{\tau_0})}}\Lambda_\rho$  by \eqref{TermLambdaPressure},  we have
$$\frac{1}{r^{3(1-\frac{5}{\tau_0})}}\left(\frac{\rho}{r}\right)^2\lambda_\rho=\frac{1}{r^{3(1-\frac{5}{\tau_0})}}\left(\frac{\rho}{r}\right)^{2}\rho^{3(1-\frac{5}{\tau_0})}\Lambda_\rho=\left(\frac{\rho}{r}\right)^{5-\frac{15}{\tau_0}}\Lambda_\rho.$$
Moreover, by using the fact that $\frac{1}{r^{3(1-\frac{5}{\tau_0})}} \frac{r}{\rho}\mathcal{P}_\rho= \left(\frac{\rho}{r}\right)^{2-\frac{15}{\tau_0}}\mathbb{P}_\rho$ by  \eqref{TermLambdaPressure}, one has
\begin{align}\label{estimatePressure}
\mathbb{P}_r\le C\biggl(\left(\frac{\rho}{r}\right)^{5-\frac{15}{\tau_0}}\Lambda_\rho+\left(\frac{\rho}{r}\right)^{2-\frac{15}{\tau_0}}\mathbb{P}_\rho\biggr).
\end{align}
Hence, we have estimated $\Lambda_r$ and $\mathbb{P}_{r}$ in terms of   $\Lambda_\rho$ and $\mathbb{P}_{\rho}$.\\[4mm]

With this information at hand, let us  study the expression $\mathbb{O}_r$ given in \eqref{IterativeO}. Thus, notice that for any $\displaystyle{0<r\le \tfrac{\rho}{2}}$, if we fix $0<\kappa <\frac 12$ such that $\kappa=\frac{r}{\rho}$ , then from the estimates (\ref{estimateA}) and  \eqref{estimatePressure} it follows  that
\begin{eqnarray*}
\mathbb{O}_r=\Lambda_r+\kappa^{6}\mathbb{P}_r&\le& C\left( \kappa^{\frac{15}{\tau_0}}\Lambda_\rho+\kappa^{-6+\frac{15}{\tau_0}}\left(\mathbb{P}_\rho^{\frac{3}{2}}+\Lambda_\rho^{\frac{3}{2}}\right)\notag+\rho^{\frac{3}{2}(-\frac12+\frac{5}{\tau_0})}\kappa^{3(-\frac32+\frac{5}{\tau_0})}\Lambda_{\rho}^{\frac{1}{2}}\right)\\
&&+C\kappa^6\left(\kappa^{-5+\frac{15}{\tau_0}}\Lambda_\rho+\kappa^{-2+\frac{15}{\tau_0}}\mathbb{P}_\rho\right).
\end{eqnarray*}
Moreover, by the definition of $\mathbb{O}_\rho$ given in \eqref{Def_QuantiteIteration}, we have $\Lambda_\rho\le \mathbb{O}_\rho$ and $\mathbb{P}_{\rho}\le k^{-6}\mathbb{O}_{\rho}$ and therefore one has
\begin{eqnarray*}
\mathbb{O}_r&\le&C\kappa^{\frac{15}{\tau_0}}\mathbb{O}_\rho+C\kappa^{-6+\frac{15}{\tau_0}}\left(\kappa^{-9}\mathbb{O}_\rho^{\frac{3}{2}}+\mathbb{O}_\rho^{\frac{3}{2}}\right)+C\rho^{\frac{3}{2}(-\frac12+\frac{5}{\tau_0})}\kappa^{3(-\frac32+\frac{5}{\tau_0})}\mathbb{O}_{\rho}^{\frac{1}{2}}\notag\\
&&+C\kappa^{1+\frac{15}{\tau_0}}\mathbb{O}_\rho+C\kappa^{-2+\frac{15}{\tau_0}}\mathbb{O}_\rho.
\end{eqnarray*}
By using the Young inequality $(1=\frac23+\frac13)$, we have
\begin{eqnarray*}
\mathbb{O}_r&\le&C\kappa^{\frac{15}{\tau_0}}\mathbb{O}_\rho+C\kappa^{-6+\frac{15}{\tau_0}}\left(\kappa^{-9}\mathbb{O}_\rho^{\frac{3}{2}}+\mathbb{O}_\rho^{\frac{3}{2}}\right)+\rho^{\frac{9}{4}(-\frac12+\frac{5}{\tau_0})}+C\kappa^{9(-\frac32+\frac{5}{\tau_0})}\mathbb{O}_{\rho}^{\frac{3}{2}}\notag\\
&&+C\kappa^{1+\frac{15}{\tau_0}}\mathbb{O}_\rho+C\kappa^{-2+\frac{15}{\tau_0}}\mathbb{O}_\rho.
\end{eqnarray*}
Rearranging the previous expression in a more convenient way, we obtain
\begin{eqnarray*}
\mathbb{O}_r&\le& C\left(\kappa^{\frac{15}{\tau_0}}+\kappa^{1+\frac{15}{\tau_0}}+\kappa^{-2+\frac{15}{\tau_0}}\right)\mathbb{O}_{\rho}+C\left(\kappa^{-15+\frac{15}{\tau_0}}+\kappa^{-6+\frac{15}{\tau_0}}+\kappa^{9(-\frac32+\frac{5}{\tau_0})}\right)\mathbb{O}_{\rho}^{\frac{1}{2}}\mathbb{O}_{\rho}\notag\\&&\quad +\rho^{\frac{9}{4}(-\frac12+\frac{5}{\tau_0})}\\
&\le& C\left(\kappa^{\frac{15}{\tau_0}}+\kappa^{1+\frac{15}{\tau_0}}+\kappa^{-2+\frac{15}{\tau_0}}\right)\mathbb{O}_{\rho}+C\kappa^{-15}\mathbb{O}_{\rho}^{\frac{1}{2}}\mathbb{O}_{\rho}+\rho^{\frac{9}{4}(-\frac12+\frac{5}{\tau_0})}.
\end{eqnarray*}
Moreover since $-2+\frac{15}{\tau_0}>0$, we  take  $0<\kappa\ll1$ small enough such that 
\begin{equation}\label{conditionKappa}
C(\kappa^{\frac{15}{\tau_0}}+\kappa^{1+\frac{15}{\tau_0}}+\kappa^{-2+\frac{15}{\tau_0}})\le \frac{1}{4},
\end{equation}
 and therefore for any $0<r\le \frac \rho2$ we obtain the following estimation
\begin{equation}
\mathbb{O}_r\le  \frac14\mathbb{O}_{\rho}+C\kappa^{-15}\mathbb{O}_{\rho}^{\frac{1}{2}}\mathbb{O}_{\rho}+\rho^{\frac{9}{4}(-\frac12+\frac{5}{\tau_0})}.\label{EstimateAux}
\end{equation} 
Now by using the estimate above, we can deduce  \eqref{IterativeO}, \textit{i.e.,} we will prove for $0<\kappa\ll1$ given by the condition \eqref{conditionKappa},  there exists $0<\mathfrak{r}<R$ such that for all $n\in \mathbb{N}$ and for all $(t,x)\in Q_{\kappa^n\mathfrak{r}}(t_0,x_0)$,  we have
\begin{equation*}
\mathbb{O}_{\kappa^n\mathfrak{r}}(t,x)=\Lambda_{\kappa^n\mathfrak{r}}(t,x)+\kappa^6\mathbb{P}_{\kappa^n\mathfrak{r}}(t,x)\le C.
\end{equation*} 
Indeed, let us define $\rho=r_0=\mathfrak{r}$ and $r=r_1=\kappa \mathfrak{r}$ with $\mathfrak{r}=\kappa^\mathfrak{N}R$ where  $0<\kappa\ll 1$ is given by the condition \eqref{conditionKappa}, and $\mathfrak{N}\in \mathbb{N}$ is such that $\mathfrak{N}>240$.  Thus,  since $\displaystyle{r_1\le\frac{r_0}{2}}$, we can rewrite \eqref{EstimateAux} as follows
\begin{eqnarray*}
\mathbb{O}_{r_{1}}(t_0,x_0)&\le& \frac14\mathbb{O}_{r_0}+C\kappa^{-15}\mathbb{O}_{r_0}^{\frac{1}{2}}\mathbb{O}_{r_0}\notag+r_0^{\frac{9}4(-\frac12+\frac5{\tau_0})}.
\end{eqnarray*}
Since $r_0=\kappa^\mathfrak{N}R$,  $R<1$ and $-\frac12+\frac5{\tau_0}>0$ we have
\begin{eqnarray}
\mathbb{O}_{r_{1}}(t_0,x_0)&\le& \frac14\mathbb{O}_{r_0}+C\kappa^{-15}\mathbb{O}_{r_0}^{\frac{1}{2}}\mathbb{O}_{r_0}+\kappa^{\frac{9\mathfrak{N}}4(-\frac12+\frac5{\tau_0})}.\label{FirstIterationAux}
\end{eqnarray} 
In order to close the iterative argument, we need to study each term  of the right-hand side above.  First, notice that since $0<\kappa\ll 1$, the expression $C\kappa^{-15}$ can be large, nevertheless since $\kappa>0$  is a fixed parameter, we may consider a parameter $0<\varepsilon_\ast\ll 1$ small enough (to be defined later on)  such that we have
\begin{equation}\label{NegativePowers}
C\kappa^{-15}\varepsilon_\ast^\frac12= \frac 14.
\end{equation}
On the other hand by  \eqref{Def_Invariants}, \eqref{Def_QuantiteIteration} and since $Q_{(\kappa^\mathfrak{N}R)}(t_0,x_0)\subset Q_{R}(t_0,x_0)$, we have
\begin{eqnarray*}
\mathbb{O}_{r_0}(t_0,x_0)&=&\frac{1}{{(\kappa^\mathfrak{N}R)}^{3(1-\frac{5}{\tau_0})}}\left(\frac{1}{(\kappa^\mathfrak{N}R)^2} \int _{Q_{(\kappa^\mathfrak{N}R)}(t_0,x_0)}|\vu|^3dyds+\frac{\kappa^{6}}{(\kappa^\mathfrak{N}R)^2} \int _{Q_{(\kappa^\mathfrak{N}R)}(t_0,x_0)}|p |^{\frac{3}{2}}dyds\right)\\
&= & \frac {R^{\frac{15}{\tau_0}}}{\kappa^{5\mathfrak{N}}}\frac{1}{{R}^{3}}\left(\frac{1}{ R^2} \int _{Q_{R}(t_0,x_0)}|\vu|^3dyds+\frac{1}{R^2} \int _{Q_{R}(t_0,x_0)}|p|^{\frac{3}{2}}dyds\right)\\
&\le & \frac {1}{\kappa^{5\mathfrak{N}}}\frac{1}{{R}^{3}}\left(\frac{1}{ R^2} \int _{Q_{R}(t_0,x_0)}|\vu|^3+|p|^{\frac{3}{2}}dyds\right).
\end{eqnarray*}
Now, recall that by the hypothesis \eqref{TeohyL3}, there exists  $0<\varepsilon\ll1$ such that $\displaystyle{\frac{1}{R^2}\int_{Q_R(t_0,x_0)}| \vu|^3+|p|^\frac{3}{2}dy ds<\varepsilon}$ and thus by setting $\varepsilon$ such that ${0<\varepsilon\le \kappa^{5\mathfrak{N}} R^3\varepsilon_\ast}$ (where $\varepsilon_\ast$ was given by the condition (\ref{NegativePowers})),  it follows that  
\begin{equation}\label{casen=0}
\mathbb{O}_{r_0}(t_0,x_0)\le \varepsilon_\ast.
\end{equation}
Then, by \eqref{NegativePowers} and the expression above, it follows that
\begin{equation}\label{ConditionSecondTerm}
C\kappa^{-15}\mathbb{O}_{r_0}^{\frac{1}{2}}\mathbb{O}_{r_0}\le C\kappa^{-15}\varepsilon_\ast^{\frac{1}{2}}\mathbb{O}_{r_0}\le \frac14\mathbb{O}_{r_0}.
\end{equation}
Now, let us study the last term in the right-hand side of \eqref{FirstIterationAux}, notice that  by \eqref{NegativePowers},  we can write $\kappa=C\varepsilon_\ast^\frac1{30}$ and therefore we have
\begin{equation*}
\kappa^{\frac{9\mathfrak{N}}4(-\frac12+\frac5{\tau_0})}=C\varepsilon_\ast^{\frac{9\mathfrak{N}}{120}(-\frac12+\frac5{\tau_0})}
\end{equation*}
Thus, since  $\mathfrak{N}>240$, recalling that $-\frac12+\frac5{\tau_0}>0$ and $\frac1{\tau_0}\ge \frac {2}{15}$ we have $\frac{9\mathfrak{N}}{120}(-\frac12+\frac5{\tau_0})>3$, and therefore we can take $\varepsilon_\ast\ll1 $ such that
\begin{equation}\label{conditionKappaAlone}
\kappa^{\frac{9\mathfrak{N}}4(-\frac12+\frac5{\tau_0})}=C\varepsilon_\ast^{\frac{9\mathfrak{N}}{120}(-\frac12+\frac5{\tau_0})}\le \frac{\varepsilon_\ast}2.
\end{equation}
Hence, by  using  \eqref{ConditionSecondTerm} and \eqref{conditionKappaAlone} in \eqref{FirstIterationAux} we obtain
\begin{equation*}
\mathbb{O}_{r_{1}}(t_0,x_0)\le  \frac14\mathbb{O}_{r_0}+\frac14\mathbb{O}_{r_0}+\frac{\varepsilon_\ast}{2}.
\end{equation*}
Furthermore, by \eqref{casen=0}, it follows that
\begin{equation}\label{casen=1}
\mathbb{O}_{r_{1}}(t_0,x_0)\le  \frac{\varepsilon_\ast}{4}+\frac{\varepsilon_\ast}{4}+\frac{\varepsilon_\ast}{2}=\varepsilon_\ast.
\end{equation}
Now, let us study the case $r_2=\kappa^2\mathfrak{r}=\kappa^{\mathfrak{N}+2}R$. Since $r_2\le \frac{\kappa^\mathfrak{\mathfrak{N}}R}{2}= \frac{r_1}{2}$  we can apply the estimate \eqref{EstimateAux}  and hence we have
\begin{eqnarray*}
\mathbb{O}_{r_{2}}(t_0,x_0)&\le& \frac14\mathbb{O}_{r_1}+C\kappa^{-15}\mathbb{O}_{r_1}^{\frac{1}{2}}\mathbb{O}_{r_1}+r_1^{\frac{9}4(-\frac12+\frac5{\tau_0})}.
\end{eqnarray*}
Notice that since $\mathbb{O}_{r_{1}}(t_0,x_0)\le \varepsilon_\ast$ by \eqref{NegativePowers}, we have $ C\kappa^{-15}\mathbb{O}_{r_1}^{\frac{1}{2}}\le C\kappa^{-15}(\varepsilon_\ast)^{\frac{1}{2}}\le \frac14$. Moreover, since $r_1=\kappa^{\mathfrak{N}+1}R\le\kappa^\mathfrak{N}$ (recall $R<1$) by \eqref{conditionKappaAlone} we have $r_1^{\frac{9}4(-\frac12+\frac5{\tau_0})}\le \frac{\varepsilon}{2}.$ Then, one has
\begin{eqnarray*}
\mathbb{O}_{r_{2}}(t_0,x_0)&\le&\frac14\mathbb{O}_{r_1}+\frac14\mathbb{O}_{r_1}+\frac{\varepsilon_\ast}{2}.
\end{eqnarray*} 
Again since $\mathbb{O}_{r_{1}}(t_0,x_0)\le \varepsilon_\ast$ by \eqref{casen=1}we have
\begin{equation*}
\mathbb{O}_{r_{2}}(t_0,x_0)\le \frac{\varepsilon_\ast}{4}+\frac{\varepsilon_\ast}{4}+\frac{\varepsilon_\ast}{2}=\varepsilon_\ast.
\end{equation*}
Finally, let us consider the case $r_n=\kappa^n\mathfrak{r}=\kappa^{\mathfrak{N}+n}R$ and we assume that $\mathbb{O}_{r_n}\le \varepsilon_\ast$. Thus, since $r_{n+1}\le  \frac{r_n}2$, by using  \eqref{EstimateAux} we have
\begin{equation*}
\mathbb{O}_{r_{n+1}}(t_0,x_0)\le \frac14\mathbb{O}_{r_n}+C\kappa^{-15}\mathbb{O}_{r_n}^{\frac{1}{2}}\mathbb{O}_{r_n}+r_n^{\frac{9}4(-\frac12+\frac5{\tau_0})}.
\end{equation*}
Similarly as we mentioned before, from \eqref{NegativePowers}, we have $C\kappa^{-15}\mathbb{O}_{r_n}^{\frac{1}{2}}\le C\kappa^{-15}(\varepsilon_\ast)^{\frac{1}{2}}\le \frac14$ and since $r_n=\kappa^{\mathfrak{N}+n}R\le\kappa^\mathfrak{N}$ (recall $R<1$) by \eqref{conditionKappaAlone} we have $r_n^{\frac{9}4(-\frac12+\frac5{\tau_0})}\le \frac{\varepsilon}{2}.$ Then, we have
\begin{eqnarray*}
	\mathbb{O}_{r_{n+1}}(t_0,x_0)&\le&\frac14\mathbb{O}_{r_n}+\frac14\mathbb{O}_{r_n}+\frac{\varepsilon_\ast}{2}\le \frac{\varepsilon_\ast}{4}+\frac{\varepsilon_\ast}{4}+\frac{\varepsilon_\ast}{2}=\varepsilon_\ast. 
\end{eqnarray*} 
Thus,  we have proved that for all $n\in \mathbb{N}$, we have
\begin{equation*}
\mathbb{O}_{r_n}(t_0,x_0)\le \varepsilon_\ast,
\end{equation*}
which is the wished control \eqref{IterativeO}, but centered in the point $(t_0,x_0)$. In order to treat the general case $(t,x)\in Q_{r_n}(t_0,x_0)$, notice that since $Q_{r_n}(t,x)\subset Q_{2r_n}(t_0,x_0)$, we have 
\begin{equation*}
\mathbb{O}_{r_n}(t,x)\le 2^{3-\frac{15}{\tau_0}} \mathbb{O}_{2r_n}(t_0,x_0)<C.
\end{equation*}
Then, we have proved that for $0<\kappa\ll1$, there exists $0<\mathfrak{r}<R$ such that for all $n\in \mathbb{N}$, we have $\mathbb{O}_{\kappa^n \mathfrak{r}}(t,x)<C,$ and as it was mentioned before in \eqref{IterativeO1}, this implies  that for all $0<r\le\mathfrak{r}$ and for all $(t,x)\in Q_{r}(t_0,x_0)$, we have
\begin{equation}\label{ControlOr}
\Lambda_r(t,x)+\kappa^6\mathbb{P}_r(t,x)=\mathbb{O}_r(t,x)\le C.
\end{equation}
Hence, by  \eqref{Def_Invariants}, \eqref{Def_QuantiteIteration} and the previous estimate we have  for all $0<r\le\mathfrak{r}$ and for all $(t,x)\in Q_{r}(t_0,x_0)$,   $$\frac{1}{r^{5(1-\frac{3}{\tau_0})}}\int_{{\bf Q}_r(t,x)}\mathds{1}_{Q_{\mathfrak{r}}(t_0,x_0)}|\vu|^3 dy ds =\Lambda_r(t,x)\le C$$ which means that
 $$\mathds{1}_{Q_{\mathfrak{r}}(t_0,x_0)}\vu\in \M_{t,x}^{3,\tau_0}(\mathbb{R}\times \R).$$  Moreover, again by \eqref{Def_Invariants} and \eqref{Def_QuantiteIteration} we have 
$\displaystyle{\frac{1}{r^{5(1-\frac{5}{\tau_0})}} \int _{{\bf Q}_r(t,x)}\mathds{1}_{Q_{\mathfrak{r}}(t_0,x_0)}|p (s,y)|^{\frac{3}{2}}dyds=\mathcal{P}_r(t,x)}$, and by \eqref{ControlOr}, (recall that $\kappa$ is a fixed parameter) we conclude  
\begin{equation*}
\mathds{1}_{Q_{\mathfrak{r}}(t_0,x_0)}p\in \M_{t,x}^{\frac{3}{2},\frac{\tau_0}{2}}(\mathbb{R}\times \R).
\end{equation*}
and thus the proof of Proposition \ref{Pro_partialRegualrity} is finished.
\hfill$\blacksquare $\\
%%%%%%%%%%%%%%%%%%%%%%%%%%%%%%%%%%%%%%%%%%%%%%%%%%%
\begin{Corollary}\label{corolarioMorrey}
Under the hypotheses of Proposition \ref{Pro_partialRegualrity}, we have the following local control:
\begin{equation*}
\mathds{1}_{Q_{\frac{\mathfrak{r}}{2}}(t_0,x_0)}\grad\otimes \vu\in \mathcal{M}^{2,\tau_1}_{t,x}(\mathbb{R}\times \R)\quad \mbox{with}\quad\frac{1}{\tau_1}=\frac{1}{\tau_0}+\frac{1}{5}.
\end{equation*}
\end{Corollary}
%%%%%%%%%%%%%%%%%%%%%%%%%%%%%%%%%%%%%%%%%%%%%%%%%%%
\noindent\textbf{Proof.}  Let $0<\mathfrak{r}<R$ be the radius given in  Proposition \ref{Pro_partialRegualrity}. By using the definition of Morrey spaces given in \eqref{DefMorreyparabolico}, we have to show that  for all $\displaystyle{0<r\le \frac{\mathfrak{r}}2}$ and for all $(t,x)\in Q_r(t_0,x_0)$ we have
\begin{equation*}
\int_{{\bf Q}_r(t,x)}\mathds{1}_{Q_{\frac{\mathfrak{r}}{2}}(t_0,x_0)}|\grad \otimes\vu|^2dyds\leq C r^{5(1-\frac{2}{\tau_1})}, \quad \text{with} \quad \frac{1}{\tau_1}=\frac{1}{\tau_0}+\frac15.
\end{equation*}
For this, notice that  by the definition of the quantity $\lambda_r$ given in \eqref{Def_Invariants} and Lemma \ref{lema_FirstEstimate}, it follows for any $\displaystyle{0<r\le\frac{\mathfrak{r}}{2}}$ and any $(t,x)\in Q_r(t_0,x_0)$ that we have
\begin{eqnarray}
\frac{1}{r}\int_{{\bf Q}_r(t,x)}\mathds{1}_{Q_{\frac{\mathfrak{r}}{2}}(t_0,x_0)}|\grad \otimes\vu|^2dyds&\le& \mathcal{A}_{r}(t,x)+\alpha_r(t,x)\notag\\
&\leq& C\Big( \lambda_{2r}^{\frac{2}{3}}(t,x)+\mathcal{P}_{2r}(t,x)+\lambda_{2r}(t,x)+(2r)^{\frac{1}{2}}\lambda_{2r}^{\frac{1}{3}}(t,x)\Big).\qquad \label{gradientEstimate}
\end{eqnarray}
Let us study in more detail the terms  $\lambda_{2r}$ and $\mathcal{P}_{2r}$ above.  For the first one, since $\displaystyle{\lambda_{r}=\frac{1}{r^2}\int_{Q_{r}}|\vu|^3dyds}$ (see \eqref{Def_Invariants}), 
we have
\begin{align*}
\lambda_{2r}(t,x)=\frac{1}{(2r)^2}\int_{Q_{2r}(t,x)}|\vu|^3dyds=\frac{1}{(2r)^2}\frac{(2r)^{5(1-\frac{3}{\tau_0})}}{(2r)^{5(1-\frac{3}{\tau_0})}}\int_{Q_{2r}(t,x)}|\vu|^3dyds\\
=\frac{(2r)^{(3-\frac{15}{\tau_0})}}{(2r)^{5(1-\frac{3}{\tau_0})}}\int_{Q_{2r}(t,x)}|\vu|^3dyds.
\end{align*}
Since $2r\le \mathfrak{r}$ and $\mathds{1}_{Q_{\mathfrak{r}}(t_0,x_0)}\vu\in \M_{t,x}^{3,\tau_0}(\mathbb{R}\times \R)$, we obtain 
\begin{align}\label{EstimateLambda}
\lambda_{2r}(t,x)\le (2r)^{(3-\frac{15}{\tau_0})}\|\mathds{1}_{Q_{\mathfrak{r}}}\vu\|^3_{\M_{t,x}^{3,\tau_0}}\le C (2r)^{(3-\frac{15}{\tau_0})}.
\end{align}
Let us get a similar estimate for the term $\mathcal{P}_{2r}$, indeed by \eqref{Def_Invariants} and since $ \mathds{1}_{Q_{\mathfrak{r}}}p\in \M_{t,x}^{\frac{3}{2},\frac{\tau_0}{2}}(\mathbb{R}\times \R),$ we have 
\begin{eqnarray}
\mathcal{P}_{2r}(t,x)&=&\frac{1}{2r^2}\int_{Q_{2r}(t,x)}|p|^\frac32dyds=\frac{1}{(2r)^2}\frac{(2r)^{5(1-\frac{3}{\tau_0})}}{(2r)^{5(1-\frac{3}{\tau_0})}}\int_{Q_{2r}}|p|^\frac32dyds\notag\\
&\le& (2r)^{(3-\frac{15}{\tau_0})}\|\mathds{1}_{Q_{\mathfrak{r}}}p\|^{\frac32}_{\M_{t,x}^{\frac32,\frac{\tau_0}2}}\le C(2r)^{(3-\frac{15}{\tau_0})}.\label{estimatePress}
\end{eqnarray}
Thus, by using the estimates \eqref{EstimateLambda} and \eqref{estimatePress} in \eqref{gradientEstimate}, one has
\begin{eqnarray*}
\frac{1}{r}\int_{{\bf Q}_r(t,x)}\mathds{1}_{Q_{\frac{\mathfrak{r}}{2}}(t_0,x_0)}|\grad \otimes\vu|^2dyds &\leq & C\Big((2r)^{\frac23(3-\frac{15}{\tau_0})}+(2r)^{(3-\frac{15}{\tau_0})}+(2r)^{\frac{1}{2}}(2r)^{\frac13(3-\frac{15}{\tau_0})}\Big)\\
&\le &C\Big( r^{(2-\frac{10}{\tau_0})}+r^{(3-\frac{15}{\tau_0})}+r^{(\frac32-\frac{5}{\tau_0})}\Big).
\end{eqnarray*}
Notice, that  since $\frac{2}{15}\le \frac{1}{\tau_0}<\frac15$ we have ${0<2-\frac{10}{\tau_0}\le  {3-\frac{15}{\tau_0}}}\le  {\frac32-\frac{5}{\tau_0}}$, and since $0<r<1$ we get
\begin{equation*}
\frac{1}{r}\int_{{\bf Q}_r(t,x)}\mathds{1}_{Q_{\frac{\mathfrak{r}}{2}}(t_0,x_0)}|\grad \otimes\vu|^2dyds\le C r^{(2-\frac{10}{\tau_0})}.
\end{equation*}
Now, using the fact that $\frac{1}{\tau_0}=\frac{1}{\tau_1}-\frac15$, it follows that  for any $0<r\le\frac{\mathfrak{r}}{2}$ and  $(t,x)\in Q_{r}(t_0,x_0)$, we have 
\begin{equation*}
\int_{{\bf Q}_r(t,x)}\mathds{1}_{Q_{\frac{\mathfrak{r}}{2}}(t_0,x_0)}|\grad \otimes\vu|^2dyds\leq C r^{(5-\frac{10}{\tau_1})}=C r^{5(1-\frac{2}{\tau_1})},
\end{equation*}
which implies that
 ${\mathds{1}_{Q_{\frac{\mathfrak{r}}2}(t_0,x_0)}\grad\otimes \vu\in \mathcal{M}^{2,\tau_1}_{t,x}(\mathbb{R}\times \R)}$ and this finishes the proof of Corollary \ref{corolarioMorrey}.\hfill$\blacksquare $\\
 
\noindent \textbf{Proof of Theorem \ref{Theo_partialRegularity}.}
 Let $(\vu,p, \vw)$ be a partial suitable solution in the sense of Definition \ref{Def_PartialSuitable} for the micropolar equations \eqref{MicropolarFluidsEquationsEqua1} and \eqref{MicropolarFluidsEquationsEqua2} in $Q_1(t_0,x_0)$.   Recall that we want to show that there exists some $0<r<R<1$ such that $\vu,\vw\in L_{t,x}^{\infty}(Q_r(t_0,x_0))$.\\
 
First, notice that since by hypothesis we have
\begin{equation*}
\frac{1}{R^2}\int_{t_0-R^2}^{t_0}\int _{B_{x_0,R}}| \vu|^3+|p|^\frac{3}{2}dx ds<\varepsilon,
\end{equation*}
 for some $\varepsilon\ll1$, we can apply  Proposition \ref{Pro_partialRegualrity} and Corollary \ref{corolarioMorrey} and therefore there exists $0<\mathfrak{r}<R$ such that
 for $5<\tau_0<\frac{15}2$ and $\frac{1}{\tau_0}=\frac{1}{\tau_1}-\frac15$ we have
\begin{equation}\label{FirtGain}
 \begin{aligned}
&{\mathds{1}_{Q_{\mathfrak{r}}(t_0,x_0)}\vu\in\M_{t,x}^{3,\tau_0}(\mathbb{R}\times \R)},\quad\mathds{1}_{Q_{\mathfrak{r}}(t_0,x_0)}p\in \M_{t,x}^{\frac{3}{2},\frac{\tau_0}{2}}(\mathbb{R}\times \R)\\
&\quad  \text{ and} \quad \mathds{1}_{Q_{\frac{\mathfrak{r}}2}(t_0,x_0)}\grad\otimes \vu\in \mathcal{M}^{2,\tau_1}_{t,x}(\mathbb{R}\times \R).
\end{aligned}
\end{equation}
 Note that the upper bound for $\tau_0$ comes from the fact that we have the term $\rot\vw$ in the equation \eqref{MicropolarFluidsEquationsEqua1}. \\
 
 Now, for simplicity sake, we assume $\tau_0=6$. Hence, by Proposition \ref{Theo_MorreyTransfer}, it follows that for some $0<r_1<\mathfrak{r}$, we have
 \begin{equation}\label{SecondGain}
 \mathds{1}_{Q_{r_1}(t_0,x_0)}\vu\in L^6_{t,x}(\mathbb{R}\times \R)\quad \text{and} \quad \mathds{1}_{Q_{r_1}(t_0,x_0)}\vw\in L^6_{t,x}(\mathbb{R}\times \R).
 \end{equation} 
It is worth noting that the integrability we have obtained lies within the framework of the Serrin criterion  $(\frac2p+\frac3q\le 1)$. However, instead of deducing directly the boundedness of the solution, we will apply the strategy given in \cite{ChLl23}.
Indeed, since we have \eqref{FirtGain} and \eqref{SecondGain} we can apply \cite[Proposition 2]{ChLl23}, and thus  we can improve the Morrey information of the velocity $\vu$ \textit{i.e.,} for any $r_2<r_1<R$, we have that 
\begin{equation*}
\mathds{1}_{Q_{\mathfrak{r}}(t_0,x_0)}\vu\in\M_{t,x}^{3,60}.
\end{equation*}
Moreover, following the same steps as in Section 6 of \cite{ChLl23}, we can deduce that $(\vu,\vw)$ are Holder continuous in time and space in $Q_r(t_0,x_0)$ for some $0<r<r_2<r_1<R$. Since  $Q_r(t_0,x_0)$ is a bounded set, the boundedness of $(\vu,\vw)$ follows immediately, thus completing the proof of Theorem \ref{Theo_partialRegularity}. \hfill $\blacksquare$\\

Finally we present the following characterization of partial singular points, 
which is just a consequence of the $\varepsilon $-regularity theory.
%%%%%%%%%%%%%%%%%%%%%%%%%%%%%%%%%%%%%%%%%%%%%%%%%%%
\begin{Proposition}\label{Propo_CaracSingular}
Let $(\vu,p,\vw)$ be a partial suitable solution on $Q_1$. Then, for any $(t_0,x_0)\in Q_1(t,x)$ we have
\begin{itemize}
\item either $(t_0,x_0)$ is partially singular and then for any $0<r<1,$
\begin{equation*}
\varepsilon\le\frac{1}{r^2}\int_{Q_r(t_0,x_0)}|\vu|^3+|p|^{\frac{3}{2}}dyds,
\end{equation*}
\item either $(t_0,x_0)$ is a partial regular point and then
\begin{equation*}
\lim_{r\longrightarrow 0}\frac{1}{r^2}\int_{Q_r(t_0,x_0)}|\vu|^3+|\vw|^3dyds=0.
\end{equation*}
\end{itemize}
\end{Proposition}
%%%%%%%%%%%%%%%%%%%%%%%%%%%%%%%%%%%%%%%%%%%%%%%%%%%
\noindent {\bf Proof.}
Let us prove the first point by contradiction. Hence, assume that $(t_0,x_0)$ is a partial singular point in the sense of Definition \ref{Def_PartialRegularPoints} such that there exists $0<r<R$ with
\begin{equation*}
\frac{1}{r^2}\int_{Q_r(t_0,x_0)}|\vu|^3+|p|^{\frac{3}{2}}dyds< \varepsilon.
\end{equation*}
Since $(\vu,p,\vw)$ is a partial suitable solution, we can use Theorem \ref{Theo_partialRegularity}, and therefore there exists $0<\rho<r$ such that $(\vu,\vw)$ is bounded on $Q_\rho(t_0,x_0)$, and hence $(t_0,x_0)$ has to be a partial regular point which is a contradiction. 
For the second point, since $(t_0,x_0$) is a partial regular point, there exists some $R>0$ such that $\vu,\vw \in L^\infty_{t,x}(Q_R(t_0,x_0))$. Hence, it is easy to see that for all $r<R$
\begin{equation*}
\frac{1}{r^2}\int_{Q_r(t_0,x_0)}|\vu|^3+|\vw|^3 dyds\le C(\|\vu\|_{L^\infty(Q_R)}^3+\|\vw\|_{L^\infty(Q_R)}^3)r^3.
\end{equation*}
The proof is completed by taking the limit when $r$ goes to zero.
\hfill $\blacksquare $


\begin{thebibliography}{2}
%%%%%%%%%%%%%%%%%%%%
\bibitem{Admas}
D. R. \textsc{Adams}, J. \textsc{Xiao}. \emph{Morrey spaces in harmonic analysis}. Ark. Mat. Vol. 50(2):201-230, (2012).
%%%%%%%%%%%%%%%%%%%%
\bibitem{BarPran19} 
T. \textsc{Barker}, C. \textsc{Prange}. \emph{Localized Smoothing for the Navier–Stokes Equations and Concentration of Critical Norms Near Singularities}. Arch. Rational. Mech. Anal. 236:1487–1541, (2020). 
%%%%%%%%%%%%%%%%%%%%%
\bibitem{BarPran20} 
T. \textsc{Barker}, C. \textsc{Prange}. \emph{Quantitative regularity for the Navier-Stokes equations via spatial concentration}. Commun. Math. Phys. 385:717–792, (2021). 
%%%%%%%%%%%%%%%%%%%%%
\bibitem{BegBahr08}
O. A. \textsc{Beg}, R. \textsc{Bhargava}, S. \textsc{Rawat}, K. \textsc{Halim} \& H. S. \textsc{Takhar}. \emph{Computational modeling of biomagnetic micropolar blood flow and heat transfer in a two-dimensional non-Darcian porous medium}. Meccanica. 43(4):391–410, (2008).
%%%%%%%%%%%%%%%%%%%%%
\bibitem{SerrinMicropolar}
D. \textsc{Bo-Qing}, C. \textsc{Zhin-Min}. \emph{Regularity criteria of weak solutions to the three-dimensional micropolar flows}. J. Math. Phys. 50:103525, (2009).
%%%%%%%%%%%%%%%%%%%%%
\bibitem{BradTsai}
Z. \textsc{Bradshaw}, T. P. \textsc{Tsai}. \emph{Global existence, regularity, and uniqueness of infinite energy solutions to the Navier-Stokes equations}. Comm. Partial Differential Equations, (2020)
%%%%%%%%%%%%%%%%%%%%
\bibitem{Brezis}
H. \textsc{Brezis}, P. \textsc{Mironescu}. \emph{Gagliardo-Nirenberg inequalities and non-inequalities: the full story}. Annales de l’Institut Henri Poincaré (C) Non Linear Analysis. Elsevier. 35(5):1355-1376, (2018).
%%%%%%%%%%%%%%%%%%%%
\bibitem{CKN} 
L. \textsc{Caffarelli}, R. \textsc{Kohn} \& L. \textsc{Nirenberg}. \emph{Partial regularity of suitable weak solutions of the Navier-Stokes equations}. Comm. Pure Appl. Math. 35(6):771-831, (1982).
%%%%%%%%%%%%%%%%%%%%%
\bibitem{Caz98}
T. \textsc{Cazenave}, A. \textsc{Haraux}. \emph{An introduction to semilinear evolution equations}. Oxford Lecture Series in Mathematics and its Applications 13. Oxford University Press. Oxford, (1998). 
%%%%%%%%%%%%%%%%%%%%%%%%%%%%
\bibitem{ChHe21}
D. \textsc{Chamorro}, J. \textsc{He}. \emph{On the partial regularity theory for the {MHD} equations}. J. Math. Anal. Appl. 494:124449, (2021). 
%%%%%%%%%%%%%%%%%%%%%
\bibitem{ChLl21}
D. \textsc{Chamorro}, D. \textsc{Llerena}. \emph{Interior e-regularity theory for the solutions of the magneto-micropolar equations with a perturbation term}. J. Elliptic Parabol. Equ. 8:555–616, (2022). 
%%%%%%%%%%%%%%%%%%%%%
\bibitem{ChLl22}
D. \textsc{Chamorro}, D. \textsc{Llerena}. \emph{A crypto-regularity result for the micropolar fluid equations}. J. Math. Anal. Appl. 520(2):126922, (2023).
%%%%%%%%%%%%%%%%%%%%%
\bibitem{ChLl23}
D. \textsc{Chamorro}, D. \textsc{Llerena}. \emph{Partial suitable solutions for the micropolar equations and regularity properties}. Prepint, arXiv:2302.02675, (2023).
%%%%%%%%%%%%%%%%%%%%
\bibitem{Eri66}
A. C. \textsc{Eringen}. \emph{Theory of micropolar fluids}. J. Math. Mech. 16:1–18, (1966).
%%%%%%%%%%%%%%%%%%%%
\bibitem{Eus03}
L. \textsc{Escauriaza}, G. \textsc{Seregin} \& W. \textsc{Sverak}. \emph{L$^{\infty,3}$-solutions of the Navier-Stokes equations and backward uniqueness}. Uspekhi Mat. Nauk. 58(350):3-44, (2003).
%%%%%%%%%%%%%%%%%%%%
\bibitem{Evans}
L. \textsc{Evans}. \emph{Partial differential equations}.  American Mathematical Society, (2010).
%%%%%%%%%%%%%%%%%%%%
\bibitem{Gal97}
G. P. \textsc{Galdi}, S. \textsc{Rionero}. \emph{A note on the existence and uniqueness of solutions of the micropolar fluid equations}. Internat. J. Eng. Sci. 15:105–108, (1977).
%%%%%%%%%%%%%%%%%%%%
\bibitem{Gay13}
F. \textsc{Gay-Balmaz}, T. S. \textsc{Ratiu}, \& C. \textsc{Tronci}. \emph{Equivalent theories of liquid crystal dynamics}. Arch. Ration. Mech. Anal. 210(3):773–811, (2013).
%%%%%%%%%%%%%%%%%%%%%
\bibitem{KanMiuTsai20} K. \textsc{Kang}, H. \textsc{Miura} \& T. \textsc{Tsai}. \emph{Short time regularity of Navier-Stokes flows with locally $L^3$ initial data and applications}. International Mathematics Research Notices 11:8763–8805, (2020).
%%%%%%%%%%%%%%%%%%%%%
\bibitem{KanMiuTsai19} K. \textsc{Kang}, H. \textsc{Miura} \& T. \textsc{Tsai}. \emph{An $\varepsilon $-regularity criterion and estimates of the regular set for Navier-Stokes flows in terms of initial data}. Pure Appl. Anal. 3:567-594, (2021).
%%%%%%%%%%%%%%%%%%%%%
\bibitem{Kukavica}
I. \textsc{Kukavica}. \emph{On partial regularity for the Navier-Stokes equations}. Discrete Contin. Dyn. Syst. 21(3):717-728, (2008).
%%%%%%%%%%%%%%%%%%%%
\bibitem{LoMelo}
J. \textsc{Lorenz}, W. G. \textsc{Melo} \& S. C. P. \textsc{de Souza}. \emph{Regularity criteria for weak solutions of the Magneto-micropolar equations}. Electronic Research Archive. 29(1):1625–1639, (2021).
%%%%%%%%%%%%%%%%%%%
\bibitem{PGLR0} P.-G. \textsc{Lemarié-Rieusset}. \emph{Recent Developments in the Navier-Stokes problem}. Chapman \& Hall/CRC, (2002).
%%%%%%%%%%%%%%%%%%%%
\bibitem{PGLR1} P.-G. \textsc{Lemarié-Rieusset}. \emph{The Navier-Stokes problem in the 21st century}. Chapman \& Hall/CRC, (2016).
%%%%%%%%%%%%%%%%%%%%%
\bibitem{Mayea}
Y. \textsc{Maekawa}, H. \textsc{Miura} \& C. \textsc{Prange}. \emph{Local Energy Weak Solutions for the Navier-Stokes Equations in the Half-Space}. Commun. Math. Phys. 367:517–580, (2019).
%%%%%%%%%%%%%%%%%%%%
\bibitem{Mit02}
N. \textsc{Mitarai}, H. \textsc{Hayakawa} \& H. \textsc{Nakanishi}. \emph{Collisional Granular Flow as a Micropolar Fluid}. Phys. Rev. Lett. 88(17):174301, (2002).
%%%%%%%%%%%%%%%%%%%%
\bibitem{OLeary}
M. \textsc{O'Leary}. \emph{Conditions for the local boundedness of solutions of the Navier-Stokes system in three dimensions}. Comm. Partial Differential Equations. 28:617-636, (2003).
%%%%%%%%%%%%%%%%%%%%%
\bibitem{Robinson} 
J. \textsc{Robinson}, J. \textsc{Rodrigo} \& W. \textsc{Sadowski}. \emph{The three dimensional Navier-Stokes equations}. Cambridge Studies in advanced Mathematics, (2016).
%%%%%%%%%%%%%%%%%%%%%
\bibitem{She77}
V. \textsc{Scheffer}. \emph{Hausdorff measure and the Navier-Stokes equation}. Comm. Math. Phys. 55:97-112, (1977).
%%%%%%%%%%%%%%%%%%%%%%%
\bibitem{Seregin14}
G. \textsc{Seregin}. \emph{Lecture Notes on Regularity Theory for the Navier-Stokes Equations}. World Scientific Publishing, (2014).
%%%%%%%%%%%%%%%%%%%%
\bibitem{Serr62}
J. \textsc{Serrin}. \emph{On the interior regularity of weak solutions of the Navier-Stokes equations}. Arch. Rat. Mech. Anal. 9:187-195, (1962).
%%%%%%%%%%%%%%%%%%%%
\bibitem{Serr63}
J. \textsc{Serrin}. \emph{The initial value problem for the Navier-Stokes equations}. In Nonlinear Problems, Rudolph E. Langer ed., Univ. Wisconsin Press. Madison, (1963).
%%%%%%%%%%%%%%%%%%%%
\bibitem{Struwe88}
M. \textsc{Struwe}. \emph{On partial regularity results for the Navier-Stokes equations.} American Mathematical Soc., (2018).
%%%%%%%%%%%%%%%%%%%
\bibitem{Takahashi90}
S. \textsc{Takahashi}. \emph{On partial regularity results for the Navier-Stokes equations}. Comm. Pure Appl. Math. 41:437–458, (1988).
%%%%%%%%%%%%%%%%%%%
\bibitem{Tsai19}
T. \textsc{Tsai}. \emph{Lectures on Navier-Stokes Equations}. Graduate Studies in Mathematics. American Mathematical Society, (2018).
%%%%%%%%%%%%%%%%%%%%%
\bibitem{Yam05}
N. \textsc{Yamaguchi}. \emph{On interior regularity criteria for weak solutions of the Navier-Stokes equations}. Manuscripta Math. 69:237–254, (1990).
%%%%%%%%%%%%%%%%%%%%%
\bibitem{Yua08}
B. \textsc{Yuan}. \emph{Regularity of weak solutions to magneto-micropolar fluid equations}. Math. Sci. Ser. 30:1469–1480, (2010).
\end{thebibliography}
\end{document}